# Syntactic approaches to opetopes

Pierre-Louis Curien

Cédric Ho Thanh[1]

Samuel Mimram

Research Institute for Foundations of Computer Science (IRIF), Paris Diderot University, Paris, France
*Email address*: curien@irif.fr
*URL*: www.irif.fr/~curien

Research Institute for Foundations of Computer Science (IRIF), Paris Diderot University, Paris, France
*Email address*: cedric.hothanh@irif.fr
*URL*: hothanh.fr/cedric

Computer Science Laboratory of École Polytechnique (LIX), Palaiseau, France
*Email address*: samuel.mimram@lix.polytechnique.fr
*URL*: www.lix.polytechnique.fr/Labo/Samuel.Mimram

[1]This author has received funding from the European Union's Horizon 2020 research and innovation program under the Marie Sklodowska-Curie grant agreement №665850.

Abstract. Opetopes are algebraic descriptions of shapes corresponding to compositions in higher dimensions. As such, they offer an approach to higher-dimensional algebraic structures, and in particular, to the definition of weak $\omega$-categories, which was the original motivation for their introduction by Baez and Dolan. They are classically defined inductively (as free operads in Leinster's approach, or as zoom complexes in the formalism of Kock et al.), using abstract constructions making them difficult to manipulate with a computer.

In this paper, we present two purely syntactic descriptions of opetopes as sequent calculi, the first using variables to implement the compositional nature of opetopes, the second using a calculus of higher addresses. We prove that well-typed sequents in both systems are in bijection with opetopes as defined in the more traditional approaches. Additionally, we propose three variants to describe opetopic sets. We expect that the resulting structures can serve as natural foundations for mechanized tools based on opetopes.

# Contents

CHAPTER 1

# Introduction

## 1.1. Opetopes

Opetopes were originally introduced by Baez and Dolan in order to formulate a definition of weak $\omega$-categories [**Baez and Dolan, 1998**]. Their name reflects the fact that they encode the possible shapes for higher-dimensional operations: they are *ope*ration poly*topes*. Over the recent years, they have been the subject of many investigations in order to provide good working definitions of opetopes allowing to explore their combinatorics [**Cheng, 2003, Hermida et al., 2000, Leinster, 2004**]. One of the most commonly used nowadays is the formulation based on polynomial functors and the corresponding graphical representation using "zoom complexes" [**Kock et al., 2010**].

In order to grasp quickly the nature of opetopes, consider a sequence of four composable arrows

$$a \xrightarrow{f} b \xrightarrow{g} c \xrightarrow{h} d \xrightarrow{i} e$$

There are various ways we can compose them. For instance, we can compose $f$ with $g$, as well as $h$ with $i$, and then $g \circ f$ with $i \circ h$. Or we can compose $f$, $g$ and $h$ together all at once, and then the result with $i$. These two schemes for composing can respectively be pictured

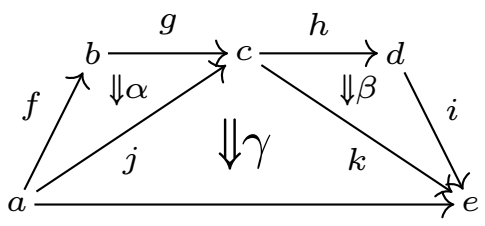

and

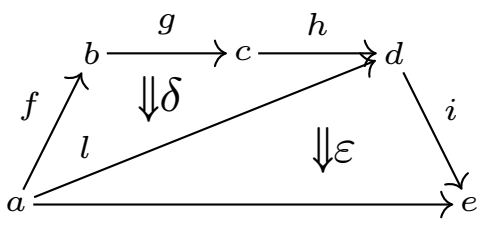

(1.1.1)

From there, the general idea of getting "higher-dimensional" is that we should take these compositions as "2-operations", which can be composed in various ways. For instance, in the first case, we can compose $\alpha$ with $\gamma$, and then $\beta$ with the result, or all three at once, and so on. The opetopes describe all the ways in which these compositions can be meaningfully specified, in arbitrary dimension.

We can expect (and it is indeed the case) that the combinatorics of these objects is not easy to describe. In particular, a representation which is adapted to computer manipulations and proofs is desirable: we can for instance mention the *Opetopic* proof assistant for higher categories [**Finster, 2016**], which is based on opetopes. Recently, Mimram and Finster have used sequent calculus as a convenient tool to describe globular weak $\omega$-categories [**Finster and Mimram, 2017**]. Our goal is to achieve a similar type-theoretic presentation for opetopic weak $\omega$-categories [**Baez and Dolan, 1998, Hermida et al., 2002**], and we begin here by defining a representation of opetopes and opetopic sets of type-theoretic flavor.

Let us now informally define opetopes. At the basis of the architecture, there is a unique 0-opetope (i.e. opetope of dimension 0), drawn as a point. An $(n+1)$-opetope (i.e. an opetope of dimension $n+1$) is made out of $n$-opetopes, and has a source and a target. In small dimensions, opetopes can be described using drawings of the kind above. For instance, the following are, from left to right, the unique 0-opetope, the unique 1-opetope, drawn as an arrow, and three 2-opetopes, respectively,

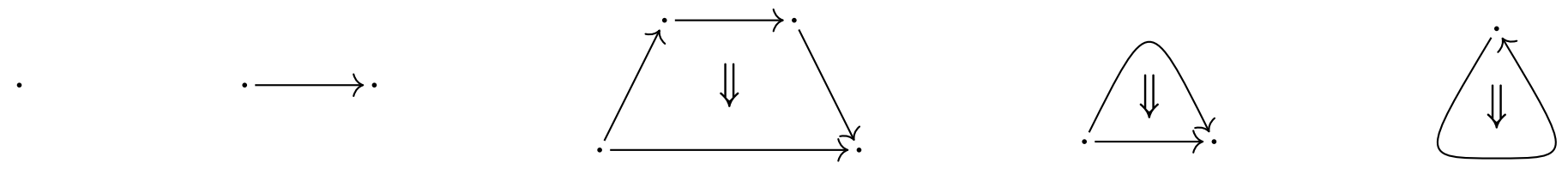

while the following drawing represents a 3-opetope:

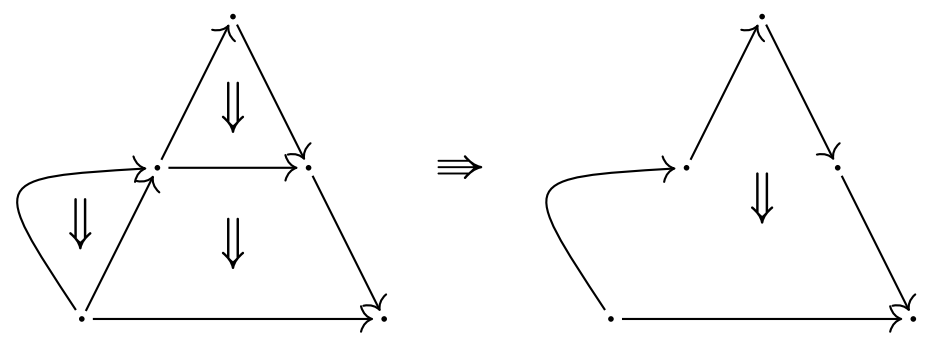

In these drawings, the target of the 1-opetope is a point, the target of the 2-opetopes is the arrow at the bottom, and the target of the 3-opetope is the 2-opetope on the right, while the source of the 1-opetope is again a point, the source of the 2-opetopes is a diagram made of arrows, and the source of the 3-opetope is the diagram made of 2-opetopes on the left of the central arrow.

At this stage, we can already make the following observations:

(1) each drawing contains only one top-dimensional cell, whose dimension determines the dimension of the opetope,
(2) if an opetope has a non empty source, then its source and target have the same source (which might be empty) and the same target,
(3) if an opetope has an empty source, then its target has the same source and target,
(4) an $(n+1)$-cell can have multiple $n$-cells in its source (including none) but always has exactly one $n$-cell in its target.

By the above remarks, an $n$-opetope has a unique top-dimensional cell $\alpha$, whose target consists of one cell $\beta$ whose source and target are the same as the ones of the source of $\alpha$. It follows recursively that the opetope is entirely determined by the source of $\alpha$, which we call the associated *pasting scheme*, which is of dimension $n-1$. The opetope associated to a pasting scheme of dimension $n$ can be obtained by adding a single cell $\beta$ of dimension $n$ parallel to the pasting scheme, and an $(n+1)$-cell $\alpha$ from the pasting scheme to $\beta$.

This article aims to provide syntactic tools making opetopes easier to manipulate than their classical definitions. The *unnamed* (or *anonymous*) approaches of chapter 4 on page 47 are purely relying on a calculus of higher addresses, and on a syntactic notion of *preopetope*. Opetopes are then well-formed preopetopes according to the derivation system $\text{Opt}^?$ presented in figure 4.1.1 on page 50, while opetopic sets are derivable contexts in system $\text{OptSet}^?$. The *named* approaches of chapter 3 on page 15, while slightly more complicated, leverage the idea of cell naming to produce user friendlier tools and results. As such it comes in two variants: $\text{Opt}^!$ for describing opetopes, and $\text{OptSet}^!$ for opetopic sets, introduced in figures 3.1.1 and 3.5.1 on page 17 and on page 30 respectively. The Python implementation of [**Ho Thanh, 2018b**] is also discussed.

## 1.2. Generating opetopes

Our two opetope derivation systems $\text{Opt}^!$ (figure 3.1.1 on page 17) and $\text{Opt}^?$ (figure 4.1.1 on page 50) are based on the observation that opetopes are precisely all the shapes one can generate with the following operations.

(1) *Introduction of a point.* There is a unique 0-opetope (the *point*).

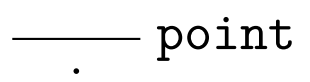

(2) *Shift to the next dimension.* Given an $n$-opetope $\omega$, we can form the $(n+1)$-globe whose source and target are $\omega$, as illustrated below. It can geometrically be thought of as the "extrusion" of $\omega$.

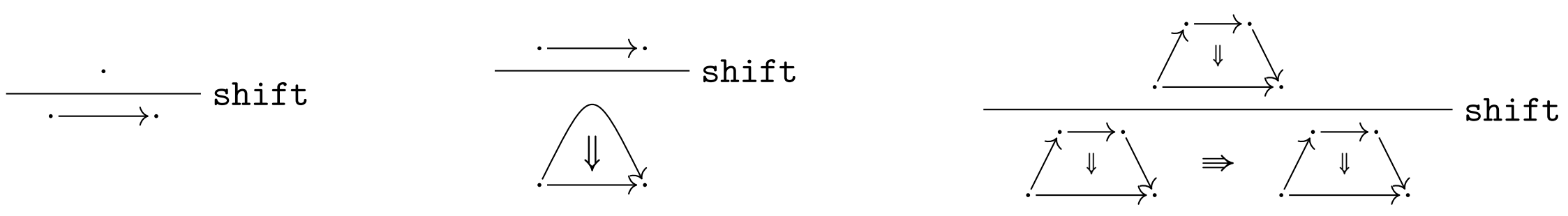

(3) *Introduction of degeneracies.* Given an $n$-opetope $\omega$, we can build an $(n+2)$-opetope with empty source, whose target is the globe at $\omega$, as illustrated below for $n = 0$ and $n = 1$:

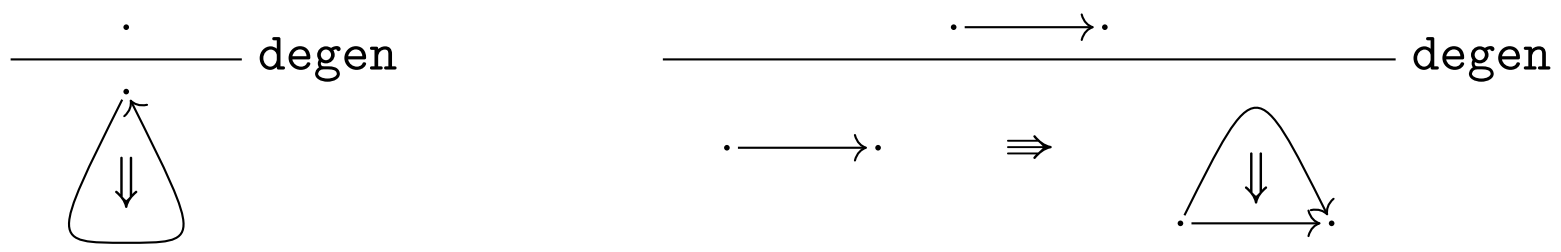

(4) *Grafting.* Given an $(n+1)$-opetope $\alpha$ and an $(n+1)$-pasting scheme $\beta$ such that the source of $\beta$ contains an $n$-cell of the same shape as the target of $\alpha$, we can *graft* $\alpha$ to $\beta$:

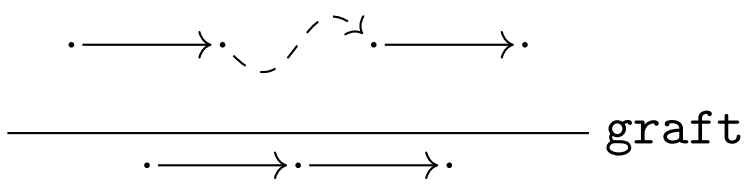

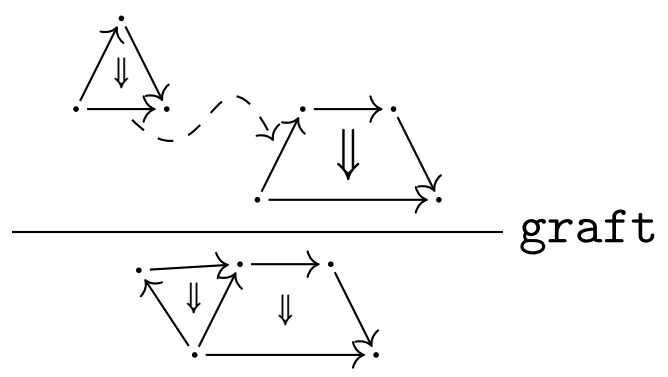

Ill-formed graftings may occur with $n$-pasting diagrams, for $n \geq 3$, and the side condition is necessary to rule them out. Here is an example the `graft` rule will not allow: we deal with a 3-pasting diagram on the right of the dashed arrow (that comprises a unique 3-opetope), and the dashed arrow indicated that we attempt to graft the 3-opetope on the left (whose target shape is a trapezoid) onto the triangle shaped cell on the right

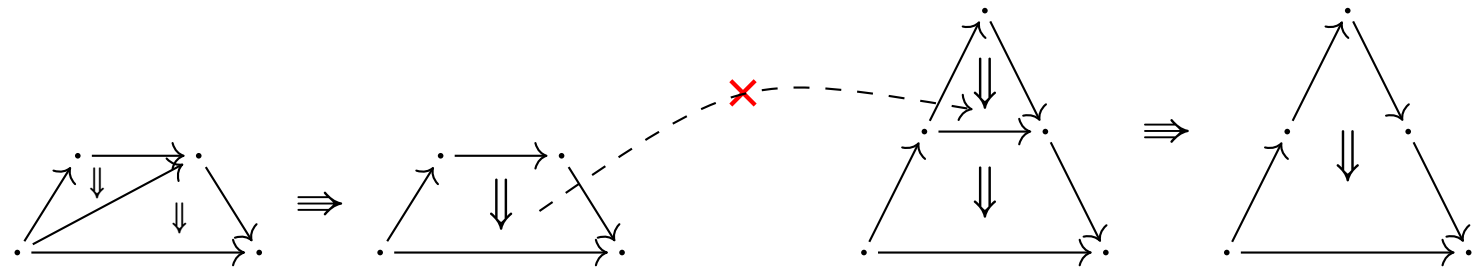

However, grafting onto the lower trapezoid of the right opetope is possible:

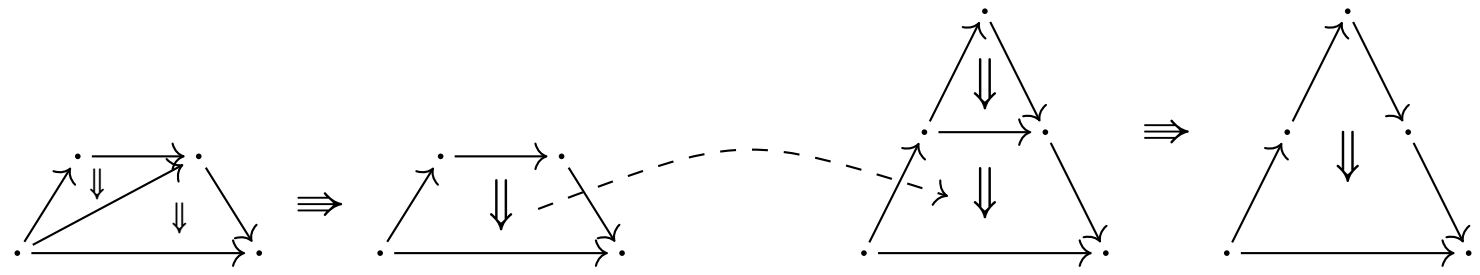

As previously mentioned, an opetope is completely determined by its source pasting scheme, i.e. "arrangement" of source faces (the dichotomy between opetopes and pasting schemes is more thoroughly discussed in section 2.2.2 on page 11). We can reformulate rules `shift` and `graft` with this point of view to respectively obtain:

(1) *Filling of pasting schemes.* Given an $n$-pasting scheme, we may "fill" it by adding a target $n$-cell, and a top dimensional $(n+1)$-cell. We illustrate an instance of this rule on the left, and invite the reader to compare it with the instance of `shift` on the right:

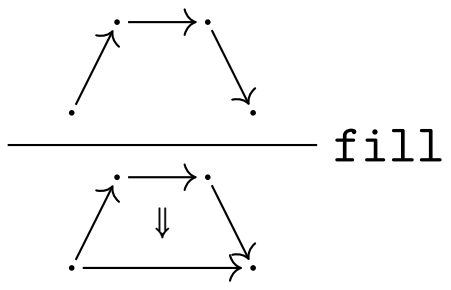

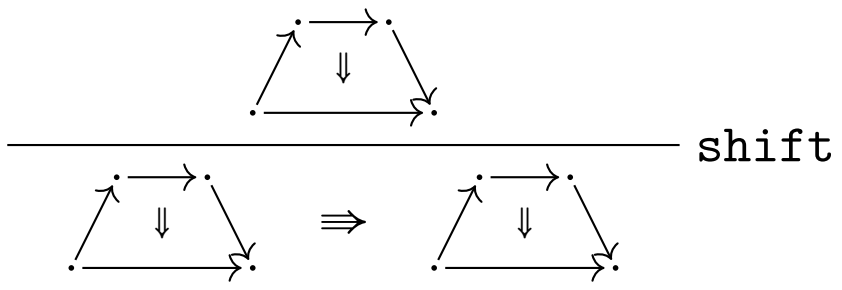

(2) *Substitution*, which consists in replacing a cell in a pasting scheme by another "parallel" pasting scheme. As before, we illustrate an instance if this rule on the left, and invite the reader to compare it with the instance of `graft` on the right:

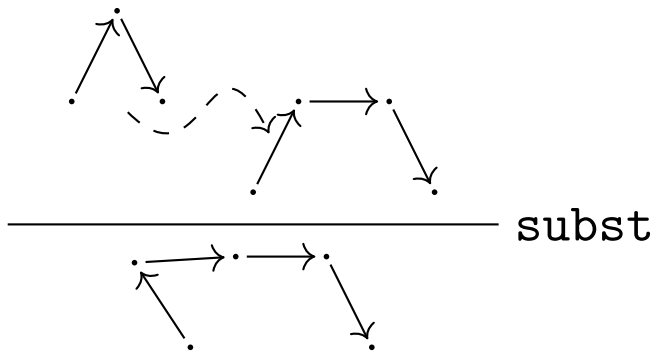

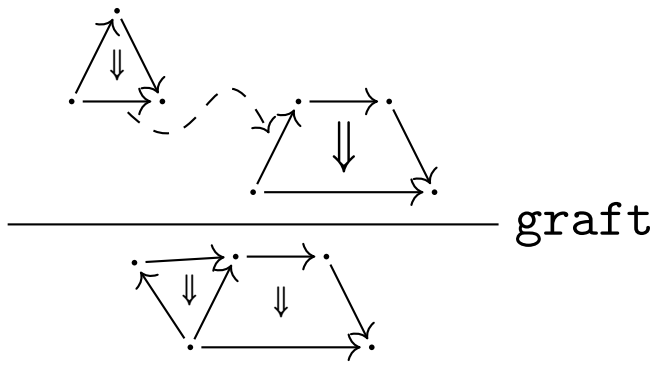

In the case of opetopes, substitution can be illustrated as follows:

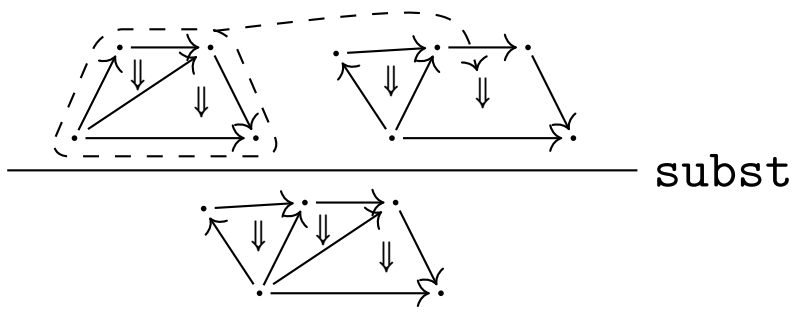

Rules `point` and `degen` can also be reformulated in term of pasting schemes, but representing them graphically is not as insightful.

## 1.3. Plan

In chapter 2 on the facing page, we expound some prerequisites and general ideas motivating our syntactic approaches. We then present the "unnamed" and "named" systems for opetopes and opetopic sets in chapters 3 and 4 on page 15 and on page 47 respectively. Those two chapters can be read independently, but have a similar structure:

(1) firstly, we present the syntactic constructs and inference rules of systems $\mathrm{O{\scriptstyle PT}}^!$ (named) and $\mathrm{O{\scriptstyle PT}}^?$ (unnamed) for opetopes, in sections 3.1 and 4.1 respectively;
(2) we then prove the equivalence with opetopes as defined in [**Kock et al., 2010**] (which we shall refer to as *polynomial opetopes*) in sections 3.2 and 4.2;
(3) we showcase the $\mathrm{O{\scriptstyle PT}}^!$ and $\mathrm{O{\scriptstyle PT}}^?$ systems in examples in sections 3.3 and 4.3;
(4) we discuss the Python implementation of $\mathrm{O{\scriptstyle PT}}^!$ and $\mathrm{O{\scriptstyle PT}}^?$ of [**Ho Thanh, 2018b**] in sections 3.4 and 4.5;
(5) in the second part of those chapters, we present the syntactic constructs and inference rules of systems $\mathrm{O{\scriptstyle PT}S{\scriptstyle ET}}^!$ (named) and $\mathrm{O{\scriptstyle PT}S{\scriptstyle ET}}^?$ (unnamed) for opetopic sets, in sections 3.5 and 4.6 respectively;
(6) we then prove the equivalence with opetopic sets of [**Ho Thanh, 2018a**] in sections 3.6 and 4.7;
(7) we showcase the $\mathrm{O{\scriptstyle PT}S{\scriptstyle ET}}^!$ and $\mathrm{O{\scriptstyle PT}S{\scriptstyle ET}}^?$ systems in examples in sections 3.7 and 4.8;
(8) we discuss the Python implementation of $\mathrm{O{\scriptstyle PT}S{\scriptstyle ET}}^!$ and $\mathrm{O{\scriptstyle PT}S{\scriptstyle ET}}^?$ sections 3.8 and 4.9;
(9) for the named approach of chapter 3, we present the additional "mixed" system $\mathrm{O{\scriptstyle PT}S{\scriptstyle ET}}^!_{\mathrm{M}}$ in section 3.9 on page 39, with the same plan as previously presented.

The equivalence proofs with polynomial opetopes of sections 3.2 and 4.2 on page 20 and on page 52 rely on a precise definition of the Baez–Dolan $(-)^+$ construction [**Kock et al., 2010, Baez and Dolan, 1998**], which is presented in appendix A on page 71.

## 1.4. Related works

Our syntactic opetopes are shown in sections 3.2 and 4.2 on page 20 and on page 52 to be equivalent to the polynomial opetopes (or "zoom complexes") of Kock et al. [**Kock et al., 2010**], which are themselves equivalent to Leinster's opetopes [**Leinster, 2004**]. It is known that the latter are incompatible with Cheng's opetopes [**Cheng, 2003**], which should be thought of as being symmetric. There is a closely related notion of *multitope* [**Hermida et al., 2002, Harnik et al., 2008**] which is defined in terms of multicategories (whence the *multi*) instead of operads (*ope*); the two notions can be shown to be equivalent [**Ho Thanh, 2018a**]. A syntax for multitopes was proposed in [**Hermida et al., 2002**], where however not all the desired computations have been given algorithmic formulations.

The Opetopic proof assistant [**Finster, 2016**] for weak higher categories relies on the notion of higher-dimensional tree. In that system, the notion of opetope is built-in, so that we have to trust the implementation. In contrast, the present approach allows us to reason about the construction of opetopes. We moreover believe that the ability to reason by induction on the proof trees, together with the very explicit nature of our syntaxes, will allow for optimizations in the automated manipulations of opetopes.

Another proof assistant for weak higher categories, called CaTT [**Finster and Mimram, 2017**], starts from the same idea of generating well-formed pasting schemes through inference rules. However, it is based on globular shapes instead of opetopic ones, making a comparison with the present work difficult: since their introduction, people have unsuccessfully tried to compare the resulting respective categorical formalisms; we hope that their formulation in a common logical language might be of help in this task.

We should also mention here the *Globular* proof assistant [**Bar et al., 2016**], also based on globular shapes, which is quite popular, notably thanks to its nice graphical interface.

CHAPTER 2

# Preliminaries

## 2.1. Polynomial functors and trees

We recall basic knowledge on the theory of polynomial functors from e.g. [**Kock, 2011, Kock et al., 2010**]. A *polynomial endofunctor* $F$ *over* $I$ is a $\mathcal{S}\mathrm{et}$-diagram of the form:

$$I \xleftarrow{s} E \xrightarrow{p} B \xrightarrow{t} I. \tag{2.1.1}$$

Elements of $B$ are called *nodes* or *operations*, elements of the fiber $E(b) := p^{-1}(b)$ are the *inputs* or *arities* of $b$, and elements of $I$ are *colors*. For $b \in B$, if $e \in E(b)$, let $s_e(b) := s(e)$. We will sometimes refer to $I$, $B$ and $E$ as $I_F$, $B_F$, and $E_F$, respectively. A morphism $f : F \longrightarrow F'$ of polynomial functors is a diagram of the form

$$\begin{array}{ccccccc} I_F & \longleftarrow & E_F & \longrightarrow & B_F & \longrightarrow & I_F \\ {\scriptstyle f}\big\downarrow & & {\scriptstyle f}\big\downarrow & \lrcorner & \big\downarrow{\scriptstyle f} & & \big\downarrow{\scriptstyle f} \\ I_{F'} & \longleftarrow & E_{F'} & \longrightarrow & B_{F'} & \longrightarrow & I_{F'} \end{array}$$

where all three squares commute, and the middle one is cartesian. This latter condition means that for $b \in B_F$, the map $f : E_F \longrightarrow E_{F'}$ exhibits a bijection between the set $E_F(b)$ of inputs of $b$, and $E_{F'}(f(b))$ of inputs of $f(b)$. We note $\mathcal{P}\mathrm{oly}\mathcal{E}\mathrm{nd}$ the category of polynomial endofunctors and morphisms, and $\mathcal{P}\mathrm{oly}\mathcal{E}\mathrm{nd}(I)$ the category of polynomial endofunctors over the set $I$, where morphisms as above are required to satisfy $f = \mathrm{id}_I$.

A *polynomial tree* is a polynomial endofunctor $T$ such that

(1) $I_T$ (also called the set of *edges* in this case), $B_T$ (also called the set of *nodes*), and $E_T$ are finite, and moreover $I_T \neq \varnothing$,
(2) $s$ and $t$ are injective, and the complement of the image of $s$ has a unique element $r \in I_T$ which we call the *root edge*,
(3) for every edge $i \in I_T$, there exists $k \in \mathbb{N}$ such that $\sigma^k(i) = r$,

where the *walk-to-root* function $\sigma : I_T \longrightarrow I_T$ is inductively defined by $\sigma(r) = r$, and, for $i \neq r \in I_T$, $\sigma(i) = tp(s^{-1}(i))$. An edge that is not in the image of $\sigma$ is called a *leaf*. Let $\mathcal{T}\mathrm{ree}$ be the full subcategory of $\mathcal{P}\mathrm{oly}\mathcal{E}\mathrm{nd}$ spanned by trees. We emphasize that morphisms in $\mathcal{T}\mathrm{ree}$ are not required to be identities on edges, unlike what would happen if $\mathcal{T}\mathrm{ree}$ was a subcategory of $\mathcal{P}\mathrm{oly}\mathcal{E}\mathrm{nd}(I)$.

For $F \in \mathcal{P}\mathrm{oly}\mathcal{E}\mathrm{nd}$, let $\mathrm{tr}\, F$, the category of $F$*-trees*, be a chosen skeleton of $\mathcal{T}\mathrm{ree}/F$. Depending on the context, it may also refer to the set of isomorphism classes of $F$-trees. An $F$-tree is thus a tree $\langle T \rangle$ and a morphism of polynomial functors $T : \langle T \rangle \to F$. For $x$ an edge of $T$ (or more precisely of $\langle T \rangle$), the morphism $T$ associates a color $i = T(x) \in I_F$, and we say that $x$ *is decorated with* $i$, and likewise for nodes.

Let $T \in \mathrm{tr}\, F$, and $\sigma$ the walk-to-root function of $\langle T \rangle$. Write $\langle T \rangle$ as

$$I_T \xleftarrow{s} E_T \xrightarrow{p} B_T \xrightarrow{t} I_T.$$

For an edge $x \in I_T$ we define inductively the *address* of $x$, noted $\&x$, which is a bracket-enclosed list of elements of $E_T$, indicating the path from the root to that edge. If $x$ is the root edge, then $\&x = []$, the empty list. Otherwise, let $\&x = (\&\sigma(x)) \cdot [x]$, where $\cdot$ is the concatenation of lists. If $y \in B_T$ is a node of $T$, let its address be given by $\&y = \&t(y)$, that is, it is the address of its output edge. Write $T^\bullet$ the set of node addresses of $T$, and for $[p] \in T^\bullet$, let $\mathsf{s}_{[p]} T = T(y)$, where $y \in B_T$ is the node such that $\&y = [p]$. In other words $\mathsf{s}_{[p]} T \in B_F$ is the decoration of the node of $T$ at address $[p]$. We let $T^|$ be the set of leaf addresses of $T$.

Let $F$ be a polynomial endofunctor as in equation (2.1.1). For $i \in I$, define $\mathsf{I}_i \in \operatorname{tr} F$ as having underlying tree

$$\{i\} \longleftarrow \varnothing \longrightarrow \varnothing \longrightarrow \{i\}, \tag{2.1.2}$$

and the obvious morphism to $F$. This corresponds uniquely to the $F$-tree with no node and a unique edge, decorated by $i$. Equivalently, it corresponds uniquely to the edge address of an edge decorated in $i$. Let now $b \in B$, and define $\mathsf{Y}_b \in \operatorname{tr} F$, the *corolla* at $b$, as having underlying tree

$$E(b) + \{t(b)\} \longleftarrow E(b) \longrightarrow \{b\} \xrightarrow{b} E(b) + \{t(b)\}, \tag{2.1.3}$$

where the left map is increasing with codomain $\{1, \ldots, n\}$, and the right one maps $b$ to $n+1$. This corresponds to a $F$-tree with a unique node, decorated by $b$, or equivalently, to the node address of such a node.

Let $S, T \in \operatorname{tr} F$ be two $F$-trees, where the root edge (that at address $[]$) of $T$ is decorated by $i \in I_F$, and let $[l] \in S^|$ be the address of a leaf of $S$ decorated by $i$. We define the *grafting* $S \circ_{[l]} T$ of $S$ and $T$ on $l$ by the following pushout (in $\operatorname{tr} F$):

$$\begin{array}{ccc} \mathsf{I}_i & \xrightarrow{r} & T \\ {\scriptstyle l}\downarrow & \ulcorner & \downarrow \\ S & \longrightarrow & S \underset{[l]}{\circ} T. \end{array} \tag{2.1.4}$$

**Proposition 2.1.5** ([**Kock, 2011**]). *Every $F$-tree is either of the form $\mathsf{I}_i$, for some $i \in I$, or obtained by graftings a corollas on an $F$-tree.*

Take $T, U_1, \ldots, U_k \in \operatorname{tr} F$, where the leaves of $T$ have addresses $[l_1], \ldots, [l_k]$, and assume the grafting $T \circ_{[l_i]} U_i$ is defined for all $i$. Then the *total grafting* will be denoted concisely by

$$T \underset{[l_i]}{\bigcirc} U_i = (\cdots(T \underset{[l_1]}{\circ} U_1) \underset{[l_2]}{\circ} U_2 \cdots) \underset{[l_k]}{\circ} U_k. \tag{2.1.6}$$

It is easy to see that the result does not depend on the order in which the graftings are performed.

2.1.1. **The polynomial Baez–Dolan $(-)^+$ construction.** We now give a brief definition of the $(-)^+$ construction, see appendix A on page 71 for details. A *polynomial monad* $(M, \mu, \eta)$ is simply a monad in the 2-category $\mathcal{P}oly\mathcal{E}nd$. This structure induces a function $\mathsf{t} : \operatorname{tr} M \longrightarrow B_M$, for each $T \in \operatorname{tr} M$ a bijection $\wp_T : T^| \longrightarrow (\mathsf{t} T)^\bullet$, called *readdressing*, see appendix A.1 on page 71.

From such a monad $M$, we may define a new polynomial monad $M^+$ whose underlying functor is given by

$$B_M \xleftarrow{\mathsf{s}} \operatorname{tr}^\bullet M \xrightarrow{u} \operatorname{tr} M \xrightarrow{\mathsf{t}} B_M$$

where $\operatorname{tr}^\bullet M$ is the set of $M$-trees with one marked node (or equivalently, equipped with one of their node addresses), $\mathsf{s}$ gives the marked node (or more precisely its decoration), $u$ forgets the marking, and $\mathsf{t}$ is induced from the monad law of $M$. The monad structure of $M^+$ is given by substitution of trees, see appendix A.2 on page 73.

## 2.2. Opetopes

2.2.1. **Definition.** Let $\mathfrak{Z}^0$ be the identity polynomial monad on $\mathcal{S}et$, on the left, and define $\mathfrak{Z}^{n+1} := (\mathfrak{Z}^n)^+$.

$$\{*\} \longleftarrow \{*\} \longrightarrow \{*\} \longrightarrow \{*\}, \qquad \mathbb{O}_n \xleftarrow{\mathsf{s}} E_{\mathfrak{Z}^{n+1}} \longrightarrow \mathbb{O}_{n+1} \xrightarrow{\mathsf{t}} \mathbb{O}_n.$$

If we write the polynomial functor $\mathfrak{Z}^n$ as on the right, an $n$-opetope $\omega$ is by definition an element of $\mathbb{O}_n$, and we say that $\omega$ has *dimension* $n$, written $\dim \omega = n$. The map $\mathsf{s}$ is called the *source* map, and $\mathsf{t}$ the *target* map.

The unique 0-opetope is written $\blacklozenge$ and called the *point*, whereas the unique 1-opetope is written $\blacksquare$ and called the *arrow*. For $n \geq 2$, we have $\mathbb{O}_n = \operatorname{tr} \mathfrak{Z}^{n-2}$, hence an $n$-opetope is a tree whose nodes and edges are decorated with $(n-1)$-opetopes and $(n-2)$-opetopes, respectively. For $\omega \in \mathbb{O}_n$, the set $E_{\mathfrak{Z}^n}(\omega)$ is precisely $\omega^\bullet$, i.e. is the set of node addresses of $\omega$. The readdressing of $\omega$ is then of the form $\wp_\omega : \omega^| \longrightarrow (\mathsf{t}\omega)^\bullet$.

2.2.2. **Opetopes vs. pasting schemes.** Opetopes are closely related to the notion of "pasting scheme" commonly used in higher category theory to describe composition of higher cells. Informally, a pasting scheme of dimension $n$ is a tree whose nodes are decorated with $n$-"cells", edges with $(n-1)$-cells, and where the output edge of a node corresponds to its target, and the input edges to the cells in its source. For instance, the figure on the left represents a 2-pasting scheme, and the corresponding tree is drawn on the right:

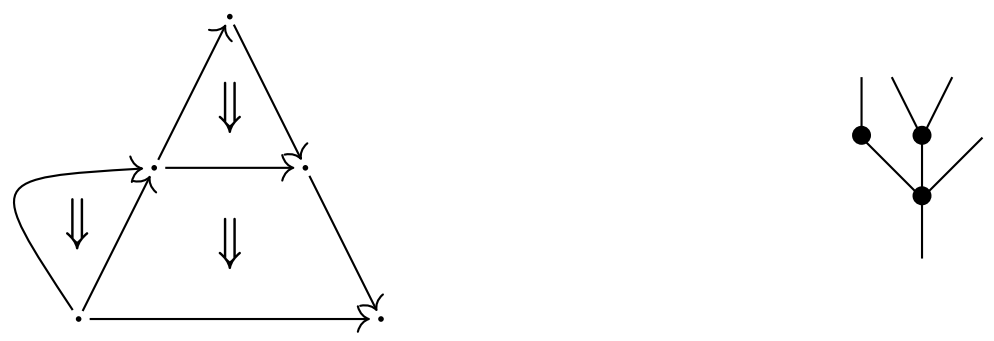

If we consider the $k$-cells as $k$-opetopes for all $k$, then an $n$-pasting scheme $P$ induces a $\mathfrak{Z}^{n-1}$-tree, thus an $(n+1)$-opetope, say $\omega$. We say that $P$ is the *source pasting scheme* of $\omega$. Further, by definition of $\mathfrak{Z}^n$, the opetope $\omega$ has a *target* $\mathsf{t}\,\omega \in \mathbb{O}_{n-1}$. In the sequel, graphical representations of opetopes include both the source pasting scheme and the target: for instance, if $P$ is the pasting scheme above, then $\omega$ is represented by

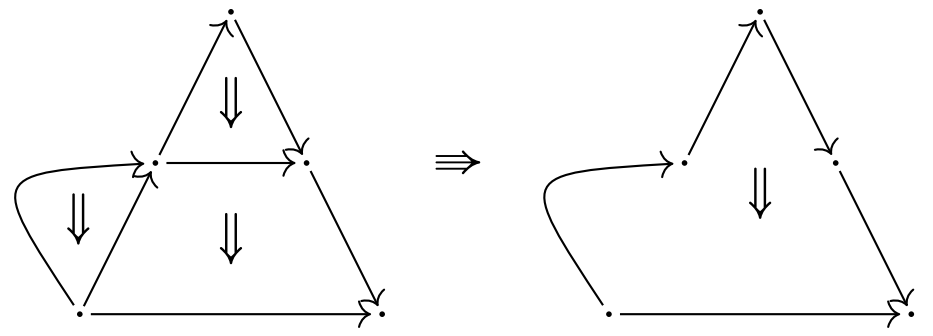

Such a graphical representation is ambiguous however: the following represents a 2-opetope (whose underlying 1-pasting scheme has three 1-cells), but also a 2-pasting scheme having a unique 2-cell. In this work, such ambiguities shall be lifted by the surrounding context.

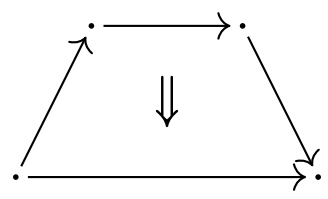

2.2.3. **Higher addresses.** We now introduce the concept of *higher dimensional address*, a refinement if that of addresses presented in section 2.1 on page 9, which will serve as a more convenient way to designate nodes and edges of opetopes. Start by defining the set $\mathbb{A}_n$ of *$n$-addresses* as follows:

$$\mathbb{A}_0 = \{*\}, \qquad \mathbb{A}_{n+1} = \operatorname{lists} \mathbb{A}_n,$$

where $\operatorname{lists} X$ is the set of finite lists (or words) on the alphabet $X$. Explicitly, the unique 0-address is $*$ (also written $[]$ by convention), while an $(n+1)$-address is a sequence of $n$-addresses. Such sequences are enclosed by brackets. Note that the address $[]$, associated to the empty word, is in $\mathbb{A}_n$ for all $n \geq 0$. However, the surrounding context will almost always make the notation unambiguous[1]. Here are examples of the words we shall use:

$$[] \in \mathbb{A}_1, \qquad [* * * *] \in \mathbb{A}_1, \qquad [[][*][]] \in \mathbb{A}_2, \qquad [[[[*]]]] \in \mathbb{A}_4 \qquad \ldots$$

For $\omega \in \mathbb{O}$ an opetope, nodes of $\omega$ can be uniquely specified using higher addresses as we now show. In $\mathfrak{Z}^0$, set $E_{\mathfrak{Z}^1}(\blacksquare) = \{*\}$, so that the unique node address of $\blacksquare$ is $* \in \mathbb{A}_0$. For $n \geq 2$, recall that an opetope $\omega \in \mathbb{O}_n$ is a $\mathfrak{Z}^{n-2}$-tree:

$$\omega : \langle\omega\rangle \longrightarrow \mathfrak{Z}^{n-2}.$$

A node $b \in B_{\langle\omega\rangle}$ has an address $\&b \in \operatorname{lists} E_{\langle\omega\rangle}$, which by $\omega : E_{\langle\omega\rangle} \longrightarrow E_{\mathfrak{Z}^{n-1}}$ is mapped to an element of $\operatorname{lists} E_{\mathfrak{Z}^{n-1}}$. By induction, elements of $E_{\mathfrak{Z}^{n-1}}$ are $(n-2)$-addresses, whence $E_{\langle\omega\rangle}(\&b) \in \operatorname{list} \mathbb{A}_{n-2} = \mathbb{A}_{n-1}$. For the induction step, elements of $E_{\mathfrak{Z}^n}(\omega)$ are nodes of $\langle\omega\rangle$, which we identify by their aforementioned $(n-1)$-addresses. Consequently, for all $n \geq 1$ and $\omega \in \mathbb{O}_n$, elements of $E_{\mathfrak{Z}^n}(\omega) = \omega^\bullet$ can be seen as the set of $(n-1)$-addresses of the nodes of $\omega$, and similarily, $\omega^|$ can be seen as the set of $(n-1)$-addresses of edges of $\omega$.

We denote the concatenation of higher addresses by $\cdot$, i.e. if $[p_1]$, $[p_2]$ are $n$-addresses, then $[p_1]\cdot[p_2] := [p_1 p_2]$. Let $\sqsubseteq$ be the prefix order on the set of $n$-addresses: $[p_1] \sqsubseteq [p_2]$ if and only if as sequences, $p_1$ is a prefix of $p_2$. Let

[1] Ambiguity with addresses could be lifted altogether by hinting the dimension as e.g. $[]^1 \in \mathbb{A}_1$, $[* * * *]^1 \in \mathbb{A}_1$, $[[]^1[*]^1[]^1]^2 \in \mathbb{A}_2$, $[[[[*]^1]^2]^3]^4 \in \mathbb{A}_4$. This however makes notations significantly heavier, so we avoid using this convention.

$\leq$ be the lexicographical order on the set of $n$-addresses. Explicitly, it is trivial on $\mathbb{A}_0$, given by the prefix order $\sqsubseteq$ on $\mathbb{A}_1$, and on $\mathbb{A}_n$, seen as the set of finite sequences of elements of $\mathbb{A}_{n-1}$, it is induced from the lexicographical order on $\mathbb{A}_{n-1}$. Remark that $\sqsubseteq$ is always a suborder of $\leq$.

**Example 2.2.1.** Consider the following 2-opetope, called **3**:

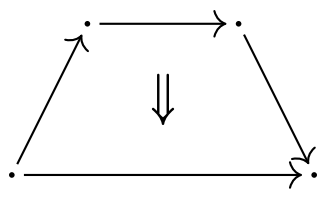

Its underlying pasting diagram consists of 3 arrows $\blacksquare$ grafted linearly. Since the only node address of $\blacksquare$ is $* \in \mathbb{A}_0$, the underlying tree of **3** can be depicted as follows:

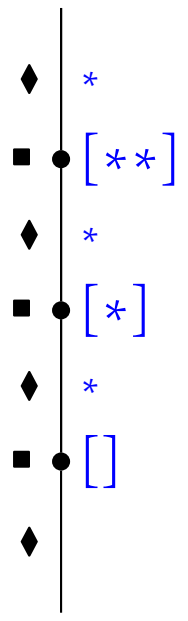

On the left are the decorations: nodes are decorated with $\blacksquare \in \mathbb{O}_1$, while the edges are decorated with $\blacklozenge \in \mathbb{O}_0$. For each node int the tree, the set of input edges of that node is in bijective correspondence with the node addresses of the decorating opetope, and that address is written on the right of each edges. In this low dimensional example, those addresses can only be $*$. Finally, on the right of each node is its 1-address, which is just a sequence of 0-addresses giving "walking instructions" to get from the root to that node.

The 2-opetope **3** can then be seen as a corolla in some 3-opetope as follows:

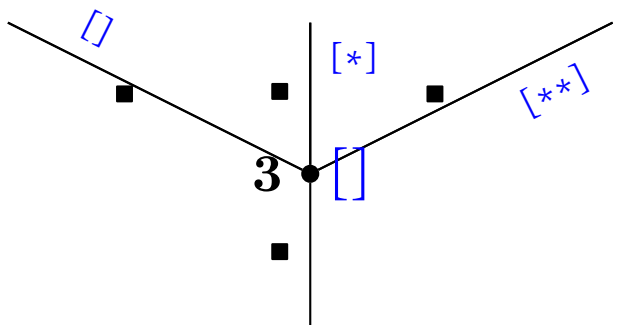

As previously mentioned, the set of input edges is in bijective correspondence with the set of node addresses of **3**. Here is now an example of a 3-opetope, with its annotated underlying tree on the right (the 2-opetopes **1** and **2** are analogous to **3**):

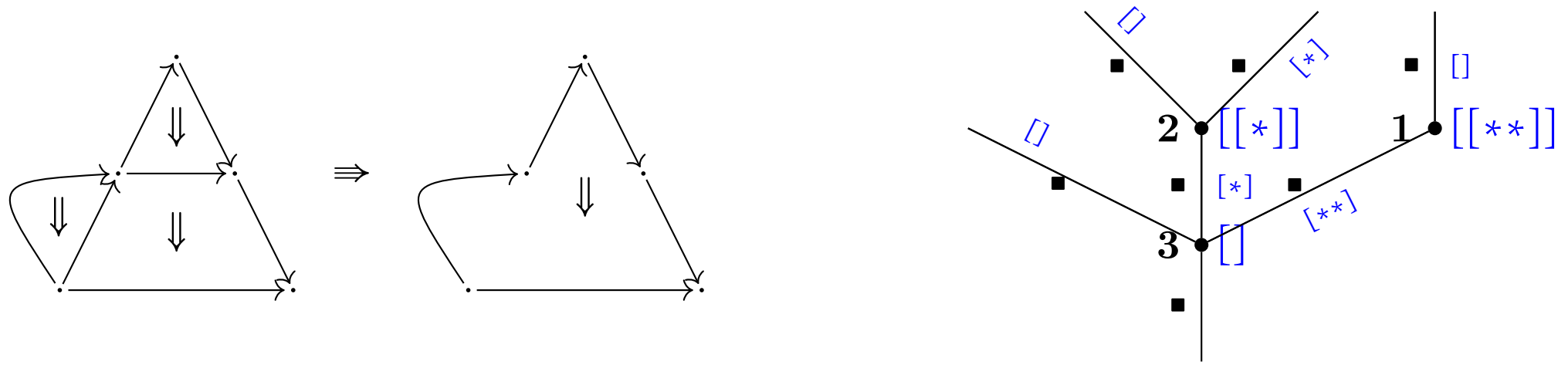

Let $[l] = [p[q]] \in \mathbb{A}_{n-1}$ be an address such that $[p] \in \omega^\bullet$ and $[q] \in \left(\mathsf{s}_{[p]}\,\omega\right)^\bullet$. Then as a shorthand, we write

$$\mathsf{e}_{[l]}\,\omega := \mathsf{s}_{[q]}\,\mathsf{s}_{[p]}\,\omega. \tag{2.2.2}$$

The source map $\mathsf{s}$ give the decoration of a node in a tree, and likewise, $\mathsf{e}$ gives the decoration of edges.

**Remark 2.2.3** (Opetopic induction)**.** The proposition 2.1.5 on page 10 allows a fundamental type of induction on opetopes that we use frequently throughout this paper. Let $\Phi(-)$ be a logical proposition about opetopes such that

(1) $\Phi$ holds for $\blacklozenge$, written $\Phi(\blacklozenge)$,

(2) $\Phi(\psi) \implies \Phi(\mathsf{I}_\psi) \wedge \Phi(\mathsf{Y}_\psi)$, for all $\psi \in \mathbb{O}$,
(3) $\Phi(\nu) \wedge \Phi(\psi) \implies \Phi(\nu \circ_{[l]} \mathsf{Y}_\psi)$, for $\nu, \psi \in \mathbb{O}$ with at least one node, $[l] \in \nu^|$, and assuming the grafting is well-defined.

Then $\Phi$ holds for all opetope.

2.2.4. **The category of opetopes.** In this subsection, we define the category $\mathbb{O}$ of opetopes. We first begin by a result whose geometrical meaning is explained later in this section.

**Lemma 2.2.4** (Opetopic identities, [**Ho Thanh, 2018a**]). *Let $\omega \in \mathbb{O}_n$ with $n \geq 2$.*

*(1) (Inner edge) For $[p[q]] \in \omega^\bullet$ (forcing $\omega$ to be non degenerate), we have $\mathsf{t}\,\mathsf{s}_{[p[q]]}\,\omega = \mathsf{s}_{[q]}\,\mathsf{s}_{[p]}\,\omega$.*
*(2) (Globularity 1) If $\omega$ is non degenerate, we have $\mathsf{t}\,\mathsf{s}_{[]}\,\omega = \mathsf{t}\,\mathsf{t}\,\omega$.*
*(3) (Globularity 2) If $\omega$ is non degenerate, and $[p[q]] \in \omega^|$, we have $\mathsf{s}_{[q]}\,\mathsf{s}_{[p]}\,\omega = \mathsf{s}_{\wp_\omega[p[q]]}\,\mathsf{t}\,\omega$.*
*(4) (Degeneracy) If $\omega$ is degenerate, we have $\mathsf{s}_{[]}\,\mathsf{t}\,\omega = \mathsf{t}\,\mathsf{t}\,\omega$.*

With those identities in mind, we define the category $\mathbb{O}$ of opetopes by generators and relations as follows.

(1) (Objects) We set $\operatorname{ob}\mathbb{O} = \sum_{n \in \mathbb{N}} \mathbb{O}_n$.
(2) (Generators) Let $\omega \in \mathbb{O}_n$ with $n \geq 1$. We introduce a generator, called *target embedding*: $\mathsf{t} : \mathsf{t}\,\omega \longrightarrow \omega$. If $[p] \in \omega^\bullet$, then we introduce a generator, called *source embedding*: $\mathsf{s}_{[p]} : \mathsf{s}_{[p]}\,\omega \longrightarrow \omega$. A *face embedding* is either a source or target embedding.
(3) (Relations) We impose 4 relations described by the following commutative squares, that are well defined thanks to theorem lemma 2.2.4. Let $\omega \in \mathbb{O}_n$ with $n \geq 2$
  (a) **(Inner)** for $[p[q]] \in \omega^\bullet$ (forcing $\omega$ to be non degenerate), the following square must commute:

$$\begin{array}{ccc} \mathsf{s}_{[q]}\,\mathsf{s}_{[p]}\,\omega & \xrightarrow{\mathsf{s}_{[q]}} & \mathsf{s}_{[p]}\,\omega \\ {\scriptstyle \mathsf{t}}\downarrow & & \downarrow{\scriptstyle \mathsf{s}_{[p]}} \\ \mathsf{s}_{[p[q]]}\,\omega & \xrightarrow{\mathsf{s}_{[p[q]]}} & \omega \end{array}$$

  (b) **(Glob1)** if $\omega$ is non degenerate, the following square must commute:

$$\begin{array}{ccc} \mathsf{t}\,\mathsf{t}\,\omega & \xrightarrow{\mathsf{t}} & \mathsf{t}\,\omega \\ {\scriptstyle \mathsf{t}}\downarrow & & \downarrow{\scriptstyle \mathsf{t}} \\ \mathsf{s}_{[]}\,\omega & \xrightarrow{\mathsf{s}_{[]}} & \omega. \end{array}$$

  (c) **(Glob2)** if $\omega$ is non degenerate, and for $[p[q]] \in \omega^|$, the following square must commute:

$$\begin{array}{ccc} \mathsf{s}_{\wp_\omega[p[q]]}\,\mathsf{t}\,\omega & \xrightarrow{\mathsf{s}_{\wp_\omega[p[q]]}} & \mathsf{t}\,\omega \\ {\scriptstyle \mathsf{s}_{[q]}}\downarrow & & \downarrow{\scriptstyle \mathsf{t}} \\ \mathsf{s}_{[p]}\,\omega & \xrightarrow{\mathsf{s}_{[p]}} & \omega. \end{array}$$

  (d) **(Degen)** if $\omega$ is degenerate, the following square must commute:

$$\begin{array}{ccc} \mathsf{t}\,\mathsf{t}\,\omega & \xrightarrow{\mathsf{t}} & \mathsf{t}\,\omega \\ {\scriptstyle \mathsf{s}_{[]}}\downarrow & & \downarrow{\scriptstyle \mathsf{t}} \\ \mathsf{t}\,\omega & \xrightarrow{\mathsf{t}} & \omega. \end{array}$$

An opetopic set is a $\mathcal{S}\mathrm{et}$-valued presheaf over $\mathbb{O}$. We write $\widehat{\mathbb{O}}$ for the category of opetopic sets and natural transformations, $O[-] : \mathbb{O} \longrightarrow \widehat{\mathbb{O}}$ for the Yoneda embedding, and $\mathcal{F}\mathrm{in}\widehat{\mathbb{O}}$ for the full subcategory of $\widehat{\mathbb{O}}$ spanned by finite opetopic sets, i.e. those $X \in \widehat{\mathbb{O}}$ such that $\sum_{\omega \in \mathbb{O}} X_\omega$ is a finite set. Equivalently, $\mathcal{F}\mathrm{in}\widehat{\mathbb{O}}$ is the completion of $\mathbb{O}$ under finite colimits.

Let us explain this definition a little more. As previously mentioned, opetopes are trees whose nodes (and edges) are decorated by opetopes. The decoration is now interpreted as a geometrical feature, namely as an embedding of a lower dimensional opetope. Further, the target of an opetope, while not an intrinsic data, is also represented as an embedding. The relations can be understood as follows.

(1) **(Inner)** The inner edge at $[p[q]] \in \omega^\bullet$ is decorated by the target of the decoration of the node "above" it (here $\mathsf{s}_{[p[q]]}\,\omega$), and in the $[q]$-source of the node "below" it (here $\mathsf{s}_{[p]}\,\omega$). By construction, those two decorations match, and this relation makes the two corresponding embeddings $\mathsf{s}_{[q]}\,\mathsf{s}_{[p]}\,\omega \longrightarrow \omega$ match as well. On the left is an informal diagram about $\omega$ as a tree (reversed gray triangle), and on the right is an example of pasting diagram represented by an opetope, with the relevant features of the **(Inner)** relation colored or thickened.

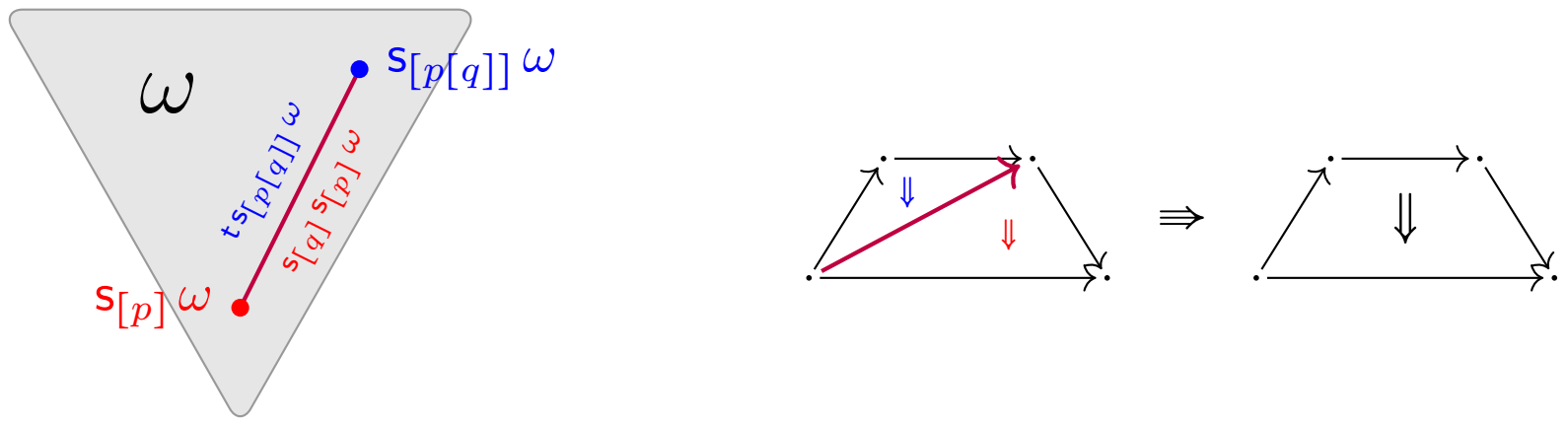

(2) **(Glob1-2)** If we consider the underlying tree of $\omega$ (which really is $\omega$ itself) as its "geometrical source", and the corolla $\mathsf{Y}_{\mathsf{t}\omega}$ as its "geometrical target", then they should be parallel. The relation **(Glob1)** expresses this idea by "gluing" the root edges of $\omega$ and $\mathsf{Y}_{\mathsf{t}\omega}$ together, while **(Glob2)** glues the leaves according to the readdressing function $\wp_\omega$.

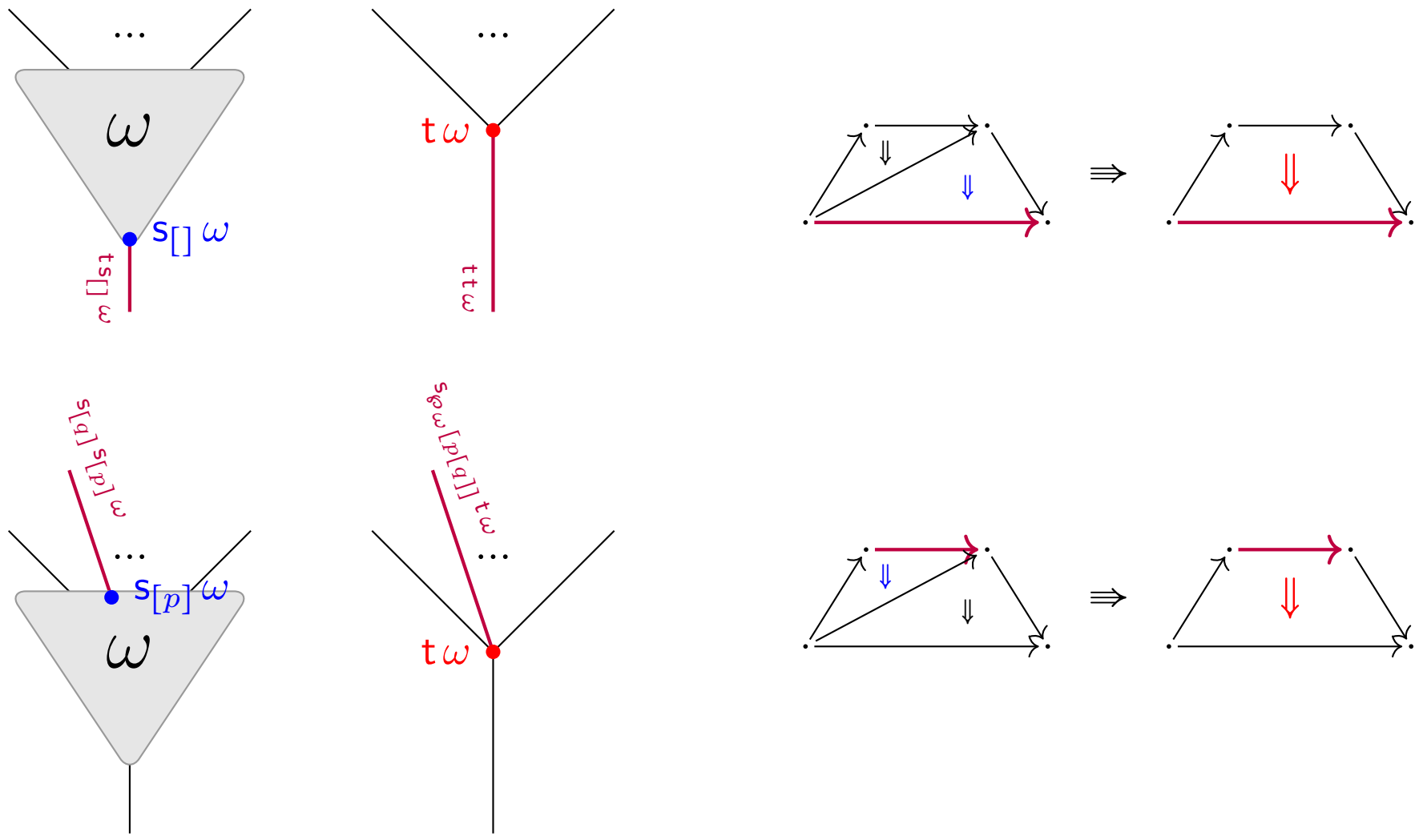

(3) **(Degen)** If $\omega$ is a degenerate opetope, depicted as on the right, then its target should be a "loop", i.e. its only source and target should be glued together.

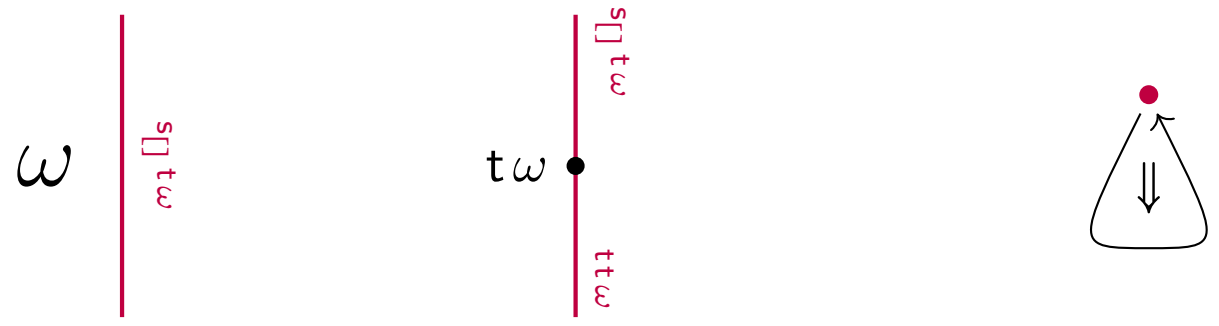

CHAPTER 3

# Named approach

## 3.1. The system for opetopes

3.1.1. **Syntax.** In this section, we define the underlying syntax of $\mathrm{OPT}^!$, our named derivation system for opetopes. As explained in the introduction, a typical pasting scheme is pictured below:

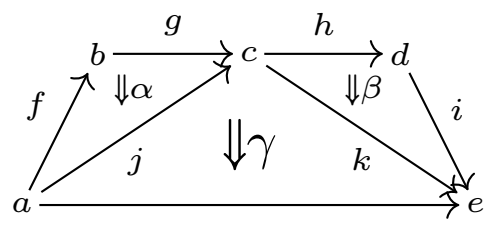

We shall use the names of the cells of this picture as variables, and encode the pasting scheme as the following expression:

$$\gamma(j \leftarrow \alpha, k \leftarrow \beta).$$

Here, $j$, $k$, $\alpha$, $\beta$, and $\gamma$ are now variables, equipped with a dimension (1 for $j$ and $k$, and 2 for $\alpha$ , $\beta$, and $\gamma$), and the notation is meant to be read as "the variable $\gamma$ in which $\alpha$ (resp. $\beta$) has been formally grafted on the input labeled $j$ (resp. $k$)". Such a term will be given a type:

$$i(t \leftarrow h(z \leftarrow g(y \leftarrow f))) \multimap a \multimap \varnothing,$$

which expresses the fact that the source is the "composite" $i \circ h \circ g \circ f$, and that the source of the source is $x$. Since the pasting scheme is 2-dimensional, there is no further iterated source, and we conclude the sequence by a $\varnothing$ symbol (which can be read as the only $(-1)$-dimensional pasting scheme). Similarly, the degenerate pasting scheme on the left below will be denoted by the typed term figured on the right:

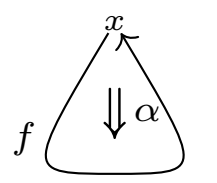

$$\alpha : \underline{x} \multimap x \multimap \varnothing$$

where the term $\underline{x}$ denotes a degenerate 1-dimensional pasting scheme with $x$ as source.

3.1.1.1. *Terms.* We have a $\mathbb{N}$-graded set $\mathbb{V}$ of variables. Elements of $\mathbb{V}_n$ represent $n$-dimensional cells. An $n$-term is constructed according to the following grammar:

$$\begin{aligned} \mathbb{T}_{-1} &::= \{\varnothing\} && \text{by convention} \\ \mathbb{T}_0 &::= \mathbb{V}_0 \\ \mathbb{T}_{n+1} &::= \underline{\mathbb{V}_n} \quad | \quad \mathbb{T}'_{n+1} \\ \mathbb{T}'_{n+1} &::= \mathbb{V}_{n+1}(\mathbb{V}_n \leftarrow \mathbb{T}'_{n+1}, \ldots) \end{aligned}$$

where the expression $(\mathbb{V}_n \leftarrow \mathbb{T}'_{n+1}, \ldots)$ signifies that there is 0 or more instances of the "$\mathbb{V}_n \leftarrow \mathbb{T}'_{n+1}$" part between the parentheses. For example, if $f, g \in \mathbb{V}_1$, and $a \in \mathbb{V}_0$, then the following is an element of $\mathbb{T}_1$:

$$g(a \leftarrow f())$$

To make notations lighter, we do not write empty parentheses "$()$", so the previous 1-term can be more concisely written as $g(a \leftarrow f)$. The terms of the form $\underline{u}$ are called *degenerate terms*, or *empty* syntactic pasting schemes.

A term of the form $g(a_1 \leftarrow f_1, \ldots, a_k \leftarrow f_k)$ will oftentimes be abbreviated as $g(\overrightarrow{a_i \leftarrow f_i})$, leaving $k$ implicit. By convention, the sequence $a_1 \leftarrow f_1, \ldots, a_k \leftarrow f_k$ above is always considered up to permutation: for $\sigma$ a bijection of the set $\{1, \ldots, k\}$, the terms $g(a_1 \leftarrow f_1, \ldots, a_k \leftarrow f_k)$ and $g(a_{\sigma(1)} \leftarrow f_{\sigma(1)}, \ldots, a_{\sigma(k)} \leftarrow f_{\sigma(k)})$ are considered equal. Note that the above definition entails that $\mathbb{T}_0 = \mathbb{V}_0$ (in particular, there are no degenerate 0-terms).

For $t \in \mathbb{T}_n$, write $t^\bullet$ for the set of $n$-variables occurring in $t$. In the previous example, $(g(a \leftarrow f))^\bullet = \{f, g\}$. Note that $x \in x^\bullet$ for all $x \in \mathbb{V}_n$. We consider terms of the form $\underline{y}$, with $y \in \mathbb{V}_k$, as empty, or *degenerate*, $(n+1)$-terms, i.e. with no $(n+1)$-variables.

3.1.1.2. *Types.* An $n$-type $T$ is a sequence of terms of the form

$$s_1 \multimapdotboth s_2 \multimapdotboth \cdots \multimapdotboth s_n \multimapdotboth \varnothing,$$

where $s_i \in \mathbb{T}_{n-i}$. As we will see (rule `shift` in figure 3.1.1 on the next page, and theorem 3.1.8 on page 19), this sequence of terms essentially describes a zoom complex in the sense of [**Kock et al., 2010**], which justifies the use of the $\multimapdotboth$ symbol. A *typing* of a term $t \in \mathbb{T}_n$ is an expression of the form $t : T$, for $T$ an $n$-type. If $T$ is as above, then $s_i$ is thought of as the $i$-th (iterated) source of $t$. We then write $\mathsf{s}\, t := s_1$, and more generally $\mathsf{s}^i\, t := s_i$. By convention, $\mathsf{s}^0\, t := t$.

3.1.1.3. *Contexts.* A context $\Gamma$ is a set of typings, more commonly written as a list. Write $\mathbb{V}_{\Gamma,k}$ for the set of $k$-variables typed in $\Gamma$, let $\mathbb{V}_\Gamma := \sum_{k\in\mathbb{N}} \mathbb{V}_{\Gamma,k}$, write $\mathbb{T}_{\Gamma,k}$ for the set of $k$-terms whose variables (in any dimension) are in $\mathbb{V}_\Gamma$, and $\mathbb{T}_\Gamma := \sum_{k\in\mathbb{N}} \mathbb{T}_{\Gamma,k}$. As we will see (inference rules in figure 3.1.1 on the next page), for a derivable context $\Gamma$, if $x$ occurs in the typing of a variable of $\Gamma$, then $x \in \mathbb{V}_\Gamma$. Note that in any context $\Gamma$, if a variable $x \in \mathbb{V}_{\Gamma,k}$ occurs in the type of $y \in \mathbb{V}_{\Gamma,l}$, then $k < l$, and thus there is no cyclic dependency among variables.

3.1.1.4. *Sequents.* Let $\Gamma$ be a context. An *equational theory* $E$ on $\mathbb{V}_\Gamma$ is a set of formal equalities between variables of $\Gamma$. We write $=_E$ for the equivalence relation on $\mathbb{V}_\Gamma$ generated by $E$. A sequent is an expression of the form

$$E \triangleright \Gamma \vdash t : T$$

where $\Gamma$ is a context, $E$ is an equational theory on $\mathbb{V}_\Gamma$, and the right hand side is a typing. We may write $\vdash_n$ to signify that $t \in \mathbb{T}_n$. The equivalence relation $=_E$ on $\mathbb{V}_\Gamma$ extends to $\mathbb{T}_\Gamma$ in an obvious way. If $x =_E y \in t^\bullet$, then, by convention $x \in t^\bullet$, so that $x$ and $y$ really are interchangeable.

If $(F \triangleright \Upsilon \vdash u : U)$ is a sequent such that there exists a bijection $\sigma : \mathbb{V}_\Upsilon \longrightarrow \mathbb{V}_\Gamma$ with

$$(E \triangleright \Gamma \vdash t : T) = (F^\sigma \triangleright \Upsilon^\sigma \vdash u^\sigma : U^\sigma),$$

where $(-)^\sigma$ is the substitution according to $\sigma$, then we say that both sequents are *equivalent* (or *$\alpha$-equivalent*), denoted

$$(E \triangleright \Gamma \vdash t : T) \simeq (F \triangleright \Upsilon \vdash u : U).$$

In the following, sequents are implicitly considered up to equivalence.

3.1.2. **Inference rules.** We present the inference rules of the $\mathrm{O{\scriptsize PT}}^!$ system in figure 3.1.1 on the facing page. Rule `graft` requires the so-called graft notation and substitution operation, respectively introduced in definitions 3.1.1 and 3.1.3 on the next page.

FIGURE 3.1.1. The $\text{OPT}^!$ system.

**Introduction of points:** This rule introduces 0-cells, also called points. If $x \in \mathbb{V}_0$, then

$$\frac{}{\triangleright x : \varnothing \vdash_0 x : \varnothing} \texttt{point}$$

**Introduction of degeneracies:** This rule derives empty pasting diagrams. If $x \in \mathbb{V}_n$, then

$$\frac{E \triangleright \Gamma \vdash_n x : T}{E \triangleright \Gamma \vdash_{n+1} \underline{x} : x \multimap T} \texttt{degen}$$

**Shift to the next dimension:** This rule takes a term $t$ and introduces a new cell $x$ having $t$ as source. If $x \in \mathbb{V}_{n+1}$ is such that $x \notin \mathbb{V}_\Gamma$, then

$$\frac{E \triangleright \Gamma \vdash_n t : T}{E \triangleright \Gamma, x : t \multimap T \vdash_{n+1} x : t \multimap T} \texttt{shift}$$

**Grafting:** This rule glues an $n$-cell $x$ onto an $n$-term $t$ along a variable $a \in s_1^\bullet := (\mathsf{s}\, t)^\bullet$. We assume that $\Gamma$ and $\Upsilon$ are compatible, in that for all $y \in \mathbb{V}$, if $y \in \mathbb{V}_\Gamma \cap \mathbb{V}_\Upsilon$, then the typing of $y$ in both contexts match modulo the equational theory $E \cup F$. Further, the only variables typed in both $\Gamma$ and $\Upsilon$ are $a$ and the variables occurring in the sources of $a$ (i.e. $\mathsf{s}^i\, a$, for $1 \leq i \leq n-1$).

If $x \in \mathbb{V}_n$, $t \in \mathbb{T}_n$ is not degenerate, $a \in (\mathsf{s}\, t)^\bullet$ is such that $\mathsf{s}\, a = \mathsf{s}\,\mathsf{s}\, x$, then

$$\frac{E \triangleright \Gamma \vdash_n t : s_1 \multimap s_2 \multimap \cdots \qquad F \triangleright \Upsilon \vdash_n x : U}{G \triangleright \Gamma \cup \Upsilon \vdash_n t(a \leftarrow x) : s_1[\mathsf{s}\, x/a] \multimap s_2 \multimap \cdots} \texttt{graft}$$

where the notations $t(a \leftarrow x)$ and $s_1[\mathsf{s}\, x/a]$ are presented below, where $G$ is the union of $E$, $F$, and potentially a set of additional equalities incurred by the substitution $s_1[\mathsf{s}\, x/a]$. We also write $\texttt{graft-}a$ to make explicit that we grafted onto $a$.

The condition $a \in (\mathsf{s}\, t)^\bullet$ ensures that $a$ hasn't been used for grafting beforehand, while the condition $\mathsf{s}\, a =_E \mathsf{s}\,\mathsf{s}\, x$ shows that $x$ may indeed be glued onto $a$.

**Definition 3.1.1** (Graft notation)**.** For a sequent $(E \triangleright \Gamma \vdash_n t : T)$, $a \in \mathbb{V}_{n-1}$, and $x \in \mathbb{V}_{\Upsilon,n}$, the graft notation $t(a \leftarrow x)$ of the `graft` rule can be simplified depending on the structure of $t$, according to the following rewriting rule: for $y \in \mathbb{V}_{\Gamma,n}$:

$$y(\overrightarrow{z_i \leftarrow v_i})(a \leftarrow x) \quad \rightsquigarrow \quad \begin{cases} y(\overrightarrow{z_i \leftarrow v_i(a \leftarrow x)}) & \text{if } a \notin (\mathsf{s}\, y)^\bullet, \\ y(\overrightarrow{z_i \leftarrow v_i}, a \leftarrow x) & \text{if } a \in (\mathsf{s}\, y)^\bullet, \end{cases} \tag{3.1.2}$$

In particular, note that if $a \notin (\mathsf{s}\, y)^\bullet$, then $y()(a \leftarrow x) \rightsquigarrow y()$, and with the "empty parentheses convention", this gives $y(a \leftarrow x) \rightsquigarrow y$.

**Definition 3.1.3** (Substitution)**.** We now explain how to evaluate $u[w/a]$ for any term $u \in \mathbb{T}_{n-1}$. We let

$$u[w/a] := \begin{cases} y(\overrightarrow{z_i \leftarrow v_i[w/a]}) & \text{if } a \neq_{E \cup F} y, \\ w(\overrightarrow{z_i \leftarrow v_i}) & \text{if } a =_{E \cup F} y. \end{cases}$$

and then, in the case $w = \underline{b}$, we perform the following actions:

(1) for subterms of the form $x(\ldots, z \leftarrow \underline{b}, \ldots)$, we remove the grafting $z \leftarrow \underline{b}$, and add the equation $b = z$ to the ambient equational theory,
(2) for subterms of the form $x(\ldots, z \leftarrow \underline{b}(b \leftarrow r), \ldots)$, we replace the grafting $z \leftarrow \underline{b}(b \leftarrow r)$ by $z \leftarrow r$ and add $b = z$ to the ambient equational theory,
(3) and finally replace any remaining subterms of the form $\underline{b}(b \leftarrow v)$ by $v$.

More explicitly:

(1) If $w$ is not an empty pasting diagram (i.e. not of the form $\underline{b}$), then writing $u = y(\overrightarrow{z_i \leftarrow v_i})$:

$$u[w/a] := \begin{cases} y(\overrightarrow{z_i \leftarrow v_i[w/a]}) & \text{if } a \neq_{E \cup F} y, \\ w(\overrightarrow{z_i \leftarrow v_i}) & \text{if } a =_{E \cup F} y. \end{cases} \tag{3.1.4}$$

(2) If $w$ is an empty pasting diagram, say $w = \underline{b}$ for $b \in \mathbb{V}_{n-2}$, Then, by the hypothesis of the `graft` rule, we have $b =_E \mathsf{s}\, a$. Then, $u[\underline{b}/a]$ is defined by cases on the form of $u$:
 (a) if $u =_{E \cup F} a$, then $u[\underline{b}/a] := \underline{b}$;
 (b) if $u$ is of the form $a(b \leftarrow r)$, then $u[\underline{b}/a] := r$;
 (c) if $u$ is of the form $y(\ldots, z \leftarrow a, \ldots)$, then

$$u[\underline{b}/a] := y(\ldots, \ldots),$$

 and we add the equality $b = z$ to the ambient equational theory;
 (d) if $u$ is of the form[1] $y(\ldots, z \leftarrow a(b \leftarrow r), \ldots)$, then

$$u[\underline{b}/a] := y(\ldots, z \leftarrow r, \ldots),$$

 as in equation (3.1.4) on the preceding page, and per equation (3.1.2) on the previous page, but we also add the equality $b = z$ to the ambient equational theory;
 (e) otherwise, if $u$ is of the form $y(\overrightarrow{z_i \leftarrow v_i})$, and if the previous cases do not apply (i.e. $a$ is not the front variable of $u$ or $v_i$ for all $i$), then

$$u[\underline{b}/a] := y(\overrightarrow{z_i \leftarrow v_i[\underline{b}/a]}),$$

 as in equation (3.1.4) on the preceding page.

From the formulation of system $\text{Opt}^!$, it is clear that a sequent that is equivalent to a derivable one is itself derivable. Let us now turn our attention to rule `shift` above. It takes a term $t$, thought of as a pasting diagram, and creates a new variable having $t$ as source. One may thus think of it as a rule creating "fillers", akin to Kan filler condition on simplicial sets.

**Example 3.1.5.** Consider the term $t = \alpha(g \leftarrow \beta)$ in a suitable context $\Gamma$:

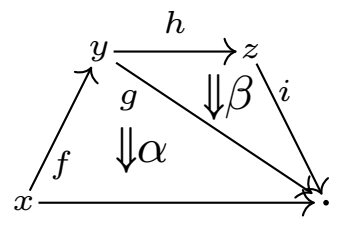

Then the graftings $t(f \leftarrow \gamma) = \alpha(f \leftarrow \gamma, g \leftarrow \beta)$ and $t(i \leftarrow \gamma') = \alpha(g \leftarrow \beta(i \leftarrow \gamma))$, for some appropriate $\gamma$ and $\gamma'$, can respectively be represented as:

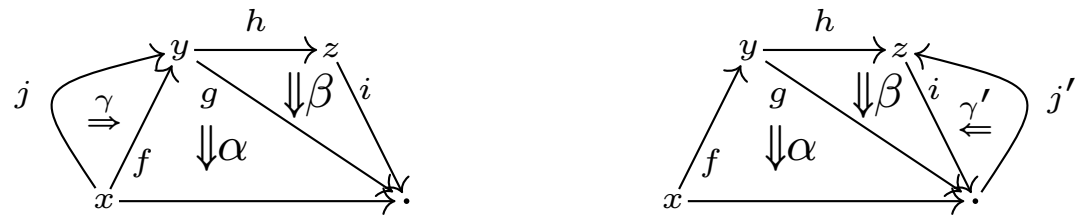

**Example 3.1.6.** Consider the variables $\alpha$, $\beta_x$, and $\beta_y$ typed as:

$$\alpha : g(y \leftarrow f) \multimapdotboth x \multimapdotboth \varnothing, \qquad \beta_x : \underline{x} \multimapdotboth x \multimapdotboth \varnothing, \qquad \beta_y : \underline{y} \multimapdotboth y \multimapdotboth \varnothing.$$

Then $\alpha$, $\alpha(f \leftarrow \beta_x)$, and $\alpha(g \leftarrow \beta_y)$ can respectively be represented as:

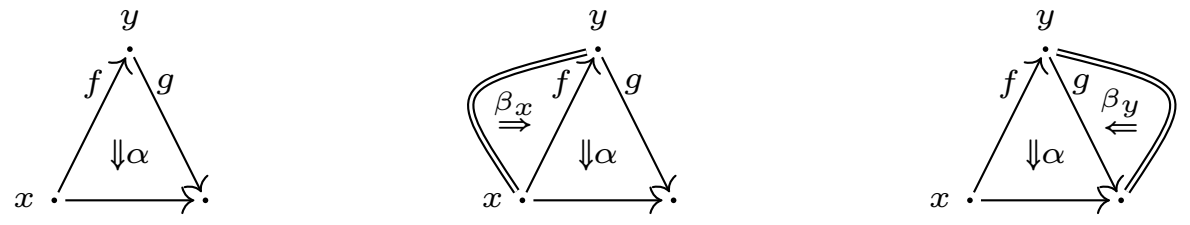

Then the sources $\alpha(f \leftarrow \beta_x)$ and $\alpha(g \leftarrow \beta_y)$ are respectively

$$g(y \leftarrow f)[\underline{x}/f] = g(y \leftarrow \underline{x}) = g, \qquad g(y \leftarrow f)[\underline{y}/g] = \underline{y}(y \leftarrow f) = f,$$

and, in the first case, the equation $x = y$ is added to the ambient equational theory.

**Remark 3.1.7.** The *degen* rule may be replaced by the following *degen-shift* rule without changing the set of derivable sequents of the form $(E \triangleright \Gamma \vdash y : T)$ with $y \in \mathbb{V}$: if $x \in \mathbb{V}_n$ and $d \in \mathbb{V}_{n+2}$ such that $d \notin \mathbb{V}_{\Gamma, n+2}$, then

$$\frac{E \triangleright \Gamma \vdash_n x : T}{E \triangleright \Gamma, d : \underline{x} \multimapdotboth x \multimapdotboth T \vdash_{n+2} d : \underline{x} \multimapdotboth x \multimapdotboth T} \; \textit{degen-shift}$$

However, note that sequents of the form $\left(E \triangleright \Gamma \vdash \underline{y} : T\right)$ are no longer derivable.

[1]The case in which $u$ is of the form $y(\ldots, z \leftarrow a(b_1 \leftarrow r_1, \ldots, b_k \leftarrow r_k), \ldots)$ with $k > 1$ does not happen in valid derivations.

3.1.3. **Uniqueness of typing.** Let $(E \triangleright \Gamma \vdash x : X)$ be a derivable sequent. We prove in theorem 3.1.8 that the type $X = (\mathsf{s}x \multimap \mathsf{s}\mathsf{s}x \multimap \cdots)$ is completely determined by $\mathsf{s}x$ and $\Gamma$.

A consequence of this result is that at any stage, a context $\Gamma$ may be replaced by its "meager form" $\bar{\Gamma}$, obtained by replacing "full typings" $y : Y$ by $y : \mathsf{s}y$, i.e. by removing all but the top term of $Y$. Using meager context comes with a cost however: checking the hypothesis of rule `graft` requires to compute the second source $\mathsf{s}\mathsf{s}x$ of $x$, which is not contained in $\bar{\Gamma}$. For clarity, we do not make use of meager forms throughout the rest of this work.

Define the function $\bar{\mathsf{s}}$ as follows:

$$\begin{aligned}
\bar{\mathsf{s}} : \mathbb{T}_\Gamma &\longrightarrow \mathbb{T}_\Gamma \\
x &\longmapsto \mathsf{s}x && x \in \mathbb{V}_\Gamma, \\
\underline{x} &\longmapsto x && x \in \mathbb{V}_\Gamma, \\
x(\overrightarrow{y_i \leftarrow u_i}) &\longmapsto (\mathsf{s}x)[\overrightarrow{\bar{\mathsf{s}}u_i / y_i}] && x, \overrightarrow{y_i} \in \mathbb{V}_\Gamma, \overrightarrow{u_i} \in \mathbb{T}_\Gamma.
\end{aligned}$$

**Theorem 3.1.8.** *Let $(E \triangleright \Gamma \vdash t : s_1 \multimap s_2 \multimap \cdots \multimap s_n \multimap \varnothing)$ be a derivable sequent. Then for $1 \leq k \leq n$ we have $s_k = \bar{\mathsf{s}}^k t$, or equivalently, for $0 \leq i \leq n$, we have $\bar{\mathsf{s}} s_i = s_{i+1}$, with $s_0 := t$.*

*Proof.* We proceed by induction on the proof tree of the sequent. For readability, we omit equational theories and contexts.

(1) If the sequent is obtained by the following proof tree:
$$\frac{}{x : \varnothing}\ \texttt{point}$$
then $\bar{\mathsf{s}}x = \mathsf{s}x = \varnothing$, since $x \in \mathbb{V}$.

(2) If the sequent is obtained by the following proof tree:
$$\frac{s_1 : s_2 \multimap \cdots \multimap s_n \multimap \varnothing}{t : s_1 \multimap s_2 \multimap \cdots \multimap s_n \multimap \varnothing}\ \texttt{degen}$$
then $s_1 \in \mathbb{V}$ and $t = \underline{s_1}$. Thus, $\bar{\mathsf{s}}t = s_1$, while for $1 \leq i \leq n$, the equality $\bar{\mathsf{s}}s_i = s_{i+1}$ holds by induction.

(3) If the sequent is obtained by the following proof tree:
$$\frac{s_1 : s_2 \multimap \cdots \multimap s_n \multimap \varnothing}{t : s_1 \multimap s_2 \multimap \cdots \multimap s_n \multimap \varnothing}\ \texttt{shift}$$
then $t \in \mathbb{V}$, so $\bar{\mathsf{s}}t = \mathsf{s}t = s_1$, while for $1 \leq i \leq n$, the equality $\bar{\mathsf{s}}s_i = s_{i+1}$ holds by induction.

(4) Assume now that the sequent is obtained by the following proof tree:
$$\frac{u : r_1 \multimap r_2 \multimap s_3 \multimap \cdots \multimap s_n \multimap \varnothing \qquad x : X}{t : s_1 \multimap s_2 \multimap \cdots \multimap s_n \multimap \varnothing}\ \texttt{graft-}a$$
with $a \in r_1^\bullet$ and $x \in \mathbb{V}$ such that $\mathsf{s}\mathsf{s}x = \mathsf{s}a$. Then $t = u(a \leftarrow x)$, $s_1 = r_1[\mathsf{s}x/a]$, and $s_i = r_i$ for $2 \leq i \leq n$. On the one hand, we have
$$\begin{aligned}
\bar{\mathsf{s}}t &= \bar{\mathsf{s}}\,(u(a \leftarrow x)) \\
&= \bar{\mathsf{s}}u[\bar{\mathsf{s}}x/a] \\
&= \bar{\mathsf{s}}u[\mathsf{s}x/a] && \text{since } x \in \mathbb{V} \\
&= s_1.
\end{aligned}$$
On the other hand, write $r_1 = v(y \leftarrow a(\overrightarrow{z_i \leftarrow w_i}))$, for some $v, \overrightarrow{w_i} \in \mathbb{T}_{n-1}$ and $y \in \mathbb{V}_{n-2}$. Then
$$\begin{aligned}
\bar{\mathsf{s}}s_1 &= \bar{\mathsf{s}}\,(r_1[\mathsf{s}x/a]) \\
&= \bar{\mathsf{s}}\left(v(y \leftarrow (\mathsf{s}x)(\overrightarrow{z_i \leftarrow w_i}))\right) \\
&= (\bar{\mathsf{s}}v)[(\bar{\mathsf{s}}\mathsf{s}x)[\overrightarrow{\bar{\mathsf{s}}w_i/z_i}]/y] \\
&= (\bar{\mathsf{s}}v)[(\mathsf{s}\mathsf{s}x)[\overrightarrow{\bar{\mathsf{s}}w_i/z_i}]/y] && \text{by ind.} \\
&= (\bar{\mathsf{s}}v)[(\mathsf{s}a)[\overrightarrow{\bar{\mathsf{s}}w_i/z_i}]/y] && \text{hyp. of } \texttt{graft-}a \\
&= \bar{\mathsf{s}}r_1 = r_2 = s_2.
\end{aligned}$$
Finally, for $1 \leq i \leq n$, the equality $\bar{\mathsf{s}}s_i = s_{i+1}$ holds by induction. □

**Corollary 3.1.9.** *Let* $(E \triangleright \Gamma \vdash t : T)$ *be a derivable sequent, and* $x : s_1 \multimapdotboth s_2 \multimapdotboth \cdots \multimapdotboth s_n \multimapdotboth \varnothing$ *be a typing in* $\Gamma$. *Then for* $1 \leq k \leq n$ *we have* $s_k = \bar{\mathsf{s}}^k t$, *or equivalently, for* $0 \leq i \leq n$, *we have* $\bar{\mathsf{s}} s_i = s_{i+1}$, *with* $s_0 := x$.

*Proof.* If $x : s_1 \multimapdotboth s_2 \multimapdotboth \cdots \multimapdotboth s_n \multimapdotboth \varnothing$ is a typing in $\Gamma$, then somewhere in the proof tree of $(E \triangleright \Gamma \vdash t : T)$ appears a sequent of the form $(F \triangleright \Upsilon \vdash x : s_1 \multimapdotboth s_2 \multimapdotboth \cdots \multimapdotboth s_n \multimapdotboth \varnothing)$, which is necessarily derivable. We conclude by applying theorem 3.1.8 on the preceding page □

By definition, $\bar{\mathsf{s}}$ extends $\mathsf{s}$ to a function $\mathbb{T}_\Gamma \longrightarrow \mathbb{T}_\Gamma$, and for convenience, we just write it as $\mathsf{s}$ in the sequel, and call it the *source* of a term.

**Example 3.1.10.** Consider the following term:

$$\alpha(f \leftarrow \gamma, g \leftarrow \beta) \qquad =$$

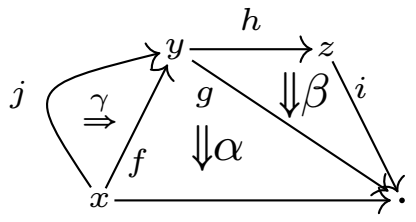

Then its source is computed as follows:

$$\begin{aligned}
\mathsf{s}\,(\alpha(f \leftarrow \gamma, g \leftarrow \beta)) &= (\mathsf{s}\,\alpha)\,[(\mathsf{s}\,\gamma)/f, (\mathsf{s}\,\beta)/g] \\
&= (g(y \leftarrow f))\,[(\mathsf{s}\,\gamma)/f, (\mathsf{s}\,\beta)/g] && \text{since } \mathsf{s}\,\alpha = g(y \leftarrow f) \\
&= (g(y \leftarrow f))\,[j/f, (\mathsf{s}\,\beta)/g] && \text{since } \mathsf{s}\,\gamma = j \\
&= (g(y \leftarrow f))\,[j/f, i(z \leftarrow h)/g] && \text{since } \mathsf{s}\,\beta = i(z \leftarrow h) \\
&= (g(y \leftarrow j))\,[i(z \leftarrow h)/g] \\
&= (i(z \leftarrow h))\,(y \leftarrow j) \\
&= i(z \leftarrow h(y \leftarrow j)) && \text{since } y \in (\mathsf{s}\,h)^\bullet.
\end{aligned}$$

The latter term indeed corresponds to the source of the pasting diagram, i.e. the arrow composition on the top.

## 3.2. Equivalence with polynomial opetopes

In this section, all sequents are assumed derivable in $\text{Opt}^!$. We show that sequents typing a variable (up to $\simeq$) are in bijective correspondence with polynomial opetopes (see section 2.2 on page 10). To this end, we define the *polynomial coding* operation $[\![-]\!]^{\text{poly}}_{n+1}$ that maps a sequent $(E \triangleright \Gamma \vdash_n t : T)$, with $t \in \mathbb{T}_n$, to an $(n+1)$-opetope $[\![E \triangleright \Gamma \vdash_n t : T]\!]^{\text{poly}}_{n+1} \in \mathbb{O}_{n+1}$, written $[\![t : T]\!]^{\text{poly}}_{n+1}$ or even $[\![t]\!]^{\text{poly}}_{n+1}$ for short, if no ambiguity arises. Also, if $\alpha \in \mathbb{V}_{n+1}$, then we set $[\![\alpha]\!]^{\text{poly}}_{n+1} := [\![\mathsf{s}\,\alpha]\!]^{\text{poly}}_{n+1}$.

The idea of the polynomial coding is to map a pasting diagram described by a term (on the left) to its underlying composition tree, and reapply the coding recursively (on the right):

$$[\![\alpha(g \leftarrow \beta)]\!] \qquad =$$

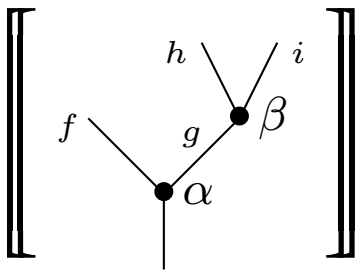

$$:=$$

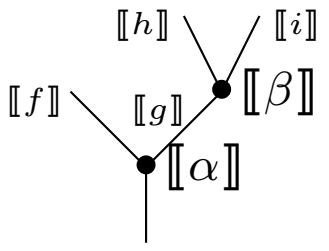

where $[\![-]\!]$ is a shorthand for $[\![-]\!]^{\text{poly}}$.

For $t = x(\overrightarrow{y_i \leftarrow u_i}) \in \mathbb{T}_n$ and $z \in t^\bullet$, the *address* $\&_t z \in \mathbb{A}_n$ of $z$ in $t$ is an $n$-address (see section 2.2.3 on page 11) that indicates "where $z$ is located in $t$".

**Definition 3.2.1** (Address)**.** Take $t \in \mathbb{T}_\Gamma$, $t = x(\overrightarrow{y_i \leftarrow u_i})$.

(1) For $z \in t^\bullet$, the *address* $\&_t z \in \mathbb{A}_n$ of $z$ in $t$ is given by

$$\&_t z := \begin{cases} [] & \text{if } z = x \text{ in the current eq. th.,} \\ [\&_{\mathsf{s}\,x} y_i] \cdot \&_{u_i} z & \text{if } z \in u_i^\bullet, \end{cases}$$

If $[p] = \&_t z$, then we write $\mathsf{s}_{[p]}\, t := z$. In particular, $\mathsf{s}_{[]}\, t = x$.

(2) For $a \in (\mathsf{s}\,t)^\bullet$, the *address* $\&_t a \in \mathbb{A}_n$ of $a$ in $t$ is given by

$$\&_t a := \begin{cases} [\&_{\mathsf{s}\,x} a] & \text{if } a \in (\mathsf{s}\,x)^\bullet, \\ [\&_{\mathsf{s}\,x} y_i] \cdot \&_{u_i} a & \text{if } a \notin (\mathsf{s}\,x)^\bullet, \text{ but } a \in (\mathsf{s}\,u_i)^\bullet. \end{cases}$$

**Example 3.2.2.** The context describing the following pasting scheme:

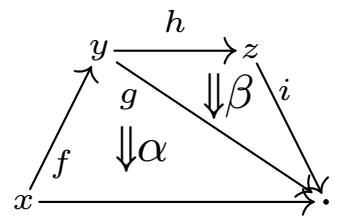

contains the following typings: $x, y, z : \varnothing$, $f : x \multimap \varnothing$, $g : y \multimap \varnothing$, $h : y \multimap \varnothing$, $i : z \multimap \varnothing$, $\alpha : g(y \leftarrow f) \multimap a \multimap \varnothing$, and $\beta : i(z \leftarrow h) \multimap b \multimap \varnothing$. Then for $t := \alpha(g \leftarrow \beta)$, we have

$$
\begin{aligned}
&\&_t \alpha = [], \\
&\&_t \beta = [\&_{\mathsf{s}\alpha} g] \cdot \&_\beta \beta = [[]] \cdot [] = [[]], \\
&\&_t i = [\&_{\mathsf{s}\alpha} g] \cdot \&_\beta i = [[]] \cdot [\&_{\mathsf{s}\beta} i] = [[]] \cdot [[]] = [[][]], \\
&\&_t h = [\&_{\mathsf{s}\alpha} g] \cdot \&_\beta h = [[]] \cdot [\&_{\mathsf{s}\beta} h] = [[]] \cdot [[\&_{\mathsf{s}i} z] \cdot \&_h h] = [[][[]]] = [[][*]], \\
&\&_t f = [\&_{\mathsf{s}\alpha} f] = [[\&_{\mathsf{s}g} y] \cdot \&_f f] = [[*]].
\end{aligned}
$$

Those addresses indeed match with that of the corresponding opetope:

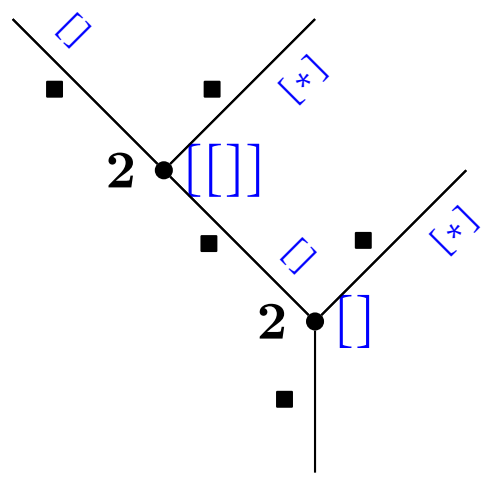

**Definition 3.2.3** (Polynomial coding)**.** The polynomial coding operation $[\![-]\!]^{\mathrm{poly}}_n$ is defined inductively by:

$$[\![x : \varnothing]\!]^{\mathrm{poly}}_0 := \blacklozenge \qquad x \in \mathbb{V}_0, \tag{3.2.4}$$

$$[\![\underline{x} : x \multimap \cdots]\!]^{\mathrm{poly}}_{n+2} := \mathsf{I}_{[\![x]\!]^{\mathrm{poly}}_n} \qquad x \in \mathbb{V}_n, \tag{3.2.5}$$

$$[\![x(\overrightarrow{y_i \leftarrow u_i}) : \cdots]\!]^{\mathrm{poly}}_{n+1} := \mathsf{Y}_{[\![x]\!]^{\mathrm{poly}}_n} \bigcirc_{[\&_{\mathsf{s}x} y_i]} [\![u_i]\!]^{\mathrm{poly}}_{n+1} \qquad x(\overrightarrow{y_i \leftarrow u_i}) \in \mathbb{T}_n. \tag{3.2.6}$$

We recall that (see equation (2.1.6) on page 10)

$$\mathsf{Y}_{[\![x]\!]^{\mathrm{poly}}_n} \bigcirc_{[\&_{\mathsf{s}x} y_i]} [\![u_i]\!]^{\mathrm{poly}}_{n+1} := \left( \cdots \left( \mathsf{Y}_{[\![x]\!]^{\mathrm{poly}}_n} \underset{[\&_{\mathsf{s}x} y_1]}{\circ} [\![u_1]\!]^{\mathrm{poly}}_{n+1} \right) \underset{[\&_{\mathsf{s}x} y_2]}{\circ} [\![u_2]\!]^{\mathrm{poly}}_{n+1} \cdots \right) \underset{[\&_{\mathsf{s}x} y_k]}{\circ} [\![u_k]\!]^{\mathrm{poly}}_{n+1}.$$

It is clear that the coding function is well defined in equations (3.2.4) and (3.2.5). We now establish a series of results to prove proposition 3.2.12 on page 23 stating that equation (3.2.6) is well defined too.

**Proposition 3.2.7.** *For $t \in \mathbb{T}_{\Gamma,n}$, $t = x(\overrightarrow{y_i \leftarrow u_i})$, $z \in t^\bullet$, and $[p] := \&_t z$, we have* $\mathsf{s}_{[p]} [\![t]\!]^{\mathrm{poly}}_{n+1} = [\![z]\!]^{\mathrm{poly}}_n$.

*Proof.* By definition we have

$$[\![t]\!]^{\mathrm{poly}}_{n+1} = \mathsf{Y}_{[\![x]\!]^{\mathrm{poly}}_n} \bigcirc_{[\&_{\mathsf{s}x} y_i]} [\![u_i]\!]^{\mathrm{poly}}_{n+1},$$

and we distinguish two cases. If $z = x$, then $[p] = []$, and the result clearly holds. Otherwise, $[p] = [\&_{\mathsf{s}x} y_j] \cdot \&_{u_j} z$, for $j$ such that $z \in u_j^\bullet$. Then,

$$
\begin{aligned}
\mathsf{s}_{[p]} [\![t]\!]^{\mathrm{poly}}_{n+1} &= \mathsf{s}_{[\&_{\mathsf{s}x} y_j] \cdot \&_{u_i} z} \left( \mathsf{Y}_{[\![x]\!]^{\mathrm{poly}}_n} \bigcirc_{[\&_{\mathsf{s}x} y_i]} [\![u_i]\!]^{\mathrm{poly}}_{n+1} \right) \\
&= \mathsf{s}_{\&_{u_j} z} [\![u_j]\!]^{\mathrm{poly}}_{n+1} \\
&= [\![z]\!]^{\mathrm{poly}}_n && \text{by induction}
\end{aligned}
$$

□

**Corollary 3.2.8.** *For $x$ as in case* (3.2.6) *above, and $[q_i] := \&_{\mathsf{s}x} y_i$ we have* $\mathsf{s}_{[q_i]} [\![x]\!]^{\mathrm{poly}}_n = [\![y_i]\!]^{\mathrm{poly}}_{n-1}$.

*Proof.* We have $\mathsf{s}_{[q_i]} [\![x]\!]^{\mathrm{poly}}_n = \mathsf{s}_{[q_i]} [\![\mathsf{s}\, x]\!]^{\mathrm{poly}}_n = [\![y_i]\!]^{\mathrm{poly}}_{n-1}$. □

**Lemma 3.2.9.** *Consider the following derivation tree:*

$$\dfrac{\cdots \vdash_n t : T \qquad \cdots \vdash_n x : X}{\cdots \vdash_n t(a \leftarrow x) : R}\ \texttt{graft-}a$$

*Writing* $r \coloneqq t(a \leftarrow x)$, *we have* $[\![r]\!]^{\mathrm{poly}}_{n+1} = [\![t]\!]^{\mathrm{poly}}_{n+1} \circ_{\&_t a} \mathsf{Y}_{[\![x]\!]^{\mathrm{poly}}_n}$.

*Proof.* We have $[\![r]\!]^{\mathrm{poly}}_{n+1} = [\![t]\!]^{\mathrm{poly}}_{n+1} \circ_{[l]} \mathsf{Y}_{[\![x]\!]^{\mathrm{poly}}_n}$, for some address $[l]$. Then, writing $t = z(\overrightarrow{y_i \leftarrow u_i})$, we have, according to definition definition 3.1.1 on page 17,

$$r = z(\overrightarrow{y_i \leftarrow u_i})(a \leftarrow x) = \begin{cases} z(\overrightarrow{y_i \leftarrow u_i}, a \leftarrow x) & \text{if } a \in (\mathsf{s}\, z)^\bullet, \\ z(\overrightarrow{y_i \leftarrow u_i(a \leftarrow x)}) & \text{if } a \notin (\mathsf{s}\, z)^\bullet. \end{cases}$$

There are two cases:

(1) if $a \in (\mathsf{s}\, z)^\bullet$, then $[l] = [\&_{\mathsf{s}\, z} a] = \&_t a$;

(2) if $a \notin (\mathsf{s}\, z)^\bullet$, then

$$[\![r]\!]^{\mathrm{poly}}_{n+1} = \mathsf{Y}_{[\![z]\!]^{\mathrm{poly}}_n} \bigcirc_{[\&_{\mathsf{s}\, z} y_i]} [\![u_i(a \leftarrow x)]\!]^{\mathrm{poly}}_{n+1}.$$

Let $j$ be the index such that $a \in (\mathsf{s}\, u_j)^\bullet$. Then, by induction,

$$[\![u_j(a \leftarrow x)]\!]^{\mathrm{poly}}_{n+1} = [\![u_j]\!]^{\mathrm{poly}}_{n+1} \underset{\&_{u_j} a}{\circ} \mathsf{Y}_{[\![x]\!]^{\mathrm{poly}}_n},$$

thus $[l] = [\&_{\mathsf{s}\, z} y_j] \cdot \&_{u_j} a = \&_t a$. □

**Lemma 3.2.10** (Named readdressing lemma)**.** *Let* $r \in \mathbb{T}_{\Gamma, n}$. *Then for* $b \in (\mathsf{s}\, r)^\bullet$, *we have* $\&_{\mathsf{s}\, r} b = \wp_{[\![r]\!]^{\mathrm{poly}}_{n+1}} \&_r b$, *where* $\wp$ *is the readdressing function introduced in section 2.1.1 on page 10.*

*Proof.* Recall the definition of the readdressing map from appendix A on page 71. If $r$ is a variable, then

$$\&_{\mathsf{s}\, r} b = \wp_{\mathsf{Y}_{[\![r]\!]^{\mathrm{poly}}_n}} [\&_{\mathsf{s}\, r} b] = \wp_{[\![r]\!]^{\mathrm{poly}}_{n+1}} \&_r b.$$

Otherwise, write the derivation tree of $r$:

$$\dfrac{\cdots \vdash_n t : T \qquad \cdots \vdash_n x : X}{\cdots \vdash_n r : R}\ \texttt{graft-}a$$

assume that $\&_{\mathsf{s}\, t} a = \wp_{[\![t]\!]^{\mathrm{poly}}_{n+1}} \&_t a$ holds by induction. Since $\mathsf{s}\, r = \mathsf{s}\, t[\mathsf{s}\, x / a]$, we have

$$\begin{aligned} \&_{\mathsf{s}\, r} b &= \begin{cases} \&_{\mathsf{s}\, t} a \cdot \&_{\mathsf{s}\, x} b & \text{if } b \in (\mathsf{s}\, x)^\bullet, \\ \&_{\mathsf{s}\, t} a \cdot \&_{\mathsf{s}\, x} c \cdot [p] & \text{if } b \in (\mathsf{s}\, t)^\bullet, \&_{\mathsf{s}\, t} a \sqsubseteq \&_{\mathsf{s}\, t} b, \\ & \text{say } \&_{\mathsf{s}\, t} b = \&_{\mathsf{s}\, t} a \cdot [\&_{\mathsf{s}\, a} c] \cdot [p], \\ \&_{\mathsf{s}\, t} b & \text{if } b \in (\mathsf{s}\, t)^\bullet, \&_{\mathsf{s}\, t} a \not\sqsubseteq \&_{\mathsf{s}\, t} b. \end{cases} \\ &= \wp_{\left([\![t]\!]^{\mathrm{poly}}_{n+1} \circ_{\&_t a} \mathsf{Y}_{[\![x]\!]^{\mathrm{poly}}_n}\right)} \&_r b = \wp_{[\![r]\!]^{\mathrm{poly}}_{n+1}} \&_r b. \end{aligned}$$

□

**Proposition 3.2.11.** *Let* $\alpha \in \mathbb{V}_n$ *be a typed* $n$*-variable, for* $n \geq 2$. *Then* $\mathsf{t}\, [\![\alpha]\!]^{\mathrm{poly}}_n = [\![\mathsf{s}\,\mathsf{s}\, \alpha]\!]^{\mathrm{poly}}_{n-1}$.

*Proof.* We proceed by induction on $\mathsf{s}\, \alpha$.

(1) If $\mathsf{s}\, \alpha = x \in \mathbb{V}_{n-1}$, then $[\![\alpha]\!]^{\mathrm{poly}}_n = \mathsf{Y}_{[\![x]\!]^{\mathrm{poly}}_{n-1}}$, and $\mathsf{t}\, [\![\alpha]\!]^{\mathrm{poly}}_n = [\![x]\!]^{\mathrm{poly}}_{n-1} = [\![\mathsf{s}\, x]\!]^{\mathrm{poly}}_{n-1} = [\![\mathsf{s}\,\mathsf{s}\, \alpha]\!]^{\mathrm{poly}}_{n-1}$.

(2) If $\mathsf{s}\, \alpha = \underline{a}$ for some $a \in \mathbb{V}_{n-2}$, then $[\![\alpha]\!]^{\mathrm{poly}}_n = \mathsf{I}_{[\![a]\!]^{\mathrm{poly}}_{n-2}}$, and $\mathsf{t}\, [\![\alpha]\!]^{\mathrm{poly}}_n = \mathsf{Y}_{[\![a]\!]^{\mathrm{poly}}_{n-2}} = [\![a]\!]^{\mathrm{poly}}_{n-1} = [\![\mathsf{s}\,\mathsf{s}\, \alpha]\!]^{\mathrm{poly}}_n$.

(3) Otherwise, $\mathsf{s}\, \alpha$ is given by the following proof tree:

$$\dfrac{\cdots \vdash_n t : T \qquad \cdots \vdash_n x : X}{\cdots \vdash_n \mathsf{s}\, \alpha : \mathsf{s}\,\mathsf{s}\, \alpha \multimapdotbothA \cdots}\ \texttt{graft-}a$$

and write $r := \mathsf{s}\,\alpha = t(a \leftarrow x)$ for shorter notations. By lemma 3.2.9 on the facing page we have

$$[\![r]\!]^{\text{poly}}_{n+1} = [\![t]\!]^{\text{poly}}_{n+1} \underset{\&_t a}{\circ} \mathsf{Y}_{[\![x]\!]^{\text{poly}}_n},$$

and finally, for $\square$ the partial multiplication of $\mathfrak{Z}^n$ as in theorem A.1.9 on page 72,

$$\begin{aligned}
[\![\mathsf{s}\,\mathsf{s}\,\alpha]\!]^{\text{poly}}_n &= [\![\mathsf{s}\,r]\!]^{\text{poly}}_n && \text{by definition} \\
&= [\![\mathsf{s}\,t]\!]^{\text{poly}}_n \underset{\&_{\mathsf{s}\,t} a}{\square} [\![\mathsf{s}\,x]\!]^{\text{poly}}_n \\
&= \mathsf{t}\,[\![t]\!]^{\text{poly}}_{n+1} \underset{\&_{\mathsf{s}\,t} a}{\square} [\![x]\!]^{\text{poly}}_n && \text{by induction} \\
&= \mathsf{t}\,[\![t]\!]^{\text{poly}}_{n+1} \underset{[p]}{\square} [\![x]\!]^{\text{poly}}_n && \text{with } [p] = \wp_{[\![t]\!]^{\text{poly}}_{n+1}}(\&_t a) \\
&= \mathsf{t}\left([\![t]\!]^{\text{poly}}_{n+1} \underset{\&_t a}{\circ} \mathsf{Y}_{[\![x]\!]^{\text{poly}}_n}\right) && \text{by lemma 3.2.10} \\
&= \mathsf{t}\,[\![\alpha]\!]^{\text{poly}}_{n+1}.
\end{aligned}$$

□

**Proposition 3.2.12.** *With variables as in equation* (3.2.6) *on page 21, we have that for all $i$*

$$\mathsf{t}\,\mathsf{s}_{[]}\,[\![u_i]\!]^{\text{poly}}_{n+1} = \mathsf{s}_{\&_{\mathsf{s}\,x} y_i}\,[\![x]\!]^{\text{poly}}_n,$$

*and the graftings are well defined.*

*Proof.* Write $u_i := a(\overrightarrow{b_j \leftarrow v_j})$. Then

$$\begin{aligned}
\mathsf{t}\,\mathsf{s}_{[]}\,[\![u_i]\!]^{\text{poly}}_{n+1} &= \mathsf{t}\,[\![a]\!]^{\text{poly}}_n && \text{by proposition 3.2.7} \\
&= [\![\mathsf{s}\,\mathsf{s}\,a]\!]^{\text{poly}}_{n-1} && \text{by proposition 3.2.11} \\
&= [\![\mathsf{s}\,y_i]\!]^{\text{poly}}_{n-1} && \text{by } \texttt{graft} \text{ rule} \\
&= [\![y_i]\!]^{\text{poly}}_{n-1} && \text{by definition} \\
&= \mathsf{s}_{\&_{\mathsf{s}\,x} y_i}\,[\![x]\!]^{\text{poly}}_n && \text{by corollary 3.2.8.}
\end{aligned}$$

□

This result concludes the proof that equations (3.2.4) to (3.2.6) of definition 3.2.3 on page 21 are well defined.

**Corollary 3.2.13.** *Let $(E \triangleright \Gamma \vdash_n t : T)$ be a derivable sequent. Then $\&_t$ exhibits a bijection between the set of $n$-variables of $t$, and $\left([\![t]\!]^{\text{poly}}_{n+1}\right)^{\bullet}$. It is also a bijection between the set of $(n-1)$-variables of $\mathsf{s}\,t$ and $\left([\![t]\!]^{\text{poly}}_{n+1}\right)^{|}$.*

The rest of this section is dedicated to prove theorem 3.2.22 on page 25 stating that $[\![-]\!]^{\text{poly}}_n$ is a bijection modulo $\simeq$. We first prove surjectivity, by defining a sequent $[\![\omega]\!]^!$ such that $\left[\!\left[[\![\omega]\!]^!\right]\!\right]^{\text{poly}}_n = \omega$, for any opetopes $\omega \in \mathbb{O}_n$. We proceed by opetopic induction (see remark 2.2.3 on page 12).

(1) Trivially, $[\![\blacklozenge]\!]^!$ is obtained by the following proof tree:

$$\frac{}{[\![\blacklozenge]\!]^!}\ \texttt{point} \tag{3.2.14}$$

with an arbitrary choice of variable (different choices lead to equivalent sequents).

(2) For $\phi \in \mathbb{O}_{n-2}$ the sequent $[\![\mathsf{I}_\phi]\!]^!$ is obtained by the following proof tree:

$$\frac{[\![\phi]\!]^!}{[\![\mathsf{I}_\phi]\!]^!}\ \texttt{degen} \tag{3.2.15}$$

(3) For $\psi \in \mathbb{O}_{n-1}$, the sequent $[\![\mathsf{Y}_\psi]\!]^!$ is obtained by the following proof tree:

$$\frac{[\![\psi]\!]^!}{[\![\mathsf{Y}_\psi]\!]^!}\ \texttt{shift} \tag{3.2.16}$$

with an arbitrary choice of fresh variable (different choices lead to equivalent sequents).

(4) Let $\nu \in \mathbb{O}_n$ having at least one node, $[l] \in \nu^|$, and $\psi \in \mathbb{O}_{n-1}$ be such that the grafting $\nu \circ_{[l]} \mathsf{Y}_\psi$ is well-defined. Then the sequent $[\![\nu \circ_{[l]} \mathsf{Y}_\psi]\!]^!$ is obtained by the following proof tree:

$$\frac{[\![\nu]\!]^! \qquad [\![\psi]\!]^!}{[\![\nu \circ_{[l]} \mathsf{Y}_\psi]\!]^!}\ \mathtt{graft}\text{-}a \tag{3.2.17}$$

where $[\![\nu]\!]^! = (E \triangleright \Gamma \vdash_n u : U)$, where the variable $a \in (\mathsf{s}\,\mathsf{s}\,u)^\bullet$ is such that $\&_{\mathsf{s}\mathsf{s}u} a = [l]$ (see corollary 3.2.13 on the preceding page), and where the adequate $\alpha$-conversion have been performed to fulfill the side conditions of $\mathtt{graft}$. We check in propositions 3.2.19 and 3.2.20 that this definition is well-founded.

**Lemma 3.2.18.** *Let* $(E \triangleright \Gamma \vdash_n w : W) := [\![\omega]\!]^!$, *with* $\omega \in \mathbb{O}_n$ *non degenerate,* $n \geq 2$. *Then, by corollary 3.2.13 on the preceding page, for* $[l] \in \omega^|$, *there is a variable* $a \in (\mathsf{s}\,\mathsf{s}\,w)^\bullet$ *such that* $\&_{\mathsf{s}\mathsf{s}w} a = [l]$.

*Proof.* By assumption, $\omega$ is either of the form $\mathsf{Y}_\psi$, or $\nu \circ_{[k]} \mathsf{Y}_\psi$, for some $\nu \in \mathbb{O}_n$ and $\psi \in \mathbb{O}_{n-1}$. From there, this is a straightforward induction. ☐

**Proposition 3.2.19.** *In proof tree (3.2.6), the instance of* $\mathtt{graft}$ *is well-defined.*

*Proof.* Write $[\![\psi]\!]^! = (F \triangleright \Upsilon \vdash_{n-1} p : P)$. By the previous lemma, we have $a \in (\mathsf{s}\,\mathsf{s}\,u)^\bullet$. Further,

$$\begin{aligned}
[\![\mathsf{s}\,a]\!]^{\text{poly}}_{n-2} &= [\![\mathsf{e}_{[l]}\,u]\!]^{\text{poly}}_{n-2} \\
&= \mathsf{e}_{[l]}\,[\![u]\!]^{\text{poly}}_{n} && \text{by proposition 3.2.7 on page 21} \\
&= \mathsf{e}_{[l]}\,\nu \\
&= \mathsf{t}\,\mathsf{s}_{[]}\,\mathsf{Y}_\psi && \text{by relation } \textbf{(Inner)}\text{, see section 2.2.4} \\
&= \mathsf{t}\,\psi \\
&= [\![\mathsf{s}\,\mathsf{s}\,p]\!]^{\text{poly}}_{n-2} && \text{by proposition 3.2.11 on page 22.}
\end{aligned}$$

By induction on $n$, the polynomial coding $[\![-]\!]^{\text{poly}}_{n-2}$ is injective modulo $\simeq$. Hence, we can assume $\mathsf{s}\,a = \mathsf{s}\,\mathsf{s}\,p$ without loss of generality, and finally, the instance of the $\mathtt{graft}$ rule is well-defined. ☐

**Proposition 3.2.20.** *Let* $\omega \in \mathbb{O}_n$ *be a non degenerate opetope, and consider two arbitrary decompositions in corollas*

$$\begin{aligned}
\omega &= \left(\cdots\left(\mathsf{Y}_{\mathsf{s}_{[p_1]}\omega} \underset{[p_2]}{\circ} \mathsf{Y}_{\mathsf{s}_{[p_2]}\omega}\right) \underset{[p_3]}{\circ} \mathsf{Y}_{\mathsf{s}_{[p_3]}\omega} \cdots\right) \underset{[p_k]}{\circ} \mathsf{Y}_{\mathsf{s}_{[p_k]}\omega} \\
&= \left(\cdots\left(\mathsf{Y}_{\mathsf{s}_{[q_1]}\omega} \underset{[q_2]}{\circ} \mathsf{Y}_{\mathsf{s}_{[q_2]}\omega}\right) \underset{[q_3]}{\circ} \mathsf{Y}_{\mathsf{s}_{[q_3]}\omega} \cdots\right) \underset{[q_k]}{\circ} \mathsf{Y}_{\mathsf{s}_{[q_k]}\omega}.
\end{aligned}$$

*Then*

$$\left[\!\!\left[\left(\mathsf{Y}_{\mathsf{s}_{[p_1]}\omega} \underset{[p_2]}{\circ} \mathsf{Y}_{\mathsf{s}_{[p_2]}\omega}\right) \cdots \underset{[p_k]}{\circ} \mathsf{Y}_{\mathsf{s}_{[p_k]}\omega}\right]\!\!\right]^! = \left[\!\!\left[\left(\mathsf{Y}_{\mathsf{s}_{[q_1]}\omega} \underset{[q_2]}{\circ} \mathsf{Y}_{\mathsf{s}_{[q_2]}\omega}\right) \cdots \underset{[q_k]}{\circ} \mathsf{Y}_{\mathsf{s}_{[q_k]}\omega}\right]\!\!\right]^!$$

*In other words,* $[\![\omega]\!]^!$ *does not depend on the decomposition of* $\omega$ *in corollas.*

*Proof.* By assumption, the sequence $[p_1], \ldots, [p_k]$ (and likewise for $[q_1], \ldots, [q_k]$) has the following property: for $1 \leq i \leq j \leq k$, either $[p_i] \leq [p_j]$ or $[p_i]$ and $[p_j]$ are $\leq$-incomparable (recall that $\leq$ is the lexicographical order on $\mathbb{A}_{n-1}$, see section 2.2.3 on page 11). Further, $\{[p_1], \ldots, [p_k]\} = \omega^\bullet = \{[q_1], \ldots, [q_k]\}$. Consequently, the sequence of addresses $[q_1], \ldots, [q_k]$ can be obtained from $[p_1], \ldots, [p_k]$ via a sequence of transpositions of consecutive $\leq$-incomparable addresses.

It is thus enough to check the following: for $\nu \in \mathbb{O}_n$, $[l], [l'] \in \nu^|$ (necessarily, neither is a prefix of the other), and $\psi, \psi' \in \mathbb{O}_{n-1}$ such that the following graftings are well defined, we have

$$\left[\!\!\left[\left(\nu \underset{[l]}{\circ} \mathsf{Y}_\psi\right) \underset{[l']}{\circ} \mathsf{Y}_{\psi'}\right]\!\!\right]^! = \left[\!\!\left[\left(\nu \underset{[l']}{\circ} \mathsf{Y}_{\psi'}\right) \underset{[l]}{\circ} \mathsf{Y}_\psi\right]\!\!\right]^!.$$

Let $[\![\nu]\!]^! = (E_\nu \triangleright \Gamma_\nu \vdash x_\nu : s_\nu \multimap X_\nu)$, and likewise for $\psi$ and $\psi'$, and $a, a' \in (\mathsf{s}\mathsf{s}\nu)^\bullet$ be such that $\&_{\mathsf{ss}\nu} a = [l]$ and $\&_{\mathsf{ss}\nu} a' = [l']$ (see corollary 3.2.13 on page 23). The sequents above are respectively obtained by the following proof trees:

$$\dfrac{\dfrac{E_\nu \triangleright \Gamma_\nu \vdash x_\nu : s_\nu \multimap X_\nu \qquad E_\psi \triangleright \Gamma_\psi \vdash x_\psi : s_\psi \multimap X_\psi}{F \triangleright \Gamma_\nu \cup \Gamma_\psi \vdash x_\nu(a \leftarrow x_\psi) : s_\nu[s_\psi/a] \multimap X_\nu}\,\texttt{graft-}a \qquad E_{\psi'} \triangleright \Gamma_{\psi'} \vdash x_{\psi'} : s_{\psi'} \multimap X_{\psi'}}{G \triangleright \Gamma_\nu \cup \Gamma_\psi \cup \Gamma_{\psi'} \vdash x_\nu(a \leftarrow x_\psi)(a' \leftarrow x_{\psi'}) : s_\nu[s_\psi/a][s_{\psi'}/a'] \multimap X_\nu}\,\texttt{graft-}a'$$

$$\dfrac{\dfrac{E_\nu \triangleright \Gamma_\nu \vdash x_\nu : s_\nu \multimap X_\nu \qquad E_{\psi'} \triangleright \Gamma_{\psi'} \vdash x_{\psi'} : s_{\psi'} \multimap X_{\psi'}}{F' \triangleright \Gamma_\nu \cup \Gamma_{\psi'} \vdash x_\nu(a' \leftarrow x_{\psi'}) : s_\nu[s_{\psi'}/a'] \multimap X_\nu}\,\texttt{graft-}a' \qquad E_{\psi} \triangleright \Gamma_{\psi} \vdash x_{\psi} : s_{\psi} \multimap X_{\psi}}{G' \triangleright \Gamma_\nu \cup \Gamma_{\psi'} \cup \Gamma_{\psi} \vdash x_\nu(a' \leftarrow x_{\psi'})(a \leftarrow x_\psi) : s_\nu[s_{\psi'}/a'][s_\psi/a] \multimap X_\nu}\,\texttt{graft-}a$$

It remains to prove that both those conclusion sequents are equal.

(1) Since by assumptions $[l], [l'] \in \nu^|$, we have $a \notin s_{\psi'}^\bullet$ and $a' \notin s_\psi^\bullet$. Thus
$$x_\nu(a \leftarrow x_\psi)(a' \leftarrow x_{\psi'}) = x_\nu(a \leftarrow x_\psi, a' \leftarrow x_{\psi'}) = x_\nu(a' \leftarrow x_{\psi'})(a \leftarrow x_\psi).$$
(2) Again, since $a \notin s_{\psi'}^\bullet$ and $a' \notin s_\psi^\bullet$, we have
$$s_\nu[s_\psi/a][s_{\psi'}/a'] = s_\nu[s_{\psi'}/a'][s_\psi/a].$$
(3) Lastly, the equational theories $G$ and $G'$ are the union of $E_\nu$, $E_\psi$, and $E_{\psi'}$, and the potential additional equalities incurred by the independent substitutions $s_\psi/a$ and $s_{\psi'}/a'$. Hence $G = G'$. □

**Corollary 3.2.21.** *For any opetope $\omega \in \mathbb{O}$, the sequent $[\![\omega]\!]^!$ is uniquely defined up to $\simeq$.*

*Proof.* Clearly, proof trees (3.2.14), (3.2.15), and (3.2.16) are well-defined. In proposition 3.2.19 on the preceding page, we show that the same holds for proof tree (3.2.17). Finally, in proposition 3.2.20 on the facing page, we show that for a non degenerate opetope $\omega \in \mathbb{O}_n$, the sequent $[\![\omega]\!]^!$ does not depend on the decomposition of $\omega$. □

**Theorem 3.2.22.** *The polynomial coding $[\![-]\!]_n^{\mathrm{poly}}$ is a bijection modulo $\simeq$, whose inverse is $[\![-]\!]^!$ restricted to $\mathbb{O}_n$.*

*Proof.* We first show that for $\omega \in \mathbb{O}_n$ we have $\left[\!\!\left[[\![\omega]\!]^!\right]\!\!\right]_n^{\mathrm{poly}} = \omega$ by opetopic induction (see remark 2.2.3 on page 12).

(1) By definition of $[\![-]\!]^{\mathrm{poly}}$, $\left[\!\!\left[[\![\blacklozenge]\!]^!\right]\!\!\right]_0^{\mathrm{poly}} = \blacklozenge$.
(2) With the same notations as in (3.2.15), and by induction, we have
$$\left[\!\!\left[[\![\mathsf{I}_\phi]\!]^!\right]\!\!\right]_n^{\mathrm{poly}} = \mathsf{I}_{[\![[\![\phi]\!]^!]\!]_{n-2}^{\mathrm{poly}}} = \mathsf{I}_\phi.$$
(3) With the same notations as in (3.2.16), and by induction, we have,
$$\left[\!\!\left[[\![\mathsf{Y}_\psi]\!]^!\right]\!\!\right]_n^{\mathrm{poly}} = \mathsf{Y}_{[\![[\![\psi]\!]^!]\!]_{n-1}^{\mathrm{poly}}} = \mathsf{Y}_\psi.$$
(4) With the same notations as in (3.2.17), and by induction, we have
$$\left[\!\!\left[\left[\!\!\left[\nu \underset{[l]}{\circ} \mathsf{Y}_\psi\right]\!\!\right]^!\right]\!\!\right]_n^{\mathrm{poly}} = \left[\!\!\left[[\![\nu]\!]^!\right]\!\!\right]_n^{\mathrm{poly}} \underset{[\&_{\mathsf{ss}u} a]}{\circ} \left[\!\!\left[[\![\mathsf{Y}_\psi]\!]^!\right]\!\!\right]_n^{\mathrm{poly}} = \left[\!\!\left[[\![\nu]\!]^!\right]\!\!\right]_n^{\mathrm{poly}} \underset{[l]}{\circ} \left[\!\!\left[[\![\mathsf{Y}_\psi]\!]^!\right]\!\!\right]_n^{\mathrm{poly}} = \nu \underset{[l]}{\circ} \mathsf{Y}_\psi.$$

Conversely, we now show that for all sequents $(E \triangleright \Gamma \vdash \alpha : T)$ (abbreviated $(\alpha : T)$ if no ambiguity arise) we have $(E \triangleright \Gamma \vdash \alpha : T) = \left[\!\!\left[[\![E \triangleright \Gamma \vdash \alpha : T]\!]_n^{\mathrm{poly}}\right]\!\!\right]^!$ up to $\simeq$.

(1) We have that $\left[\!\!\left[[\![x : \varnothing]\!]_0^{\mathrm{poly}}\right]\!\!\right]^! = [\![\blacklozenge]\!]^! \simeq (x : \varnothing)$.
(2) With the same notations as in equation (3.2.5) on page 21, $\left[\!\!\left[[\![\delta : \underline{x} \multimap x \multimap X]\!]_n^{\mathrm{poly}}\right]\!\!\right]^! = \left[\!\!\left[\mathsf{I}_{[\![x:X]\!]_n^{\mathrm{poly}}}\right]\!\!\right]^!$, and both sequents $\left[\!\!\left[\mathsf{I}_{[\![x:X]\!]_n^{\mathrm{poly}}}\right]\!\!\right]^!$ and $(\delta : \underline{x} \multimap x \multimap X)$ are obtained by applying `degen` to $(x : X)$. Thus $\left[\!\!\left[[\![\delta : \underline{x} \multimap x \multimap X]\!]_n^{\mathrm{poly}}\right]\!\!\right]^! \simeq (\delta : \underline{x} \multimap x \multimap X)$.

(3) Lastly, consider the sequent $\left(\alpha : x(\overrightarrow{y_i \leftarrow u_i}) \multimap T\right)$ as in equation (3.2.6) on page 21. Then

$$\left[\!\left[ [\![\alpha : x(\overrightarrow{y_i \leftarrow u_i}) \multimap T]\!]^{\text{poly}}_{n+1} \right]\!\right]^! = \left[\!\left[ \mathsf{Y}_{[\![x]\!]^{\text{poly}}_n} \bigcirc_{[\&_{\mathsf{s}x} y_i]} [\![u_i]\!]^{\text{poly}}_{n+1} \right]\!\right]^! \simeq \left(\alpha : x(\overrightarrow{y_i \leftarrow u_i}) \multimap T'\right).$$

Since $T$ and $T'$ are completely determined by $x(\overrightarrow{y_i \leftarrow u_i})$ (see theorem 3.1.8 on page 19), we have that $T = T'$, whence

$$\left[\!\left[ [\![\alpha : x(\overrightarrow{y_i \leftarrow u_i}) \multimap T]\!]^{\text{poly}}_{n+1} \right]\!\right]^! \simeq \left(\alpha : x(\overrightarrow{y_i \leftarrow u_i}) \multimap T\right).$$

□

## 3.3. Examples

In this section, we showcase the derivation of some low dimensional opetopes. On a scale of a proof tree, specifying the context at every step is redundant. Hence we allow omitting it, only having the equational theory on the left of $\vdash$.

**Example 3.3.1** (The arrow)**.** The unique 1-opetope, the *arrow*, is given by the following simple derivation:

$$\dfrac{\dfrac{}{\vdash_0 a : \varnothing}\,\texttt{point}}{\vdash_1 f : a \multimap \varnothing}\,\texttt{shift}$$

**Example 3.3.2** (Opetopic integers)**.** The set of $\mathbb{O}_2$ of 2-opetopes is in bijection with the set of natural numbers. Given $n \in \mathbb{N}$, we denote by $\mathbf{n}$ the 2-opetope whose pasting diagram is a sequence of $n$ arrows as follows:

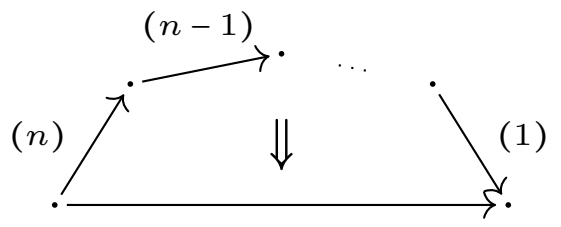

The derivation of the opetope $\mathbf{0}$ is

$$\dfrac{\dfrac{\dfrac{}{\vdash_0 a : \varnothing}\,\texttt{point}}{\vdash_1 \underline{a} : a \multimap \varnothing}\,\texttt{degen}}{\vdash_1 \mathbf{0} : \underline{a} \multimap a \multimap \varnothing}\,\texttt{shift}$$

alternatively, we could have used the `degen-shift` rule. For $n \geq 1$, the opetope $\mathbf{n}$ is derived as

$$\dfrac{\dfrac{\begin{array}{c}\vdots\\ \dfrac{\dfrac{\dfrac{\begin{array}{c}\vdots\\ \vdash_1 f_1 : a_1 \multimap \varnothing\end{array} \quad \begin{array}{c}\vdots\\ \vdash_1 f_2 : a_2 \multimap \varnothing\end{array}}{\vdash_1 f_1(a_1 \leftarrow f_2) : a_2 \multimap \varnothing}\,\texttt{graft-}a_1 \quad \begin{array}{c}\vdots\\ \vdash_1 f_3 : a_3 \multimap \varnothing\end{array}}{\vdash_1 f_1(a_1 \leftarrow f_2(a_2 \leftarrow f_3)) : a_3 \multimap \varnothing}\,\texttt{graft-}a_2}{}\\ \vdots\\ \vdash_1 f_1(a_1 \leftarrow f_2(\cdots a_{n-2} \leftarrow f_{n-1})) : a_{n-1} \multimap \varnothing\end{array} \quad \begin{array}{c}\vdots\\ \vdash_1 f_n : a_n \multimap \varnothing\end{array}}{\vdash_1 f_1(a_1 \leftarrow f_2(\cdots a_{n-1} \leftarrow f_n)) : a_n \multimap \varnothing}\,\texttt{graft-}a_{n-1}}{\vdash_1 \mathbf{n} : f_1(a_1 \leftarrow f_2(\cdots a_{n-1} \leftarrow f_n)) \multimap a_n \multimap \varnothing}\,\texttt{shift}$$

**Example 3.3.3** (A classic)**.** The 3-opetope

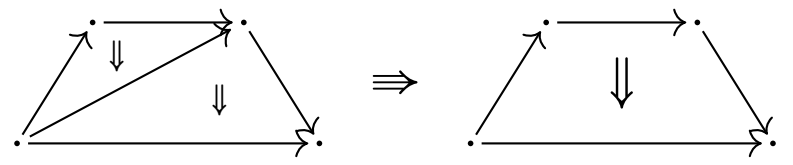

is derived as follows.

$$\dfrac{\dfrac{\dfrac{}{\vdash_0 c : \varnothing}\,\texttt{point}}{\vdash_1 h : c \multimap \varnothing}\,\texttt{shift} \quad \dfrac{\dfrac{}{\vdash_0 a : \varnothing}\,\texttt{point}}{\vdash_1 i : a \multimap \varnothing}\,\texttt{shift}}{\vdash_1 h(c \leftarrow i) : c[a/c] \multimap \varnothing}\,\texttt{graft-}c$$

and $c[a/c] = a$. Then,

$$\dfrac{\begin{array}{c}\vdots\\ \vdash_1 h(c \leftarrow i) : a \multimap \varnothing\end{array}}{\vdash_2 \beta : h(c \leftarrow i) \multimap a \multimap \varnothing}\,\texttt{shift}$$

On the other hand we have

$$\dfrac{\dfrac{\dfrac{}{\vdash_0 b : \varnothing}\,\texttt{point}}{\vdash_1 g : b \multimap \varnothing}\,\texttt{shift} \qquad \dfrac{\dfrac{}{\vdash_0 a : \varnothing}\,\texttt{point}}{\vdash_1 f : a \multimap \varnothing}\,\texttt{shift}}{\vdash_1 g(b \leftarrow f) : b[a/b] \multimap \varnothing}\,\texttt{graft-}b$$

and $b[a/b] = a$. Then,

$$\dfrac{\begin{array}{c}\vdots\\ \vdash_2 \beta : h(c \leftarrow i) \multimap a \multimap \varnothing\end{array} \qquad \dfrac{\begin{array}{c}\vdots\\ \vdash_1 g(b \leftarrow f) : \multimap \varnothing\end{array}}{\vdash_2 \alpha : g(b \leftarrow f) \multimap a \multimap \varnothing}\,\texttt{shift}}{\vdash_2 \beta(i \leftarrow \alpha) : h(c \leftarrow i)[g(b \leftarrow f)/i] \multimap a \multimap \varnothing}\,\texttt{graft-}i$$

The last grafting is well defined as $\mathsf{s}\,i = a = \mathsf{ss}\,\alpha$, and $h(c \leftarrow i)[g(b \leftarrow f)/i] = h(c \leftarrow g(b \leftarrow f))$. Finally

$$\dfrac{\begin{array}{c}\vdots\\ \vdash_2 \beta(i \leftarrow \alpha) : h(c \leftarrow g(b \leftarrow f)) \multimap a \multimap \varnothing\end{array}}{\vdash_3 A : \beta(i \leftarrow \alpha) \multimap h(c \leftarrow g(b \leftarrow f)) \multimap a \multimap \varnothing}\,\texttt{shift}$$

Finally, the complete proof tree is:

$$\dfrac{\dfrac{\dfrac{\dfrac{\dfrac{}{\vdash c : \varnothing}\,\texttt{point}}{\vdash h : c \multimap \varnothing}\,\texttt{shift} \quad \dfrac{\dfrac{}{\vdash a : \varnothing}\,\texttt{point}}{\vdash i : a \multimap \varnothing}\,\texttt{shift}}{\vdash h(c \leftarrow i) : a \multimap \varnothing}\,\texttt{graft-}c}{\vdash \beta : h(c \leftarrow i) \multimap a \multimap \varnothing}\,\texttt{shift} \quad \dfrac{\dfrac{\dfrac{\dfrac{}{\vdash b : \varnothing}\,\texttt{point}}{\vdash g : b \multimap \varnothing}\,\texttt{shift} \quad \dfrac{\dfrac{}{\vdash a : \varnothing}\,\texttt{point}}{\vdash f : a \multimap \varnothing}\,\texttt{shift}}{\vdash g(b \leftarrow f) : a \multimap \varnothing}\,\texttt{graft-}b}{\vdash \alpha : g(b \leftarrow f) \multimap a \multimap \varnothing}\,\texttt{shift}}{\dfrac{\vdash \beta(i \leftarrow \alpha) : h(c \leftarrow g(b \leftarrow f)) \multimap a \multimap \varnothing}{\vdash A : \beta(i \leftarrow \alpha) \multimap h(c \leftarrow g(b \leftarrow f)) \multimap a \multimap \varnothing}\,\texttt{shift}}\,\texttt{graft-}i$$

**Example 3.3.4** (A degenerate case)**.** The 3-opetope on the left is derived as on the right:

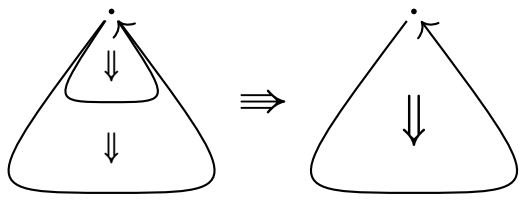

$$\dfrac{\dfrac{\dfrac{\dfrac{\dfrac{}{\vdash a : \varnothing}\,\texttt{point}}{\vdash f : a \multimap \varnothing}\,\texttt{shift}}{\vdash \alpha : f \multimap a \multimap \varnothing}\,\texttt{shift} \quad \dfrac{\dfrac{}{\vdash a : \varnothing}\,\texttt{point}}{\vdash \delta : \underline{a} \multimap a \multimap \varnothing}\,\texttt{degen}}{\vdash \alpha(f \leftarrow \delta) : \underline{a} \multimap a \multimap \varnothing}\,\texttt{graft-}f}{\vdash A : \alpha(f \leftarrow \delta) \multimap \underline{a} \multimap a \multimap \varnothing}\,\texttt{shift}$$

**Example 3.3.5** (Another degenerate case)**.** The 3-opetope

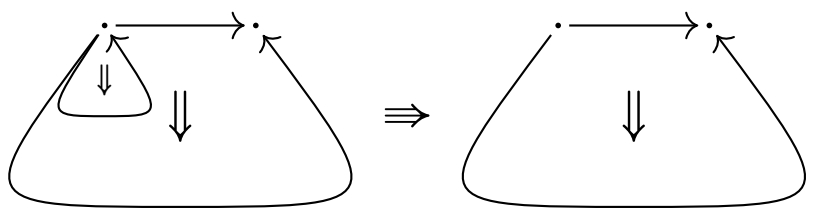

is derived as follows:

$$\dfrac{\dfrac{\dfrac{\dfrac{\dfrac{}{\vdash b : \varnothing}\,\texttt{point}}{\vdash g : b \multimap \varnothing}\,\texttt{shift} \quad \dfrac{\dfrac{}{\vdash a : \varnothing}\,\texttt{point}}{\vdash f : a \multimap \varnothing}\,\texttt{shift}}{\vdash g(b \leftarrow f) : a \multimap \varnothing}\,\texttt{graft-}b}{\vdash \beta : g(b \leftarrow f) \multimap a \multimap \varnothing}\,\texttt{shift} \quad \dfrac{\dfrac{\dfrac{}{\vdash a : \varnothing}\,\texttt{point}}{\vdash \underline{a} : a \multimap \varnothing}\,\texttt{degen}}{\vdash \alpha : \underline{a} \multimap a \multimap \varnothing}\,\texttt{shift}}{a = b \vdash \beta(f \leftarrow \alpha) : g(b \leftarrow f)[\underline{a}/f] \multimap a \multimap \varnothing}\,\texttt{graft-}f$$

and $g(b \leftarrow f)[\underline{a}/f] = g$, with the added equality $a = b$.

$$\dfrac{\begin{array}{c}\vdots \\ a = b \vdash \beta(f \leftarrow \alpha) : g \bullet\!\!-\!\!\circ\ a \bullet\!\!-\!\!\circ\ \varnothing\end{array}}{a = b \vdash A : \beta(f \leftarrow \alpha) \bullet\!\!-\!\!\circ\ g \bullet\!\!-\!\!\circ\ a \bullet\!\!-\!\!\circ\ \varnothing}\ \texttt{shift}$$

## 3.4. Python implementation

In this section, we briefly discuss the Python implementation [**Ho Thanh, 2018b**] of the present work. System $\mathrm{Opt}^{!}$ and all required syntactic constructs are implemented in module `opetopy.NamedOpetope`. The rules are represented by functions `point`, `degen`, `shift`, `graft`, as well as `degenshift` for the alternative rule presented in remark 3.1.7 on page 18. Those rules are further encapsulated in rule instance classes `Point`, `Degen`, `shift`, `Graft`, and `DegenShift`, which represent rule instances in a proof tree, so constructing a derivation amounts to writing a Python term using those four classes. If that term evaluates without raising any exception, then the proof tree is considered correct.

Figure 3.4.1. Derivation of the arrow sequent using `opetopy.NamedOpetope`

```python
from opetopy.NamedOpetope *
# We first derive the point that will act as the source of the arrow by invoking the point rule on variable "a".
a = Point("a")
# We then apply the shift rule on a by providing a fresh variable, here "f".
f = Shift(a, "f")
# Since we use names, the following sequent, while corresponding to the same opetope, is different from f
g = Shift(Point("b"), "g")
# Note that the function opetopy.NamedOpetope.Arrow can be used to concisely get a proof tree of ■.
```

Figure 3.4.2. Derivation of some opetopic integers using `opetopy.NamedOpetope`, continuation of figure 3.4.1

```python
opetopic_integer_0 = DegenShift(a, "n_0")
opetopic_integer_1 = Shift(f, "n_1")
opetopic_integer_2 = Shift(Graft(g, f, "b"), "n_2")
# Note that the function opetopy.NamedOpetope.OpetopicInteger can be used to get the proof tree of an arbitrary opetopic integer.
```

FIGURE 3.4.3. Derivation of example 3.3.3 on page 26 using `opetopy.NamedOpetope`

```
from opetopy.NamedOpetope import *
f = Shift(Point("a"), "f")
g = Shift(Point("b"), "g")
h = Shift(Point("c"), "h")
i = Shift(Point("a"), "i")
alpha = Shift(Graft(g, f, "b"), "alpha")
beta = Shift(Graft(h, i, "c"), "beta")
A = Shift(Graft(beta, alpha, "i"), "A")
```

FIGURE 3.4.4. Derivation of example 3.3.5 on page 27 using `opetopy.NamedOpetope`

```
from opetopy.NamedOpetope import *
f = Shift(Point("a"), "f")
g = Shift(Point("b"), "g")
alpha = DegenShift(Point("a"), "alpha")
beta = Shift(Graft(g, f, "b"), "beta")
D = Shift(Graft(beta, alpha, "f"), "D")
```

## 3.5. The system for opetopic sets

We now present $\text{OPTSET}^!$, a derivation system for opetopic sets that is based on $\text{OPT}^!$. We first present the required syntactic constructs and conventions in section 3.1.1 on page 15, and present the inference rules in figure 3.5.1 on the next page.

3.5.1. **Syntax.** The main distinguishing feature of the named approach above is that only source faces of whichever cell is currently being derived are specified:

$$\cdots \vdash_n x : \mathsf{s}x \multimap \mathsf{ss}x \multimap \cdots$$

Nonetheless, as proven in proposition 3.2.11 on page 22, all the information about targets remain. To adapt our previous derivation system to opetopic sets, all faces, including targets, need to be explicitly specified. This will be part of the role of rule `repr` of system $\text{OPTSET}^!$. Further, recall that a sequent in system $\text{OPT}^!$:

$$E \triangleright \Gamma \vdash t : T.$$

An opetopic set is a set of opetopic cells, and in a sequent, the context and the equational theory suffice to describe this. Thus, system $\text{OPTSET}^!$ will only deal with expressions of the form $(E \triangleright \Gamma)$, called *opetopic contexts modulo theory* (or OCMTs for short).

For example, the OCMT describing the following opetopic set (note that the graphical representation is not unique):

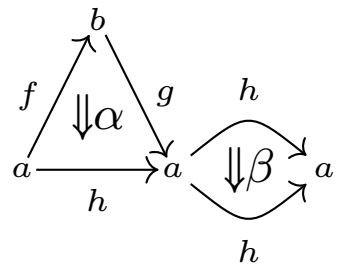

is given by

$$\left( \begin{array}{l} b = \mathsf{t}f, a = \mathsf{t}g = \mathsf{tt}\alpha = \mathsf{t}h = \mathsf{tt}\beta \\ h = \mathsf{t}\beta = \mathsf{t}\alpha \end{array} \triangleright \begin{array}{l} a : \varnothing, b : \varnothing, \mathsf{t}f : \varnothing, \mathsf{t}g : \varnothing, \mathsf{tt}\alpha : \varnothing, \mathsf{t}h : \varnothing, \mathsf{tt}\beta : \varnothing \\ f : a \multimap \varnothing, g : a \multimap \varnothing, \mathsf{t}\alpha : a \multimap \varnothing, h : a \multimap \varnothing, \mathsf{t}\beta : a \multimap \varnothing \\ \alpha : g(b \leftarrow f) \multimap a \multimap \varnothing, \beta : h \multimap a \multimap \varnothing \end{array} \right)$$

3.5.2. **Inference rules.** Our derivation system for opetopic sets, presented in figure 3.5.1, has four rules:

(1) `repr` that takes an opetope in our previous system and makes it into the representable opetopic set of that opetope;
(2) `zero` that constructs the empty OCMT;
(3) `sum` that takes the disjoint union of two opetopic sets;
(4) `glue` that identifies cells of an opetopic set.

In virtue of the fact that every finite opetopic set is a quotient of a finite sum of representables, those rules should be enough to derive all finite opetopic sets. This is proved in theorem 3.6.15 on page 35.

FIGURE 3.5.1. The OPTSET$^!$ system.

**Introduction of all targets:** This rule takes an opetope $(E \triangleright \Gamma \vdash_n x : X)$, and completes it by adding all the missing cells, i.e. targets, and making it into an OCMT. If $x \in \mathbb{V}_n$, then

$$\frac{E \triangleright \Gamma \vdash_n x : X}{E' \triangleright \Gamma'}\ \texttt{repr}$$

where

$$\Gamma' := \Gamma \cup \left\{ \mathsf{t}^k\, a : \mathsf{s}^{k+1}\, a \multimap \mathsf{s}^{k+2}\, a \multimap \cdots \mid a \in \mathbb{V}_{\Gamma,l}, 1 \leq k \leq l \leq n \right\},$$

and

$$\begin{aligned} E' := E & \\ \cup\ & \{ \mathsf{t}\, a = b \mid \text{for all } b \leftarrow a(\cdots) \text{ occurring in a type in } \Gamma \} && (3.5.1) \\ \cup\ & \left\{ \mathsf{t}\,\mathsf{t}\, a = \mathsf{t}\,\mathsf{s}_{[]}\,\mathsf{s}\, a \mid a \in \mathbb{V}_{\Gamma',k} \text{ non degen.}, 2 \leq k \leq n \right\} && (3.5.2) \\ \cup\ & \left\{ \mathsf{t}^{k+2}\, a = b \mid \text{if } \mathsf{t}^k\, a : \underline{b} \multimap b \multimap \cdots, 0 \leq k \leq n-2 \right\}. && (3.5.3) \end{aligned}$$

Here, $\mathsf{t}^k\, a = \mathsf{t} \cdots \mathsf{t}\, a$ is a tagging on the variable $a \in \mathbb{V}_l$, but for simplicity, we consider it as a variable of its own: $\mathsf{t}^k\, a \in \mathbb{V}_{l-k}$. By convention, $\mathsf{t}^0\, a := a$, and if $a = b$, then $\mathsf{t}\, a = \mathsf{t}\, b$, for all $a, b \in \Gamma'$. In line (3.5.2), the source of $a$ is assumed non degenerate, thus $\mathsf{s}\, a$ is a term of the form $x(\overrightarrow{y_i \leftarrow u_i})$, and $\mathsf{s}_{[]}\,\mathsf{s}\, a = x$ (see definition 3.2.1 on page 20).

**Zero:** This rules introduces the empty OCMT.

$$\frac{}{\triangleright}\ \texttt{zero}$$

**Binary sums:** This rule takes two disjoint opetopic sets (i.e. whose cells have different names), and produces their sum. If $\Gamma \cap \Upsilon = \varnothing$, then

$$\frac{E \triangleright \Gamma \qquad F \triangleright \Upsilon}{E \cup F \triangleright \Gamma \cup \Upsilon}\ \texttt{sum}$$

**Quotients:** This rule identifies two parallel cells in an opetopic set by extending the underlying equational theory. If $a, b \in \mathbb{V}_\Gamma$ are such that $\mathsf{s}\, a =_E \mathsf{s}\, b$ and $\mathsf{t}\, a =_E \mathsf{t}\, b$, then

$$\frac{E \triangleright \Gamma}{E \cup \{a = b\} \triangleright \Gamma}\ \texttt{glue}$$

We also write `glue-`$(a{=}b)$ to make explicit that we added $\{a = b\}$ to the theory.

**Remark 3.5.4.** Akin to OPT$^!$, in OPTSET$^!$, an OCMT that is equivalent to a derivable one is itself derivable.

**Remark 3.5.5.** The `sum` and `zero` rules may be replaced by the following `usum` rule (unbiased sum) without changing the set of derivable OCMTs: for $k \geq 0$, $(E_1 \triangleright \Gamma_1), \ldots, (E_k \triangleright \Gamma_k)$ OCMTs such that $\Gamma_i \cap \Gamma_j = \varnothing$ for all $i \neq j$, then

$$\frac{E_1 \triangleright \Gamma_1 \qquad \cdots \qquad E_k \triangleright \Gamma_k}{\sum_{i=1}^k E_i \triangleright \sum_{i=1}^k \Gamma_i}\ \texttt{usum}$$

## 3.6. Equivalence with opetopic sets

3.6.1. **Opetopic sets from OCTM.** Let $E \triangleright \Gamma$ and $F \triangleright \Upsilon$ be two OCMTs. Then we write $\Gamma / E$ for the set $\mathbb{V}_\Gamma$ quotiented by the equivalence relation generated by the equational theory $E$, and likewise for $\Upsilon / F$.

For $(E \triangleright \Gamma \vdash_n x : X)$ a derivable sequent on $\mathrm{Opt}^!$, where $x \in \mathbb{V}_n$, let $(T_x \triangleright C_x)$, the *OCMT of* $x$, be given by

$$\frac{E \triangleright \Gamma \vdash_n x : X}{T_x \triangleright C_x} \, \texttt{repr}$$

We wish to prove that $C_x / T_x$ carries a natural structure of representable opetopic set. We proceed in 4 steps:

1. we show that $C_x / T_x$ is fibered over $\mathbb{O}$ via $[\![-]\!]^{\mathrm{poly}}$ (proposition 3.6.1);
2. we construct the source and target maps, i.e. the structure maps of an opetopic set (proposition 3.6.2);
3. we show that the opetopic identities of section 2.2.4 on page 13 are satisfied (theorem 3.6.3), and that consequently, $C_x / T_x$ has the structure of an opetopic set;
4. finally, we show in proposition 3.6.10 on page 33 that $C_x / T_x$ is in fact a representable opetopic set.

From there, we define a structure of opetopic set on an arbitrary OCMT by induction on its proof tree in equations (3.6.11) and (3.6.12) on page 34.

Let $a \in \mathbb{V}_{C_x, k}$. If $a$ is not a target, then the sequent $(E|_a \triangleright \Gamma|_a \vdash_k a : A)$ is also derivable (where $(-)|_a$ denotes restriction of contexts and theories to variables occurring in type $A$, the type of $a$ in $\Gamma$), and its proof tree is a subtree of that of $x$. Thus we have a well-defined opetope $[\![a]\!]_k^{\mathrm{poly}} \in \mathbb{O}_k$. Otherwise, if $a = \mathsf{t}^{n-k} x$, then $\mathsf{s}\, a = \mathsf{s}^{n-k+1} x$, and define $[\![a]\!]_k^{\mathrm{poly}} = [\![\mathsf{s}^{n-k+1} x]\!]_k^{\mathrm{poly}}$. We thus have a map $[\![-]\!]^{\mathrm{poly}} : \mathbb{V}_{C_x} \longrightarrow \mathbb{O}$.

**Proposition 3.6.1.** *The map $[\![-]\!]^{\mathrm{poly}} : \mathbb{V}_{C_x} \longrightarrow \mathbb{O}$ factors through $C_x / T_x$.*

*Proof.* By construction, the theory $T_x$ identifies variables $a, b \in \mathbb{V}_{C_x, k}$ only if $\mathsf{s}\, a = \mathsf{s}\, b$, thus $[\![a]\!]_k^{\mathrm{poly}} = [\![\mathsf{s}\, a]\!]_k^{\mathrm{poly}} = [\![\mathsf{s}\, b]\!]_k^{\mathrm{poly}} = [\![b]\!]_k^{\mathrm{poly}}$. □

For $\psi \in \mathbb{O}_k$, write

$$(C_x / T_x)_\psi = \left\{ a \in \mathbb{V}_{C_x, k} \mid [\![a]\!]_k^{\mathrm{poly}} = \psi \right\}.$$

We now construct source and target maps between those subsets.

1. (Sources) If $[p] \in [\![a]\!]_k^{\mathrm{poly}\bullet}$, Then, by corollary 3.2.13 on page 23, there is a unique $b \in \mathbb{V}_{C_x|_a, k-1}$ such that $\&_{\mathsf{s}\, a} b = [p]$. Write then $\mathsf{s}_{[p]}\, a = b$.
2. (Target) For $a \in \mathbb{V}_{C_x, k}$, $k > 0$, we of course set $\mathsf{t}(a) = \mathsf{t}\, a$, the latter being a variable introduced by the `repr` rule.

**Proposition 3.6.2.** *Let $a \in \mathbb{V}_{C_x, k}$.*

*(1) For $[p] \in [\![a]\!]_k^{\mathrm{poly}\bullet}$ we have $[\![\mathsf{s}_{[p]}\, a]\!]_{k-1}^{\mathrm{poly}} = \mathsf{s}_{[p]} [\![a]\!]_k^{\mathrm{poly}}$.*
*(2) We have $[\![\mathsf{t}\, a]\!]_{k-1}^{\mathrm{poly}} = \mathsf{t} [\![a]\!]_k^{\mathrm{poly}}$.*

*Proof.* (1) Write $[p] = \&_{\mathsf{s}\, a} b$, for some $b \in (\mathsf{s}\, a)^\bullet$ (see corollary 3.2.13 on page 23). Then, by corollary 3.2.8 on page 21, we have $\mathsf{s}_{[p]} [\![a]\!]_k^{\mathrm{poly}} = [\![b]\!]_{k-1}^{\mathrm{poly}} = [\![\mathsf{s}_{[p]}\, a]\!]_{k-1}^{\mathrm{poly}}$.
(2) By proposition 3.2.11 on page 22, $\mathsf{t} [\![a]\!]_k^{\mathrm{poly}} = [\![\mathsf{s}\,\mathsf{s}\, a]\!]_{k-1}^{\mathrm{poly}} = [\![\mathsf{s}\,\mathsf{t}\, a]\!]_{k-1}^{\mathrm{poly}} = [\![\mathsf{t}\, a]\!]_{k-1}^{\mathrm{poly}}$. □

**Theorem 3.6.3.** *With all the structure introduced above, $C_x / T_x$ is an opetopic set.*

*Proof.* We check the opetopic identities of section 2.2.4 on page 13. Take $a \in \mathbb{V}_{C_x, k}$,

1. **(Inner)** Take $[p[q]] \in [\![a]\!]_k^{\mathrm{poly}\bullet}$, and write $d = \mathsf{s}_{[p[q]]}\, a$. In $\mathsf{s}\, a$, the variable $d$ occurs as
$$\mathsf{s}\, a = \cdots, \ b(c \leftarrow d), \ \cdots$$
for some $b \in^\bullet$ and $c \in (\mathsf{s}\, b)^\bullet$. Then, $[p[q]] = \&_{\mathsf{s}\, a} d = \&_{\mathsf{s}\, a} b \cdot [[\&_{\mathsf{s}\, b} c]]$, and thus $\&_{\mathsf{s}\, a} b = [p]$ and $\&_{\mathsf{s}\, b} c = [q]$. Finally, $\mathsf{s}_{[q]}\, \mathsf{s}_{[p]}\, a = \mathsf{s}_{[q]}\, b = c = \mathsf{t}\, d = \mathsf{t}\, \mathsf{s}_{[p[q]]}\, a$.
2. **(Glob1)** Assume that $a$ is not degenerate. Then, by definition of $T_x$, we have $\mathsf{t}\,\mathsf{t}\, a = \mathsf{t}\, \mathsf{s}_{[]}\, a$.

(3) **(Glob2)** Assume that $a$ is not degenerate, and take $[p[q]] \in \llbracket a \rrbracket_k^{\mathrm{poly}^|}$. Then $[p[q]] = \&_{\mathsf{s}\,a} b$ for some $c \in (\mathsf{s}\,\mathsf{s}\,a)^\bullet$, and further, $c = \mathsf{s}_{[q]}\,\mathsf{s}_{[p]}\,a$. Then

$$\begin{aligned}
\wp_{\llbracket a \rrbracket_k^{\mathrm{poly}}}[p[q]] &= \wp_{\llbracket \mathsf{s}\,a \rrbracket_k^{\mathrm{poly}}}[p[q]] && \text{by def. } \llbracket \mathsf{s}\,a \rrbracket_k^{\mathrm{poly}} = \llbracket a \rrbracket_k^{\mathrm{poly}} \\
&= \wp_{\llbracket \mathsf{s}\,a \rrbracket_k^{\mathrm{poly}}} \&_{\mathsf{s}\,a} c \\
&= \&_{\mathsf{s}\,\mathsf{s}\,a} c && \text{by lemma 3.2.10} \\
&= \&_{\mathsf{s}\,\mathsf{t}\,a} c && \text{since } \mathsf{s}\,\mathsf{s}\,a = \mathsf{s}\,\mathsf{t}\,a,
\end{aligned}$$

and thus $\mathsf{s}_{[q]}\,\mathsf{s}_{[p]}\,a = c = \mathsf{s}_{\&_{\mathsf{s}\,\mathsf{t}\,a} c}\,\mathsf{t}\,a = \mathsf{s}_{\wp_{\llbracket a \rrbracket_k^{\mathrm{poly}}}[p[q]]}\,\mathsf{t}\,a$.

(4) **(Degen)** Assume that $a$ is degenerate, say $\mathsf{s}\,a = \underline{b}$. Then $\mathsf{s}_{[]}\,\mathsf{t}\,a = b = \mathsf{t}\,\mathsf{t}\,a = a$, where the last equality comes from the rule `repr` (recall that by convention $a = \mathsf{t}^0\,a$). ☐

**Lemma 3.6.4.** *The opetopic set $C_x/T_x$ is a quotient of the representable $O[\llbracket x \rrbracket_n^{\mathrm{poly}}]$, where the latter is the representable at $\llbracket x \rrbracket_n^{\mathrm{poly}} \in \mathbb{O}_n$ (see section 2.2.4 on page 13).*

*Proof.* Clearly, the poset of cells $\int_{\mathbb{O}} C_x/T_x$ of $C_x/T_x$ has a unique maximum element, namely $x$ itself. Moreover, that element has shape $\llbracket x \rrbracket_n^{\mathrm{poly}}$, i.e. $x \in (C_x/T_x)_{\llbracket x \rrbracket_n^{\mathrm{poly}}}$. Then, by the Yoneda lemma, there is a map $O[\llbracket x \rrbracket_n^{\mathrm{poly}}] \longrightarrow C_x/T_x$ having cell $x$ in its image, and since $x$ is a maximum, the map is surjective. ☐

Let $(E \triangleright \Gamma \vdash_n x : X)$ be a derivable sequent, with $x \in \mathbb{V}_n$. In lemma 3.6.4, we established that $C_x/T_x$ is a quotient of the representable opetopic set $O[\llbracket x \rrbracket_n^{\mathrm{poly}}]$. We now aim to show that the two are actually isomorphic (proposition 3.6.10 on the next page) by showing that they have the same number of cells: for $\omega \in \mathbb{O}$, let

$$\#\,\omega := \sum_{\psi \in \mathbb{O}} \#\,O[\omega]_\psi = \sum_{\psi \in \mathbb{O}} \#\,\mathbb{O}(\psi, \omega),$$

which is a finite number since the slice category $\mathbb{O}/\omega$ is finite (see the definition of $\mathbb{O}$ in section 2.2.4 on page 13). The strategy of the proof of proposition 3.6.10 on the next page is to show that the number of cells in $C_x/T_x$ is precisely $\#\,\llbracket x \rrbracket_n^{\mathrm{poly}}$. We need some preliminary results first.

**Proposition 3.6.5.** *(1) We have $\#\,\blacklozenge = 1$, and $\#\,\blacksquare = 3$.*

*(2) For a non degenerate opetope $\omega \in \mathbb{O}_n$, with $n \geq 2$, we have*

$$\#\,\omega = 2 + \left( \sum_{[p] \in \omega^\bullet} \#\,\mathsf{s}_{[p]}\,\omega \right) - \left( \sum_{[p[q]] \in \omega^\bullet} \#\,\mathsf{s}_{[q]}\,\mathsf{s}_{[p]}\,\omega \right) \tag{3.6.6}$$

*(3) If $\omega$ is a degenerate opetope, say $\omega = \mathsf{I}_\phi$, then $\#\,\omega = 2 + \#\,\phi$.*

*Proof.* This is a straightforward enumeration exercise from the definition of $\mathbb{O}$ by generators and relations (see section 2.2.4 on page 13), and using opetopic induction (see remark 2.2.3 on page 12). ☐

**Corollary 3.6.7.** *Let $\omega \in \mathbb{O}_n$, for $n \geq 1$.*

*(1) If $\omega = \mathsf{Y}_\psi$ for some $\psi \in \mathbb{O}_{n-1}$ then $\#\,\omega = 2 + \#\,\psi$.*

*(2) If $\omega = \nu \circ_{[l]} \mathsf{Y}_\psi$, for some $\nu \in \mathbb{O}_n$, $[l] \in \nu^|$, and $\psi \in \mathbb{O}_{n-1}$, then*

$$\#\,\omega = \#\,\nu + \#\,\psi - \#\,\mathsf{e}_{[l]}\,\nu. \tag{3.6.8}$$

*Proof.* (1) Using proposition 3.6.5, we have

$$\#\,\omega = 2 + \left( \sum_{[p] \in \omega^\bullet} \#\,\mathsf{s}_{[p]}\,\omega \right) - \left( \sum_{[p[q]] \in \omega^\bullet} \#\,\mathsf{s}_{[q]}\,\mathsf{s}_{[p]}\,\omega \right) = 2 + \#\,\psi.$$

(2) Using proposition 3.6.5 on the preceding page, we have

$$\begin{aligned}\#\,\omega &= 2 + \left(\sum_{[p]\in\omega^\bullet} \#\,\mathsf{s}_{[p]}\,\omega\right) - \left(\sum_{[p[q]]\in\omega^\bullet} \#\,\mathsf{s}_{[q]}\,\mathsf{s}_{[p]}\,\omega\right) \\ &= 2 + \left(\#\,\psi + \sum_{[p]\in\nu^\bullet} \#\,\mathsf{s}_{[p]}\,\nu\right) - \left(\#\,\mathsf{e}_{[l]}\,\nu + \sum_{[p[q]]\in\nu^\bullet} \#\,\mathsf{s}_{[q]}\,\mathsf{s}_{[p]}\,\nu\right) \\ &= \#\,\nu + \#\,\psi - \#\,\mathsf{e}_{[l]}\,\nu.\end{aligned}$$

□

**Examples 3.6.9.** Consider the opetopic integer $\mathbf{n} \in \mathbb{O}_2$ from example 3.3.2 on page 26. We show that $\#\,\mathbf{n} = 2n+3$. If $n = 0$, then by definition, $\#\,\mathbf{0} = \#\,\mathsf{I}_\blacklozenge = 2 + \#\,\blacklozenge = 3$. If $n = 1$, then $\#\,\mathbf{1} = \#\,\mathsf{Y}_\blacksquare = 2 + \#\,\blacksquare = 5$ by corollary 3.6.7 on the preceding page. Otherwise,

$$\begin{aligned}\#\,\mathbf{n} &= \#\left((\mathbf{n}-\mathbf{1}) \underset{[*^{n-1}]}{\circ} \mathsf{Y}_\blacksquare\right) && \text{by def. of } \mathbf{n} \\ &= \#(\mathbf{n}-\mathbf{1}) + \#\,\blacksquare - \#\,\mathsf{e}_{[*^{n-1}]}(\mathbf{n}-\mathbf{1}) && \text{by eq. (3.6.8)} \\ &= (2n+1) + 3 - \#\,\blacklozenge && \text{by ind.} \\ &= (2n+1) + 3 - 1 = 2n+3.\end{aligned}$$

As an other example, consider the 3-opetope $\omega = \mathsf{Y}_\mathbf{2} \circ_{[[*]]} \mathsf{Y}_\mathbf{0}$ of example 3.3.5 on page 27:

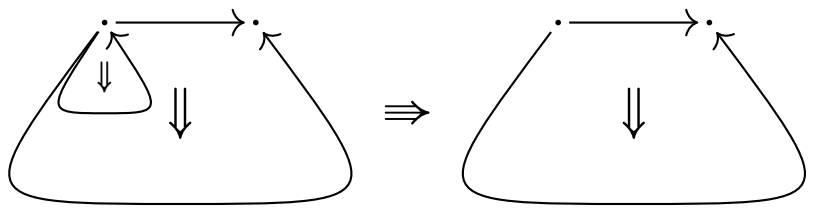

Then,

$$\begin{aligned}\#\,\omega &= \#\left(\mathsf{Y}_\mathbf{2} \underset{[[*]]}{\circ} \mathsf{Y}_\mathbf{0}\right) \\ &= \#\,\mathsf{Y}_\mathbf{2} + \#\,\mathbf{0} - \#\,\mathsf{e}_{[[*]]}\,\mathsf{Y}_\mathbf{2} && \text{by eq. (3.6.8)} \\ &= 2 + \#\,\mathbf{2} + \#\,\mathbf{0} - \#\,\blacksquare && \text{by coroll. 3.6.7} \\ &= 9 && \text{since } \#\,\mathbf{n} = 2n+3.\end{aligned}$$

**Proposition 3.6.10.** *We have an isomorphism $C_x/T_x \cong O[\![x]\!]^{\text{poly}}_n]$ of opetopic sets.*

*Proof.* If $x = \blacklozenge$, then $\mathbb{V}_{C_x} = \blacklozenge$, while $T_x = \varnothing$. Thus $\#\,C_x/T_x = 1 = \#\,\blacklozenge$ by proposition 3.6.5 on the preceding page. We know by lemma 3.6.4 on the facing page that $C_x/T_x$ is a quotient of $O[[\![x]\!]^{\text{poly}}_n]$, and we just showed that the two have the same number of cells, namely $\#\,\blacklozenge = 1$. Consequently, $C_x/T_x \cong O[[\![x]\!]^{\text{poly}}_0]$.

Likewise, if $x = \blacksquare$, then $\mathbb{V}_{C_x} = \{\blacklozenge, \blacksquare, \mathsf{t}\,\blacksquare\}$, while $T_x = \varnothing$. Thus $\#\,C_x/T_x = 3 = \#\,\blacksquare$ by proposition 3.6.5 on the preceding page, and by the same argument as above, $C_x/T_x \cong O[[\![x]\!]^{\text{poly}}_1]$.

Assume now that $x \in \mathbb{V}_n$ for $n \geq 2$. We proceed by cases on the form of $\mathsf{s}\,x$.

(1) If $\mathsf{s}\,x = y \in \mathbb{V}_{n-1}$, then $[\![x]\!]^{\text{poly}}_n = \mathsf{Y}_{[\![y]\!]^{\text{poly}}_{n-1}}$ so that $\#\,[\![x]\!]^{\text{poly}}_n = 2 + \#\,[\![y]\!]^{\text{poly}}_{n-1}$ by proposition 3.6.5 on the facing page. Then $C_x = C_y + \left\{\mathsf{t}^k\,x \mid 0 \leq k \leq n\right\}$, and $\mathsf{t}\,\mathsf{t}\,x T_x\,\mathsf{t}\,\mathsf{s}_{[]}\,x = \mathsf{t}\,y$. Consequently, $T_x$ is equivalent to the theory $T_y + \{\mathsf{t}\,\mathsf{t}\,x = \mathsf{t}\,y\}$, and thus

$$C_x/T_x = C_y/T_y + \{x, \mathsf{t}\,x\}\,.$$

By induction, $C_y/T_y \cong O[[\![y]\!]^{\text{poly}}_{n-1}]$, and $\#\,C_x/T_x = 2 + \#\,C_y/T_y = 2 + \#\,[\![y]\!]^{\text{poly}}_{n-1} = \#\,[\![x]\!]^{\text{poly}}_n$, which, by the same argument as above, proves the isomorphism $C_x/T_x \cong O[[\![x]\!]^{\text{poly}}_n]$.

(2) If $\mathsf{s}\,x = \underline{a}$ for some $a \in \mathbb{V}_{n-2}$, then $[\![x]\!]^{\text{poly}}_n = \mathsf{I}_{[\![a]\!]^{\text{poly}}_{n-2}}$ so that $\#\,[\![x]\!]^{\text{poly}}_n = 2 + \#\,[\![a]\!]^{\text{poly}}_{n-2}$ by proposition 3.6.5 on the preceding page. Then $C_x = C_a + \left\{\mathsf{t}^k\,x \mid 0 \leq k \leq n\right\}$, and $\mathsf{t}\,\mathsf{t}\,x T_x a$. Therefore $T_x$ is equivalent to the theory $T_a + \{\mathsf{t}\,\mathsf{t}\,x = a\}$, and thus

$$C_x/T_x = C_a/T_a + \{x, \mathsf{t}\,x\}\,.$$

Consequently, $\#\,C_x/T_x = 2 + \#\,C_a/T_a = 2 + \#\,[\![a]\!]^{\text{poly}}_{n-2} = \#\,[\![x]\!]^{\text{poly}}_n$.

(3) Assume $\mathsf{s}\, x = t(a \leftarrow y)$, for some $t \in \mathbb{T}_{n-1}$, $a \in \mathbb{V}_{n-2}$, and $y \in \mathbb{V}_{n-1}$. Let $z : t \multimap \cdots$ be a fresh $n$-variable. Since all cells of $z$ except $z$ and its targets are also cells of $x$, we have

$$C_x = C_y \cup \left(C_z - \left\{\mathsf{t}^k\, z \mid 0 \le k \le n\right\}\right) + \left\{\mathsf{t}^k\, x \mid 0 \le k \le n\right\},$$

while $T_x$ is equivalent to $T_y \cup T_z + \{\mathsf{t}\, y = a, \mathsf{t}\,\mathsf{t}\, x = \mathsf{t}\,\mathsf{s}_{[]}\, x\}$. Since $\mathsf{s}_{[]}\, x = \mathsf{s}_{[]}\, z$, and $\mathsf{t}\,\mathsf{s}_{[]}\, z T_z\, \mathsf{t}\,\mathsf{t}\, z$, we have

$$C_x/T_x = C_y/T_y \cup C_z/T_z + \{x, \mathsf{t}\, x\} - \{z, \mathsf{t}\, z\}.$$

By hypothesis of the `graft` rule, $C_z/T_z \cap C_y/T_y = C_a/T_a$, thus

$$\begin{aligned} \# C_x/T_x &= \# C_y/T_y + \# C_z/T_z - \# C_a/T_a = \# [\![y]\!]^{\text{poly}}_{n-1} + \# [\![z]\!]^{\text{poly}}_{n} - \# [\![a]\!]^{\text{poly}}_{n-2} \\ &= \# \left( [\![t]\!]^{\text{poly}}_{n} \underset{\&_t a}{\circ} \mathsf{Y}_{[\![y]\!]^{\text{poly}}_{n-1}} \right) = \# [\![x]\!]^{\text{poly}}_{n}. \end{aligned}$$

□

Let $(E \triangleright \Gamma)$ be a derivable OCMT in $\textsc{OptSet}^!$, $a \in \mathbb{V}_{\Gamma,k}$, and $(E|_a \triangleright \Gamma|_a)$ the restriction to $a$. Then we have an inclusion $i : \Gamma|_a/E|_a \hookrightarrow \Gamma/E$, where $i : b \longmapsto b$ for all $b \in \mathbb{V}_{\Gamma|_a}$, and thus we have a natural map $\tilde{a} : O[[\![a]\!]^{\text{poly}}_k] \longrightarrow \Gamma/E$ given by the composite

$$O[[\![a]\!]^{\text{poly}}_k] \xrightarrow{\cong} C_a/T_a \longrightarrow \Gamma|_a/E|_a \xrightarrow{i} \Gamma/E,$$

where the middle map comes from the fact that $C_a = \Gamma|_a$ up to renaming, and that $T_a \subseteq E|_a$. Let now $(E \triangleright \Gamma)$ and $(F \triangleright \Upsilon)$ be two OCMTs, and assume by induction that $\Gamma/E$ and $\Upsilon/F$ have a structure of opetopic set. Then, by definition of rule `sum`, and for $(G \triangleright \Xi)$ given by

$$\frac{E \triangleright \Gamma \qquad F \triangleright \Upsilon}{G \triangleright \Xi}\ \texttt{sum}$$

we have $\Xi = \Gamma + \Upsilon$, and $G = E + F$. We endow the quotient $\Xi/G$ with a structure of opetopic set as follows:

$$\Xi/G := \Gamma/E + \Upsilon/F. \tag{3.6.11}$$

Let $a, b \in \mathbb{V}_{\Gamma,k}$ be such that $\mathsf{s}\, a =_E \mathsf{s}\, b$ and $\mathsf{t}\, a =_E \mathsf{t}\, b$. Remark that $O[[\![a]\!]^{\text{poly}}_k] = O[[\![b]\!]^{\text{poly}}_k]$. Then, by definition of rule `glue`, and for $(F \triangleright \Gamma)$ given by

$$\frac{E \triangleright \Gamma}{F \triangleright \Gamma}\ \texttt{glue-}(a{=}b)$$

we have $F = E + \{a = b\}$. We endow the quotient $\Gamma/F$ with a structure of opetopic set defined by the following coequalizer:

$$O[[\![a]\!]^{\text{poly}}_k] \underset{\tilde{b}}{\overset{\tilde{a}}{\rightrightarrows}} \Gamma/E \longrightarrow \Gamma/F, \tag{3.6.12}$$

**Proposition 3.6.13.** *For $(E \triangleright \Gamma)$ a derivable OCMT in $\textsc{OptSet}^!$, the structure of opetopic set on $\Gamma/E$ does not depend on the proof tree of $(E \triangleright \Gamma)$.*

*Proof.* The opetopic set $\Gamma/E$ is given by the following expression that does not depend on the proof tree of $(E \triangleright \Gamma)$:

$$\Gamma/E \cong \frac{\sum_{k \in \mathbb{N}, a \in \mathbb{V}_{\Gamma,k}} O[[\![a]\!]^{\text{poly}}_k]}{a \sim b, \forall a, b \in \mathbb{V}_\Gamma \text{ s.t. } a =_E b}.$$

Indeed, it holds when $(E \triangleright \Gamma)$ is of the form $(T_a \triangleright C_a)$, and since other non empty OCMTs are given by application of rules `sum` and `glue`, we can show by induction on the proof tree of $(E \triangleright \Gamma)$ that the opetopic set $\Gamma/E$ doesn't depend on that proof tree! □

3.6.2. **Equivalence.** Recall that $\widehat{\mathbb{O}}$ is the category of opetopic sets, and that $\mathcal{F}\mathrm{in}\widehat{\mathbb{O}}$ is the full subcategory of $\widehat{\mathbb{O}}$ spanned by finite opetopic sets (see section 2.2.4 on page 13). In this subsection, we provide the last results needed to establish the equivalence between the category of derivable OCMTs and $\mathcal{F}\mathrm{in}\widehat{\mathbb{O}}$.

For $(E \triangleright \Gamma)$ an OCMT, and $a, b \in \mathbb{V}$, the substitution $\Gamma[a/b]$ is defined in the obvious manner, by applying said substitution to all typings in $\Gamma$. A morphism $f : (E \triangleright \Gamma) \longrightarrow (F \triangleright \Upsilon)$ is a (non necessarily bijective) map $f : \mathbb{V}_\Gamma \longrightarrow \mathbb{V}_\Upsilon$ compatible with $E$ and $F$, such that if $x : X$ is a typing in $\Gamma$, then $f(x) : f(X)$ is a typing in $\Upsilon$, where $f(X)$ is the result of applying $f$ to every variable in $X$. Note that this condition implies that $f$ preserves the dimension of variables. Also, if $f, g : (E \triangleright \Gamma) \longrightarrow (F \triangleright \Upsilon)$, and if for all $x \in \mathbb{V}_\Gamma$ we have $f(x) =_F g(x)$, then $f = g$. We write $\mathcal{C}\mathrm{tx}^!$ for the category of derivable OCMTs and such morphisms. In a sense, it is the syntactic category of the $\textsc{OptSet}^!$ system.

**Lemma 3.6.14.** *Let $f : (E \triangleright \Gamma) \longrightarrow (F \triangleright \Upsilon)$ be a morphism of OCMT, and $a \in \mathbb{V}_{\Gamma,k}$. Then $[\![a]\!]_k^{\mathrm{poly}} = [\![f(a)]\!]_k^{\mathrm{poly}}$. Informally, morphisms of OCMTs preserve the "shape" of variables.*

*Proof.* Since there is a unique 0-opetope and a unique 1-opetope, the result holds trivially if $k = 0, 1$. Assume now $k \geq 2$. We proceed by induction on $\mathsf{s}\, a$.

(1) If $\mathsf{s}\, a = b \in \mathbb{V}_{k-1}$, then $[\![a]\!]_k^{\mathrm{poly}} = \mathsf{Y}_{[\![b]\!]_{k-1}^{\mathrm{poly}}} = \mathsf{Y}_{[\![f(b)]\!]_{k-1}^{\mathrm{poly}}} = [\![f(a)]\!]_k^{\mathrm{poly}}$.
(2) If $\mathsf{s}\, a = \underline{b}$ for some $b \in \mathbb{V}_{k-2}$, then $[\![a]\!]_k^{\mathrm{poly}} = \mathsf{I}_{[\![b]\!]_{k-2}^{\mathrm{poly}}} = \mathsf{I}_{[\![f(b)]\!]_{k-2}^{\mathrm{poly}}} = [\![f(a)]\!]_{k-2}^{\mathrm{poly}}$.
(3) If $\mathsf{s}\, a = b(\overrightarrow{c_i \leftarrow u_i})$, we can show by a secondary induction that $[\![u_i]\!]_k^{\mathrm{poly}} = [\![f(u_i)]\!]_k^{\mathrm{poly}}$, and

$$[\![a]\!]_k^{\mathrm{poly}} = \mathsf{Y}_{[\![b]\!]_{k-1}^{\mathrm{poly}}} \bigcirc_{[\&_{\mathsf{s}\, b} c_i]} [\![u_i]\!]_k^{\mathrm{poly}} = [\![f(a)]\!]_k^{\mathrm{poly}}$$

□

The *named stratification functor* $S^! : \mathcal{C}\mathrm{tx}^! \longrightarrow \mathcal{F}\mathrm{in}\widehat{\mathbb{O}}$ is defined as follows:

$$\begin{aligned} S^! : \mathcal{C}\mathrm{tx}^! &\longrightarrow \mathcal{F}\mathrm{in}\widehat{\mathbb{O}} \\ E \triangleright \Gamma &\longmapsto \Gamma/E \\ \left((E \triangleright \Gamma) \xrightarrow{f} (F \triangleright \Upsilon)\right) &\longmapsto \left(\Gamma/E \xrightarrow{f} \Upsilon/F\right). \end{aligned}$$

By lemma 3.6.14, $f = S^! f : S^!(E \triangleright \Gamma) \longrightarrow S^!(F \triangleright \Upsilon)$ is indeed a morphism of opetopic sets.

**Theorem 3.6.15.** *The stratification functor $S^! : \mathcal{C}\mathrm{tx}^! \longrightarrow \mathcal{F}\mathrm{in}\widehat{\mathbb{O}}$ is an equivalence of categories.*

*Proof.* The category $\mathrm{im}\, D$, consisting of finite opetopic sets of the form $\Gamma/E$ for some OCMT $E \triangleright \Gamma$, contains all representables (proposition 3.6.10 on page 33), the initial object (since $S^!(\triangleright)$ is the opetopic set with no cell), and is closed under finite sums (equation (3.6.11) on the preceding page) and quotients (equation (3.6.12) on the facing page). Thus its closure under isomorphism is finitely cocomplete, and equal to the whole category $\mathcal{F}\mathrm{in}\widehat{\mathbb{O}}$, so $S^!$ is essentially surjective. By definition, $S^!$ is also faithful, and it remains to show that is it full.

Let $f : S^!(E \triangleright \Gamma) \longrightarrow S^!(F \triangleright \Upsilon)$ be a morphism of opetopic sets. Then, in particular, it is a map $f : \Gamma/E \longrightarrow \Upsilon/F$ between the set of cells of $S^!(E \triangleright \Gamma)$ and $S^!(F \triangleright \Upsilon)$. To prove that it is a morphism of OCMT, we show that $\Gamma[f(x)/x \mid x \in \mathbb{V}_\Gamma]$ is a subcontext of $\Upsilon$ modulo $F$, i.e. that for all typings $x : X$ in $\Gamma$, for some $x \in \mathbb{V}_k$, the type of $f(x)$ in $\Upsilon$ is $f(X)$ modulo $F$. If $(E \triangleright \Gamma)$ is the empty OCMT, the result is trivial, so we assume now that this is not the case. Since $f$ is a morphism of opetopic sets, we have $f(x) \in \mathbb{V}_{\Upsilon,k}$, and by lemma 3.6.14, $[\![x]\!]_k^{\mathrm{poly}} = [\![f(x)]\!]_k^{\mathrm{poly}}$. We proceed by induction on $k$.

(1) If $k = 0$, then $X = \varnothing$. Since $f(x) \in \mathbb{V}_{\Upsilon,0}$, its type is necessarily $\varnothing = f(X)$, thus $f(x) : f(X)$ is a typing in $\Upsilon$.
(2) If $k = 1$, then $X = \left(\mathsf{s}_{[]}\, x \multimap \varnothing\right)$, and, since $f$ is a morphism of opetopic sets, $f(\mathsf{s}_{[]}\, x) =_F \mathsf{s}_{[]}\, f(x)$. Thus

$$f(X) = \left(f(\mathsf{s}_{[]}\, x) \multimap \varnothing\right) =_F \left(\mathsf{s}_{[]}\, f(x) \multimap \varnothing\right),$$

the latter being the type of $f(x)$ in $\Upsilon$.

(3) Assume now that $k \geq 2$. The type $X$ of $x$ in $\Gamma$ is by definition $X = (\mathsf{s}x \multimapdotbothA \mathsf{ss}x \multimapdotbothA \cdots \multimapdotbothA \varnothing)$, and the type of $\mathsf{t}x$ in $\Gamma$ is $Y := (\mathsf{ss}x \multimapdotbothA \cdots \multimapdotbothA \varnothing)$. By induction, the type of $f(\mathsf{t}x)$ in $\Upsilon$ is $f(Y)$, and since $f(\mathsf{t}x) = \mathsf{t}f(x)$, and the type of the latter is by definition $\mathsf{ss}f(x) \multimapdotbothA \cdots \multimapdotbothA \varnothing$, we have

$$(f(\mathsf{ss}x) \multimapdotbothA \cdots \multimapdotbothA \varnothing) = f(Y) =_F (\mathsf{ss}f(x) \multimapdotbothA \cdots \multimapdotbothA \varnothing),$$

or in other words, $\mathsf{s}^i f(x) =_F f(\mathsf{s}^i x)$, for $2 \leq i \leq k$. It remains to show that the latter formula holds in the case $i = 1$ (the case $i = 0$ is tautological). Towards a contradiction, assume $\mathsf{s}f(x) \neq_F f(\mathsf{s}x)$. Then there exists $[p] \in [\![x]\!]_k^{\mathrm{poly}^\bullet} = [\![f(x)]\!]_k^{\mathrm{poly}^\bullet}$ such that $\mathsf{s}_{[p]} f(x) \neq_F f(\mathsf{s}_{[p]} x)$, which contradicts the fact that $f$ is a morphism of opetopic sets. Consequently, $\mathsf{s}f(x) =_F f(\mathsf{s}x)$, and $f(X)$ is the type of $f(x)$ in $\Upsilon$ modulo $F$.

Finally, the underlying map $f : \Gamma/E \longrightarrow \Upsilon/F$ of the morphism of opetopic sets $f : S^!(E \triangleright \Gamma) \longrightarrow S^!(F \triangleright \Upsilon)$ is a morphism of OCMT, and $S^!$ is full. □

**Proposition 3.6.16.** *The category* $\mathcal{C}\mathrm{tx}^!$ *has finite colimits, and* $S^!$ *preserves them.*

*Proof.* The empty OCMT is certainly an initial object of $\mathcal{C}\mathrm{tx}^!$ (it is in fact the only one), and the OCMT $(E + F \triangleright \Gamma + \Upsilon)$ obtained by the rule `sum` is clearly a coproduct of $(E \triangleright \Gamma)$ and $(F \triangleright \Upsilon)$. Let now $f, g : (E \triangleright \Gamma) \longrightarrow (F \triangleright \Upsilon)$ be two parallel morphism in $\mathcal{C}\mathrm{tx}^!$. Consider the following:

$$(F \triangleright \Upsilon) \xrightarrow{q : x \longmapsto x} (\underbrace{F \cup \{f(x) = g(x) \mid x \in \mathbb{V}_\Gamma\}}_{G} \triangleright \Upsilon).$$

Then, $G \triangleright \Upsilon$ is derived from $F \triangleright \Upsilon$ by repeated application of the `glue` rule, and by universal property, it is easy to see that it is a coequalizer of $f$ and $g$ in $\mathcal{C}\mathrm{tx}^!$. So, to summarize, $\mathcal{C}\mathrm{tx}^!$ contains all finite sums and coequalizers, and it is thus finitely cocomplete. The fact that $S^!$ preserves finite colimits is trivial in the case of initial objects, a consequence of equation (3.6.11) on page 34 for binary sums, and of equation (3.6.12) on page 34 for coequalizers. □

**Theorem 3.6.17.** *We have an equivalence* $\widehat{\mathbb{O}} \simeq \mathcal{M}\mathrm{od}((\mathcal{C}\mathrm{tx}^!)^{\mathrm{op}}, \mathcal{S}\mathrm{et})$*, where the latter is the category of* $(\mathcal{C}\mathrm{tx}^!)^{\mathrm{op}}$*-models in* $\mathcal{S}\mathrm{et}$*, i.e. finitely continuous functors* $(\mathcal{C}\mathrm{tx}^!)^{\mathrm{op}} \longrightarrow \mathcal{S}\mathrm{et}$ *and natural transformations*[2]*.*

*Proof.* We apply Gabriel–Ulmer duality [**Adámek and Rosický, 1994, Lack and Power, 2009**]. The category $\widehat{\mathbb{O}}$ is finitely locally presentable as it is the presheaf category of a small category. Finite opetopic sets are the finitely presentable objects, thus

$$\mathcal{M}\mathrm{od}\left((\mathcal{C}\mathrm{tx}^!)^{\mathrm{op}}, \mathcal{S}\mathrm{et}\right) \simeq \mathcal{M}\mathrm{od}\left(\mathcal{F}\mathrm{in}\widehat{\mathbb{O}}^{\mathrm{op}}, \mathcal{S}\mathrm{et}\right) \simeq \widehat{\mathbb{O}}.$$

□

## 3.7. Examples

In this section, we give example derivations in system OptSet$^!$. For clarity, we do not repeat the type of previously typed variables in proof trees.

**Example 3.7.1.** The opetopic set

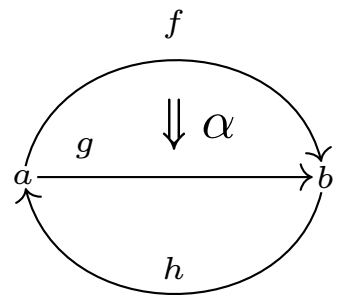

is derived as follows. First, we derive the cells $\alpha$, $g$, and $h$ as opetopes (i.e. in OptSet$^!$) to obtain the following sequents:

$$\triangleright a : \varnothing, f : a \multimapdotbothA \varnothing, \alpha : f \multimapdotbothA a \multimapdotbothA \varnothing \vdash_2 \alpha : f \multimapdotbothA a \multimapdotbothA \varnothing$$

$$\triangleright c : \varnothing, g : c \multimapdotbothA \varnothing \vdash_1 g : c \multimapdotbothA \varnothing$$

$$\triangleright b : \varnothing, h : b \multimapdotbothA \varnothing \vdash_1 h : b \multimapdotbothA \varnothing$$

[2]The category of models $\mathcal{M}\mathrm{od}$ is sometimes written $\mathcal{L}\mathrm{ex}$, for *left exact functors*.

and applying the `repr` rule yields respectively:

$$\triangleright\, a : \varnothing, f : a \multimap \varnothing, \alpha : f \multimap a \multimap \varnothing, \mathsf{t}\, f : \varnothing, \mathsf{t}\,\alpha : a \multimap \varnothing, \mathsf{t}\,\mathsf{t}\,\alpha : \varnothing$$
$$\triangleright\, c : \varnothing, g : c \multimap \varnothing, \mathsf{t}\, g : \varnothing$$
$$\triangleright\, b : \varnothing, h : b \multimap \varnothing, \mathsf{t}\, h : \varnothing.$$

The proof tree then reads:

$$\dfrac{\dfrac{\dfrac{\dfrac{\vdots}{\triangleright a, f, \alpha \vdash_2 \alpha}}{\mathsf{tt}\alpha = \mathsf{t}f \triangleright a, f, \alpha, \mathsf{t}f, \mathsf{t}\alpha, \mathsf{tt}\alpha}\,\texttt{repr} \qquad \dfrac{\dfrac{\vdots}{\triangleright c, g \vdash_1 g}}{\triangleright c, g, \mathsf{t}g}\,\texttt{repr}}{\mathsf{tt}\alpha = \mathsf{t}f \triangleright a, f, \alpha, \mathsf{t}f, \mathsf{t}\alpha, \mathsf{tt}\alpha, c, g, \mathsf{t}g}\,\texttt{sum} \qquad \dfrac{\dfrac{\vdots}{\triangleright b, h \vdash_1 h}}{\triangleright b, h, \mathsf{t}h}\,\texttt{repr}}{\mathsf{tt}\alpha = \mathsf{t}f \quad \triangleright \begin{array}{l} a, b, c, \mathsf{t}f, \mathsf{t}g, \mathsf{t}h, \mathsf{tt}\alpha \\ f, g, h, \mathsf{t}\alpha \\ \alpha \end{array}}\,\texttt{sum}$$

$$\dfrac{}{\mathsf{tt}\alpha = \mathsf{t}f, a = c \quad \triangleright \begin{array}{l} a, b, c, \mathsf{t}f, \mathsf{t}g, \mathsf{t}h, \mathsf{tt}\alpha \\ f, g, h, \mathsf{t}\alpha \\ \alpha \end{array}}\,\texttt{glue-}(a = c)$$

$$\dfrac{}{\mathsf{tt}\alpha = \mathsf{t}f, a = c, b = \mathsf{t}f \quad \triangleright \begin{array}{l} a, b, c, \mathsf{t}f, \mathsf{t}g, \mathsf{t}h, \mathsf{tt}\alpha \\ f, g, h, \mathsf{t}\alpha \\ \alpha \end{array}}\,\texttt{glue-}(b = \mathsf{t}f)$$

$$\dfrac{}{\mathsf{tt}\alpha = \mathsf{t}f, a = c, b = \mathsf{t}f, b = \mathsf{t}g \quad \triangleright \begin{array}{l} a, b, c, \mathsf{t}f, \mathsf{t}g, \mathsf{t}h, \mathsf{tt}\alpha \\ f, g, h, \mathsf{t}\alpha \\ \alpha \end{array}}\,\texttt{glue-}(b = \mathsf{t}g)$$

$$\dfrac{}{\mathsf{tt}\alpha = \mathsf{t}f, a = c, b = \mathsf{t}f, b = \mathsf{t}g, a = \mathsf{t}h \quad \triangleright \begin{array}{l} a, b, c, \mathsf{t}f, \mathsf{t}g, \mathsf{t}h, \mathsf{tt}\alpha \\ f, g, h, \mathsf{t}\alpha \\ \alpha \end{array}}\,\texttt{glue-}(a = \mathsf{t}h)$$

$$\dfrac{}{\begin{array}{r} \mathsf{tt}\alpha = \mathsf{t}f, a = c, b = \mathsf{t}f, b = \mathsf{t}g, a = \mathsf{t}h \\ g = \mathsf{t}\alpha \end{array} \quad \triangleright \begin{array}{l} a, b, c, \mathsf{t}f, \mathsf{t}g, \mathsf{t}h, \mathsf{tt}\alpha \\ f, g, h, \mathsf{t}\alpha \\ \alpha \end{array}}\,\texttt{glue-}(g = \mathsf{t}\alpha)$$

**Example 3.7.2.** The opetopic set

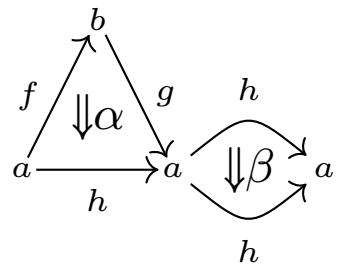

is derived as follows. First, we derive the cells $\alpha$ and $\beta$ as opetopes to obtain the following sequents:

$$\triangleright \begin{array}{l} a : \varnothing, b : \varnothing \\ f : a \multimap \varnothing, g : a \multimap \varnothing \\ \alpha : g(b \leftarrow f) \multimap a \multimap \varnothing \end{array} \quad \vdash_2 \alpha : g(b \leftarrow f) \multimap a \multimap \varnothing$$

$$\triangleright \begin{array}{l} a' : \varnothing \\ h : a' \multimap \varnothing \\ \beta : h \multimap a' \multimap \varnothing \end{array} \quad \vdash_2 \beta : h \multimap a' \multimap \varnothing$$

Applying the `repr` rule yields respectively:

$$
b = \mathtt{t}f, \mathtt{t}g = \mathtt{tt}\alpha \;\triangleright\; \begin{array}{l} a : \varnothing, b : \varnothing, \mathtt{t}f : \varnothing, \mathtt{t}g : \varnothing, \mathtt{tt}\alpha : \varnothing \\ f : a \multimap \varnothing, g : a \multimap \varnothing, \mathtt{t}\alpha : a \multimap \varnothing \\ \alpha : g(b \leftarrow f) \multimap a \multimap \varnothing \end{array}
$$

$$
\mathtt{t}h = \mathtt{tt}\beta \;\triangleright\; \begin{array}{l} a' : \varnothing, \mathtt{t}h : \varnothing, \mathtt{tt}\beta : \varnothing \\ h : a' \multimap \varnothing, \mathtt{t}\beta : a' \multimap \varnothing \\ \beta : h \multimap a' \multimap \varnothing \end{array}
$$

The proof tree then reads:

$$
\dfrac{\dfrac{\dfrac{\vdots}{\triangleright \begin{array}{l} a : \varnothing, b : \varnothing \\ f : a \multimap \varnothing, g : a \multimap \varnothing \\ \alpha : g(b \leftarrow f) \multimap a \multimap \varnothing \end{array} \quad \vdash_2 \alpha : g(b \leftarrow f) \multimap a \multimap \varnothing}{b = \mathtt{t}f, \mathtt{t}g = \mathtt{tt}\alpha \;\triangleright\; \begin{array}{l} a, b, \mathtt{t}f : \varnothing, \mathtt{t}g : \varnothing, \mathtt{tt}\alpha : \varnothing \\ f, g, \mathtt{t}\alpha : a \multimap \varnothing \\ \alpha \end{array}}\,\mathtt{repr}}{b = \mathtt{t}f, a = \mathtt{t}g = \mathtt{tt}\alpha \;\triangleright\; \begin{array}{l} a, b, \mathtt{t}f, \mathtt{t}g, \mathtt{tt}\alpha \\ f, g, \mathtt{t}\alpha \\ \alpha \end{array}}\,\mathtt{glue\text{-}}(a = \mathtt{t}g)
\qquad
\dfrac{\dfrac{\dfrac{\dfrac{\vdots}{\triangleright \begin{array}{l} a' : \varnothing \\ h : a' \multimap \varnothing \\ \beta : h \multimap a' \multimap \varnothing \end{array} \quad \vdash_2 \beta : h \multimap a' \multimap \varnothing}{\mathtt{t}h = \mathtt{tt}\beta \;\triangleright\; \begin{array}{l} a', \mathtt{t}h : \varnothing, \mathtt{tt}\beta : \varnothing \\ h, \mathtt{t}\beta : a' \multimap \varnothing \\ \beta \end{array}}\,\mathtt{repr}}{\begin{array}{l} \mathtt{t}h = \mathtt{tt}\beta \\ h = \mathtt{t}\beta \end{array} \;\triangleright\; \begin{array}{l} a', \mathtt{t}h, \mathtt{tt}\beta \\ h, \mathtt{t}\beta \\ \beta \end{array}}\,\mathtt{glue\text{-}}(h = \mathtt{t}\beta)}{\begin{array}{l} a' = \mathtt{t}h = \mathtt{tt}\beta \\ h = \mathtt{t}\beta \end{array} \;\triangleright\; \begin{array}{l} a', \mathtt{t}h, \mathtt{tt}\beta \\ h, \mathtt{t}\beta \\ \beta \end{array}}\,\mathtt{glue\text{-}}(a' = \mathtt{t}h)
$$

$$
\dfrac{\dfrac{\dfrac{\cdots}{\begin{array}{l} b = \mathtt{t}f, a = \mathtt{t}g = \mathtt{tt}\alpha, a' = \mathtt{t}h = \mathtt{tt}\beta \\ h = \mathtt{t}\beta \end{array} \;\triangleright\; \begin{array}{l} a, b, \mathtt{t}f, \mathtt{t}g, \mathtt{tt}\alpha, a', \mathtt{t}h, \mathtt{tt}\beta \\ f, g, \mathtt{t}\alpha, h, \mathtt{t}\beta \\ \alpha, \beta \end{array}}\,\mathtt{sum}}{\begin{array}{l} b = \mathtt{t}f, a = \mathtt{t}g = \mathtt{tt}\alpha = a' = \mathtt{t}h = \mathtt{tt}\beta \\ h = \mathtt{t}\beta \end{array} \;\triangleright\; \begin{array}{l} a, b, \mathtt{t}f, \mathtt{t}g, \mathtt{tt}\alpha, a', \mathtt{t}h, \mathtt{tt}\beta \\ f, g, \mathtt{t}\alpha, h, \mathtt{t}\beta \\ \alpha, \beta \end{array}}\,\mathtt{glue\text{-}}(a = a')}{\begin{array}{l} b = \mathtt{t}f, a = \mathtt{t}g = \mathtt{tt}\alpha = a' = \mathtt{t}h = \mathtt{tt}\beta \\ h = \mathtt{t}\beta = \mathtt{t}\alpha \end{array} \;\triangleright\; \begin{array}{l} a, b, \mathtt{t}f, \mathtt{t}g, \mathtt{tt}\alpha, a', \mathtt{t}h, \mathtt{tt}\beta \\ f, g, \mathtt{t}\alpha, h, \mathtt{t}\beta \\ \alpha, \beta \end{array}}\,\mathtt{glue\text{-}}(\mathtt{t}\alpha = h)
$$

Write $(E \triangleright \Gamma)$ for the final OCMT of this proof tree. At the beginning of section 3.5.1 on page 29, we gave a slightly different OCMT for the same opetopic set:

$$
(F \triangleright \Upsilon) := \left( \begin{array}{l} b = \mathtt{t}f, a = \mathtt{t}g = \mathtt{tt}\alpha = \mathtt{t}h = \mathtt{tt}\beta \\ h = \mathtt{t}\beta = \mathtt{t}\alpha \end{array} \;\triangleright\; \begin{array}{l} a : \varnothing, b : \varnothing, \mathtt{t}f : \varnothing, \mathtt{t}g : \varnothing, \mathtt{tt}\alpha : \varnothing, \mathtt{t}h : \varnothing, \mathtt{tt}\beta : \varnothing \\ f : a \multimap \varnothing, g : a \multimap \varnothing, \mathtt{t}\alpha : a \multimap \varnothing, h : a \multimap \varnothing, \mathtt{t}\beta : a \multimap \varnothing \\ \alpha : g(b \leftarrow f) \multimap a \multimap \varnothing, \beta : h \multimap a \multimap \varnothing \end{array} \right).
$$

However, there exists an isomorphism $f : (E \triangleright \Gamma) \longrightarrow (F \triangleright \Upsilon)$, defined by mapping $a'$ to $a$, and variables in $\mathbb{V}_\Gamma - \{a'\}$ to those in $\mathbb{V}_\Upsilon$ with the same name. Its inverse it given by mapping $a$ to either $a$ or $a'$, and similarly for the over variables of $\Upsilon$.

## 3.8. Python implementation

The system $\text{OptSet}^{!}$ is implemented in module `opetopy.NamedOpetopicSet` of [**Ho Thanh, 2018b**]. The rules are represented by functions `repres` (since `repr` a Python standard function), `sum`, `glue`, and `zero`, and are further encapsulated in rule instance classes `Repr`, `Sum`, `Glue`, and `Zero`. We do not discuss the implementation, but we give an example derivation in figure 3.8.1 on the next page.

FIGURE 3.8.1. Derivation of example 3.7.1 on page 36 using `opetopy.NamedOpetopicSet`

```python
from opetopy.NamedOpetopicSet import *
# We first define all relevant variables using system Opt!.
f = Shift(Point("a"), "f")
g = Shift(Point("c"), "g")
h = Shift(Point("b"), "h")
alpha = Shift(f, "alpha")
# We then take the sum of all the representables we need. Note that the new target variables added by the repr rule have "t" prepended to their name e.g. the target variable of f is "tf", while that of α is "talpha".
example_unglued = Sum(Sum(Repr(alpha), Repr(g)), Repr(h))
example = Glue(Glue(Glue(Glue(Glue(example_unglued,
                                   "a",
                                   "c"),
                              "b",
                              "tf"),
                         "b",
                         "tg"),
                    "a",
                    "th"),
               "g",
               "talpha")
```

## 3.9. The mixed system for opetopic sets

The $\mathrm{OptSet}^!$ system, presented in section 3.5 on page 29, suffers from the following drawback: derivations of opetopic sets start with instances of rules `zero` or `repr`, the latter requiring a full opetope derivation in system $\mathrm{Opt}^!$ (presented in section 3.1 on page 15). This makes derivations somewhat unintuitive, since for an opetopic set $X \in \widehat{\mathbb{O}}$ written as

$$X = \frac{\sum_i O[\omega_i]}{\sim}$$

where $\sim$ represents some quotient, the opetopes $\omega_i$ have to be derived in $\mathrm{Opt}^!$ first, then the `repr` rule has to be used on each one to produce the corresponding representables $O[\omega_i]$, and only then can the sums and gluing be performed.

In this section, we present system $\mathrm{OptSet}^!_{\mathrm{M}}$ (the M standing for "mixed") for opetopic sets, which does not depend on $\mathrm{Opt}^!$, and allows to perform introductions of new cells, sums, and gluings in any sound order. This is done by introducing new cells along with all their targets, effectively rendering $\mathrm{OptSet}^!$'s `repr` rule superfluous, and removing the "barrier" between $\mathrm{Opt}^!$ and $\mathrm{OptSet}^!$ in the above schema.

3.9.1. **Syntax.** The syntax of system $\mathrm{OptSet}^!_{\mathrm{M}}$ uses sequents from $\mathrm{Opt}^!$ (see section 3.1.1 on page 15) and OCMTs from $\mathrm{OptSet}^!$ (see section 3.5.1 on page 29). Specifically, we use two types of judgments.

(1) $E \triangleright \Gamma$, stating that $E \triangleright \Gamma$ is a well formed OCMT.
(2) $E \triangleright \Gamma \vdash t : T$, stating that in OCMT $E \triangleright \Gamma$, the term $t$ is well formed, and has type $T$. We may also write $E \triangleright \Gamma \vdash_n t : T$ if $t \in \mathbb{T}_n$.

3.9.2. **Inference rules.** We present the inference rules of system $\mathrm{OptSet}^!_{\mathrm{M}}$ in figure 3.9.1 on the following page. It uses rules `point`, `degen`, and `graft` from system $\mathrm{Opt}^!$, and rules `zero`, `sum`, and `glue` from system $\mathrm{OptSet}^!$. Rule `pd` is a counterpart to `degen`, introducing a non degenerate pasting diagram from an OCMT. Lastly, rule `shift` is a variant of that of system $\mathrm{Opt}^!$, and introduces a new cell from a pasting diagram, along with all its targets. It can be viewed as a fusion of $\mathrm{Opt}^!$'s `shift` rule and $\mathrm{OptSet}^!$'s `repr` rule.

FIGURE 3.9.1. The $\text{OptSet}^!_{\mathrm{M}}$ system.

**Introduction of points:** This rule introduces 0-cells, also called points. If $x \in \mathbb{V}_0$, then

$$\frac{}{\triangleright x : \varnothing}\ \texttt{point}$$

**Introduction of degenerate pasting diagrams:** This rule creates a new degenerate pasting diagram. If $x \in \mathbb{V}_{\Gamma,k}$, then

$$\frac{E \triangleright \Gamma, x : X}{E \triangleright \Gamma, x : X \vdash_{k+1} \underline{x} : x \multimap X}\ \texttt{degen}$$

**Introduction of non degenerate pasting diagrams:** This rule creates a new non-degenerate pasting diagram consisting of a single cell. It can then be extended using the `graft` rule. If $x \in \mathbb{V}_{\Gamma,k}$, then

$$\frac{E \triangleright \Gamma, x : X}{E \triangleright \Gamma, x : X \vdash_k x : X}\ \texttt{pd}$$

**Grafting:** This rule extends a previously derived non degenerate pasting diagram by grafting a cell. With the same conditions as rule `graft` of system $\text{Opt}^!$ (see section 3.1 on page 15): for $x \in \mathbb{V}_n$, $t \in \mathbb{T}_n$ is not degenerate, $a \in (\mathsf{s}\, t)^\bullet$ such that $\mathsf{s}\, a = \mathsf{s}\,\mathsf{s}\, x$, then

$$\frac{E \triangleright \Gamma \vdash_n t : s_1 \multimap s_2 \multimap \cdots \qquad F \triangleright \Upsilon \vdash_n x : X}{G \triangleright \Gamma \cup \Upsilon \vdash_n t(a \leftarrow x) : s_1[\mathsf{s}\, x/a] \multimap s_2 \multimap \cdots}\ \texttt{graft}$$

where $G$ is the union of $E$, $F$, and potentially a set of additional equalities incurred by the substitution $s_1[\mathsf{s}\, x/a]$. We also write `graft-`$a$ to make explicit that we grafted onto $a$.

**Shifting of pasting diagrams:** This rule takes a previously derived pasting diagram (degenerate or not), and introduces a new cell having this pasting diagram as source. It also introduces the targets of all its iterated sources, and extends the ambient equational theory with the required identities. If $x \in \mathbb{V}_{n+1}$ is such that $x \notin \mathbb{V}_\Gamma$, then

$$\frac{E \triangleright \Gamma \vdash_n t : T}{F \triangleright \Upsilon}\ \texttt{shift}$$

where

$$\Upsilon := \Gamma \cup \{x : t \multimap T\} \cup \left\{\mathsf{t}^i\, x : \mathsf{s}^{i+1}\, x \multimap \mathsf{s}^{i+2}\, x \multimap \cdots \mid 0 < i \leq n\right\}$$

by convention, we let $\mathsf{t}^0\, x = x$, and where $F$ is defined as follows:

(1) if $t$ is a degenerate term, say $t = \underline{a}$, then

$$F := E \cup \left\{\mathsf{t}^{i+2}\, x = \mathsf{t}^i\, a \mid 0 \leq i \leq n-1\right\} \tag{3.9.1}$$

(2) if $t$ is not degenerate, say $t = y(\overrightarrow{z_i \leftarrow u_i})$, for some $y \in \mathbb{V}_n$, $\overrightarrow{z_i} \in \mathbb{V}_{n-1}$, and $\overrightarrow{u_i} \in \mathbb{T}_n$, then

$$F := E \cup \left\{\mathsf{t}^2\, x = \mathsf{t}\, y \mid \text{if } n \geq 1\right\} \cup \{\mathsf{t}\, a = b \mid \text{for all } b \leftarrow a(\cdots) \text{ occurring in } t\}.$$

**Zero:** This rules introduces the empty OCMT.

$$\frac{}{\triangleright}\ \texttt{zero}$$

**Binary sums:** This rule takes two disjoint OCMTs (i.e. whose cells have different names), and produces their sum. If $\Gamma \cap \Upsilon = \varnothing$, then

$$\frac{E \triangleright \Gamma \qquad F \triangleright \Upsilon}{E \cup F \triangleright \Gamma \cup \Upsilon}\ \texttt{sum}$$

**Quotients:** This rule identifies two parallel cells in an opetopic set by extending the underlying equational theory. If $a, b \in \mathbb{V}_\Gamma$ are such that $\mathsf{s}\, a =_E \mathsf{s}\, b$ and $\mathsf{t}\, a =_E \mathsf{t}\, b$, then

$$\frac{E \triangleright \Gamma}{E \cup \{a = b\} \triangleright \Gamma}\ \texttt{glue}$$

We also write `glue-`$(a{=}b)$ to make explicit that we added $\{a = b\}$ to the theory.

**Remark 3.9.2.** The *sum* and *zero* rules may be replaced by the following *usum* rule (unbiased sum) without changing the set of derivable OCMTs: for $k \geq 0$, $(E_1 \triangleright \Gamma_1), \ldots, (E_k \triangleright \Gamma_k)$ OCMTs such that $\Gamma_i \cap \Gamma_j = \varnothing$ for all $i \neq j$, then

$$\frac{E_1 \triangleright \Gamma_1 \qquad \cdots \qquad E_k \triangleright \Gamma_k}{\sum_{i=1}^{k} E_i \triangleright \sum_{i=1}^{k} \Gamma_i} \; \mathit{usum}$$

**Remark 3.9.3.** Akin to $\mathrm{Opt}^!$ and $\mathrm{OptSet}^!$, in $\mathrm{OptSet}^!_{\mathrm{M}}$ a sequent or an OCMT that is equivalent to a derivable one is itself derivable.

## 3.10. Equivalence with opetopic sets

The aim of this section is to prove theorem 3.10.4 on page 44, stating that system $\mathrm{OptSet}^!_{\mathrm{M}}$ precisely derives opetopic sets, in the sense of theorems 3.6.15 and 3.6.17 on page 35 and on page 36. In other words, we prove that the set of derivable OCMTs of systems $\mathrm{OptSet}^!_{\mathrm{M}}$ and $\mathrm{OptSet}^!$ are the same. This is done by rewriting proof trees in $\mathrm{OptSet}^!$ to proof trees in $\mathrm{OptSet}^!_{\mathrm{M}}$ (see proposition 3.10.1) and conversely (see proposition 3.10.3 on page 44).

Throughout this section, the rules of systems $\mathrm{Opt}^!$ and $\mathrm{OptSet}^!$ will be decorated by a prime, e.g. `shift'`, in order to differentiate them from the rules of system $\mathrm{OptSet}^!_{\mathrm{M}}$. Further, to make notations lighter and the demonstrations more graphical, we write proof trees as actual trees, whose nodes are decorated by rules, and edges by sequents or OCMTs. For instance, derivation of the arrow ■ (see example 3.3.1 on page 26) in system $\mathrm{Opt}^!$ is represented as on the left, and more concisely as on the right:

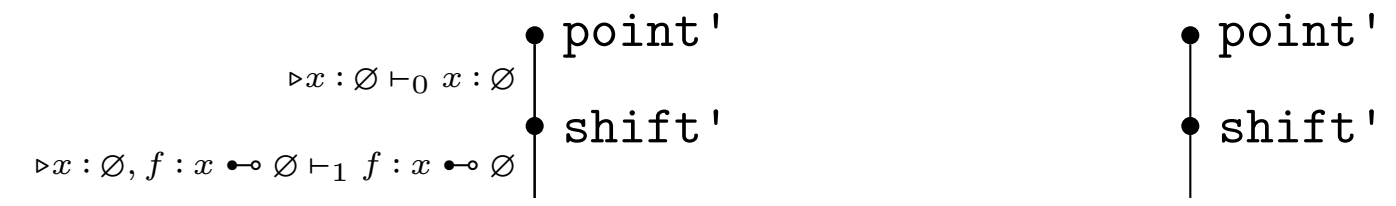

If no uncertainty arise, we leave the decoration of the edges implicit, as on the right.

**Proposition 3.10.1.** *Every OCMT derivable in system* $\mathrm{OptSet}^!$ *is also derivable in system* $\mathrm{OptSet}^!_{\mathrm{M}}$.

*Proof.* Recall that a proof tree in system $\mathrm{OptSet}^!$ has the following structure:

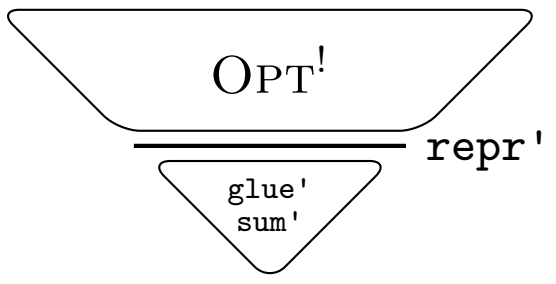

meaning that it begins with derivations in system $\mathrm{Opt}^!$, followed by instances of the `repr'` rule, followed by a derivation in system $\mathrm{OptSet}^!$. Remark that rule `glue'` is exactly `glue`, and likewise for `sum`, so that the bottom part of the proof tree is already a derivation in system $\mathrm{OptSet}^!_{\mathrm{M}}$.

We now show that we can rewrite the top part to a proof in system $\mathrm{OptSet}^!_{\mathrm{M}}$ by "lifting" the instances of `repr`, and replacing the other rule instances by those of $\mathrm{OptSet}^!_{\mathrm{M}}$. This rewriting procedure is defined by the following cases.

(1) If we have a proof tree as on the left, we rewrite it as on the right:

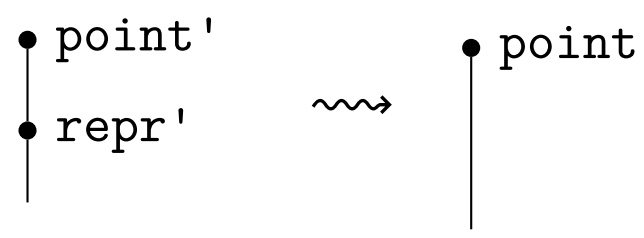

(2) Likewise, if we have a derivation as on the left, where $\Pi$ is a proof tree in system $\textsc{Opt}^!$ or $\textsc{OptSet}^!$, we rewrite it as on the right:

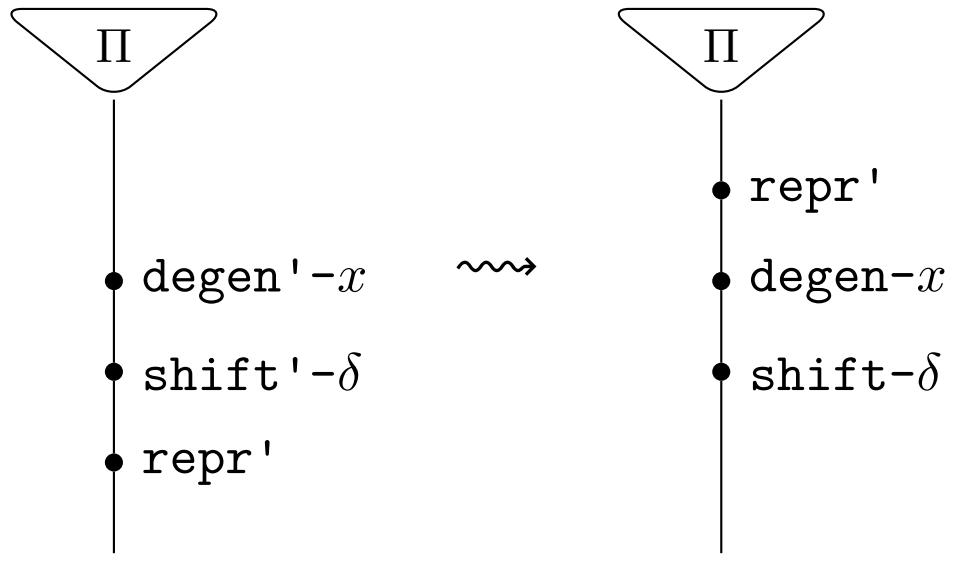

(3) Likewise, if we have a proof tree as on the left of figure 3.10.1, we rewrite it as on the right.

FIGURE 3.10.1. Proof of proposition 3.10.1 on the preceding page: commuting `repr'` with a `graft'` and `shift'` sequence

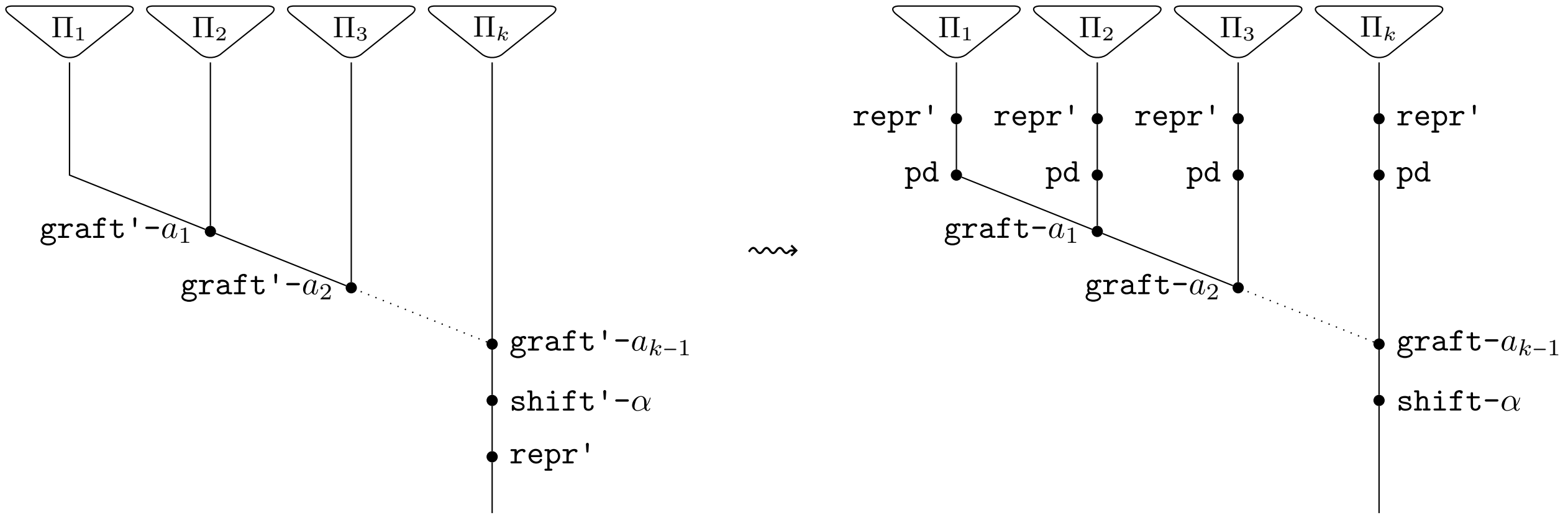

It is routine verification to check that the conclusion OCMT on the left and the right of any of those cases are the same. This rewriting procedure is convergent (i.e. confluent and terminating), and the normal form of a proof tree in system $\textsc{OptSet}^!$ is a proof tree in system $\textsc{OptSet}^!_{\mathrm{M}}$ that derives the same OCMT. □

**Lemma 3.10.2.** *Let $(E \triangleright \Gamma)$ be a derivable OCMT in system $\textsc{OptSet}^!_{\mathrm{M}}$. Then it admits a proof tree of the following form*

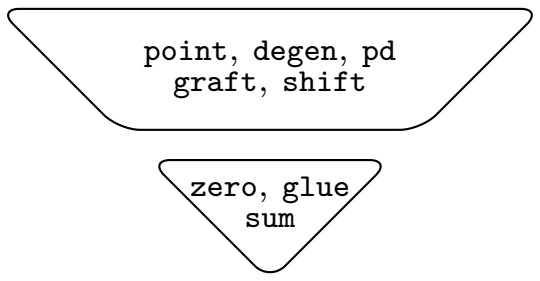

*meaning a proof tree starting with derivation in the fragment of system $\textsc{OptSet}^!_{\mathrm{M}}$ containing only rules* `point`*,* `degen`*,* `pd`*,* `graft`*, and* `shift`*, followed by a derivation in the complementary fragment.*

*Proof.* If we have a proof tree consisting only of an instance of rule `zero`, then the result trivially holds. Otherwise, we proceed by stating rewriting steps of proof trees in system $\textsc{OptSet}^!_{\mathrm{M}}$, as in the proof of proposition 3.10.1 on the previous page.

(1) If we have a proof tree as on the left, we rewrite it as on the right:

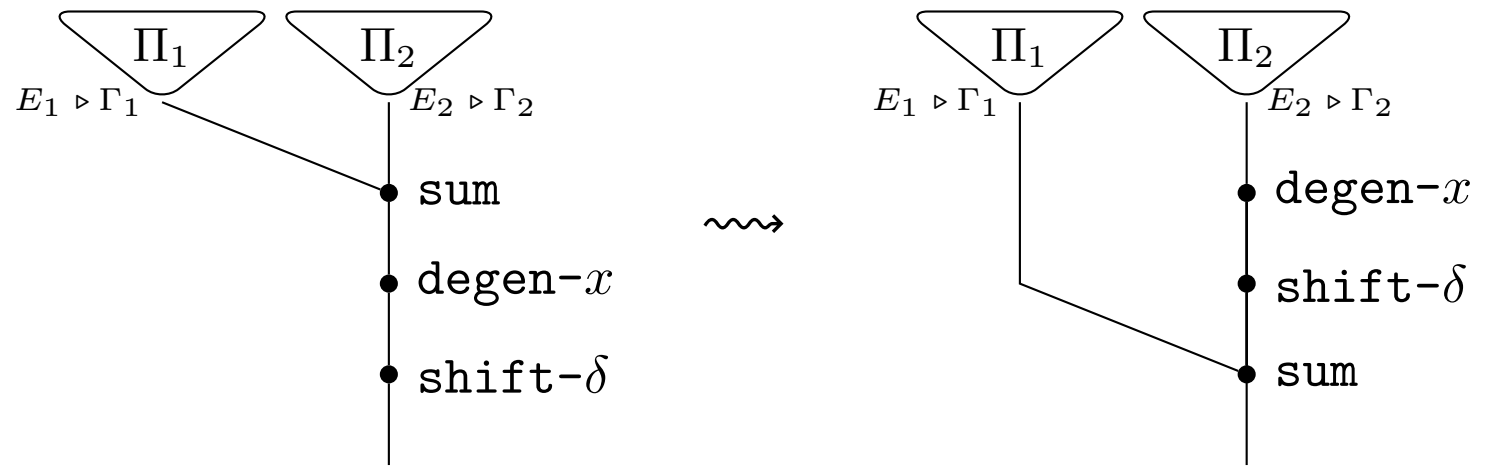

(2) Likewise, consider a proof tree as on the left. Then, by assumption on rule `sum`, either $x \in \mathbb{V}_{\Gamma_{i,1}}$ or $x \in \mathbb{V}_{\Gamma_{i,2}}$. Without loss of generality, assume the latter holds. Then we rewrite the proof tree as on the right:

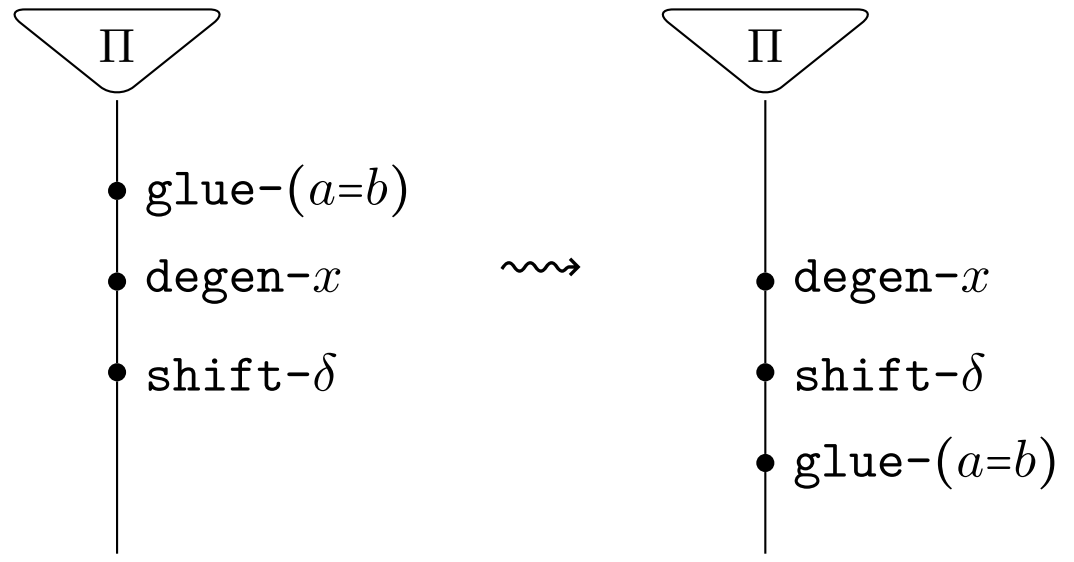

(3) Likewise, consider a proof tree as on the left of figure 3.10.2. Then, by assumption on rule `sum`, either $a_{i-1} \in \mathbb{V}_{\Gamma_1}$ or $a_{i-1} \in \mathbb{V}_{\Gamma_2}$. Without loss of generality, assume the latter holds. Then we rewrite the proof tree as on the right of the figure.

FIGURE 3.10.2. Proof of lemma 3.10.2 on the facing page: commuting `sum` with a `graft` and `shift` sequence

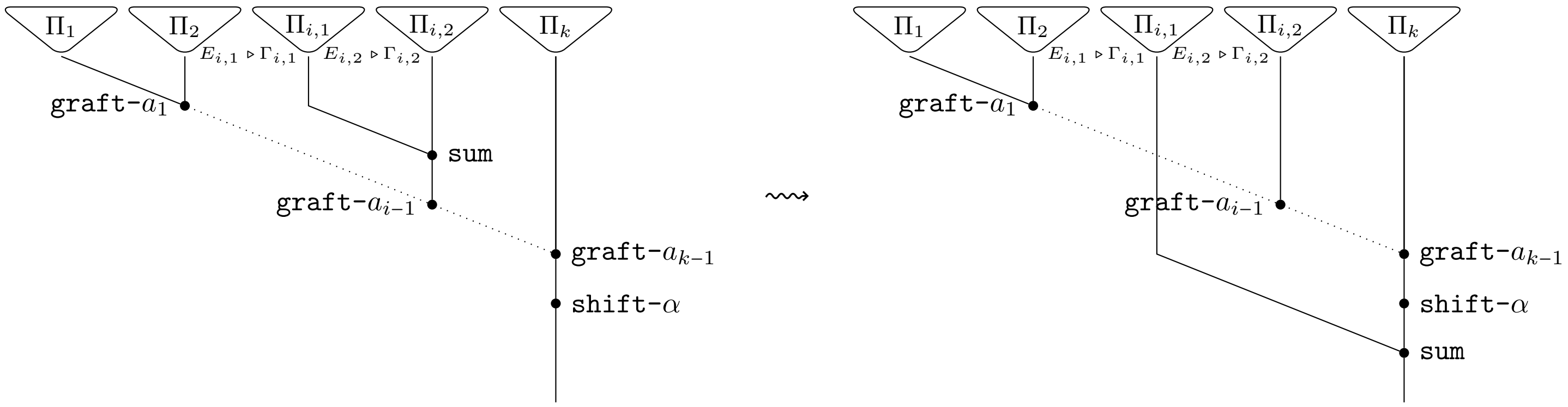

(4) Likewise, consider a proof tree as on the left of figure 3.10.3 on the next page. Then for $1 \leq j \leq k$, we have $x \in \mathbb{V}_{\Gamma_j}$ if and only if $y \in \mathbb{V}_{\Gamma_j}$. Assume $x, y \in \mathbb{V}_{\Gamma_j}$, and by induction, assume that the last rule instance of $\Pi_j$ is `glue-`$(x{=}y)$. Then we rewrite the proof tree as on the right of the figure, where $\Pi'_j$ is $\Pi_j$ with the instance of `glue-`$(x{=}y)$ removed, if $x, y \in \mathbb{V}_{\Gamma_j}$, or otherwise $\Pi_j$.

FIGURE 3.10.3. Proof of lemma 3.10.2 on page 42: commuting `glue` with a `graft` and `shift` sequence

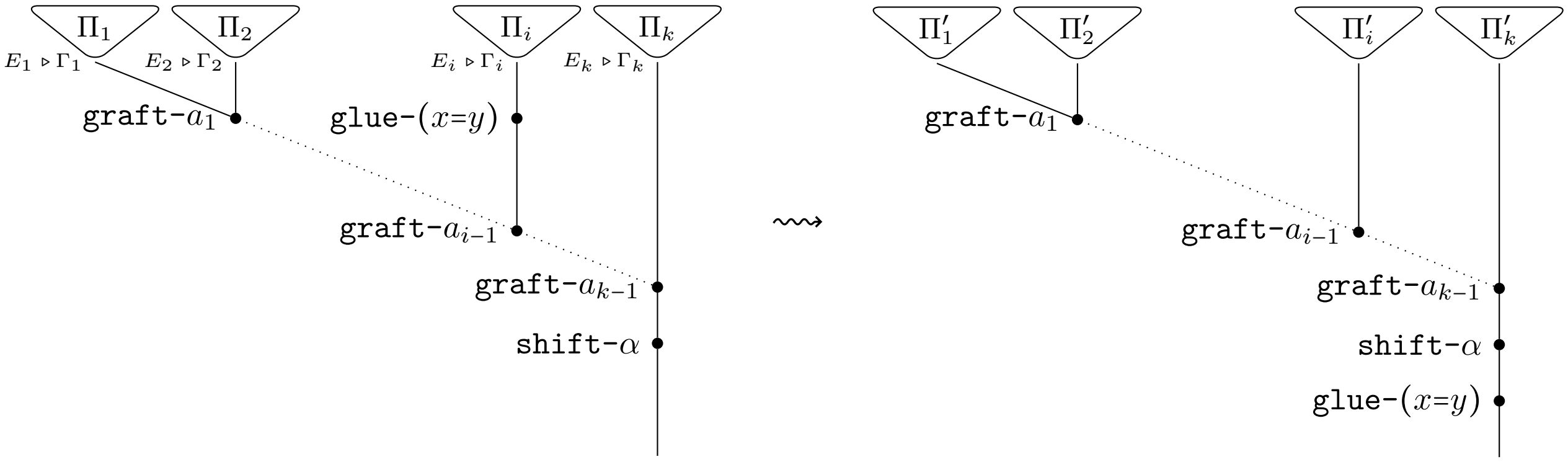

□

**Proposition 3.10.3.** *Every OCMT derivable in system* $\mathrm{OptSet}^!_{\mathrm{M}}$ *is also derivable in system* $\mathrm{OptSet}^!$.

*Proof.* Consider a proof tree in system $\mathrm{OptSet}^!_{\mathrm{M}}$. Then, without loss of generality, it has the structure described in lemma 3.10.2 on page 42. Applying the rewriting steps of proposition 3.10.1 on page 41 in reverse direction yields a proof tree in systems $\mathrm{Opt}^!$ and $\mathrm{OptSet}^!$ that derives the same OCMT. □

**Theorem 3.10.4.** *System* $\mathrm{OptSet}^!_{\mathrm{M}}$ *derives opetopic sets, in the sense of theorems 3.6.15 and 3.6.17 on page 35 and on page 36.*

*Proof.* By propositions 3.10.1 and 3.10.3 on page 41 and on the current page, the OCMTs derived by system $\mathrm{OptSet}^!_{\mathrm{M}}$ and $\mathrm{OptSet}^!$ are the same. □

## 3.11. Examples

In this section, we give example derivations in system $\mathrm{OptSet}^!_{\mathrm{M}}$. For clarity, we do not repeat the type of previously typed variables in proof trees.

**Example 3.11.1.** The opetopic set

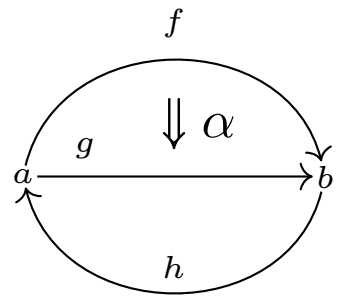

of example 3.7.1 on page 36 can be derived as follows. The first half of the proof tree is on the left, and the second half on the right. Moreover, for clarity, we do not repeat the typing of previously typed variables

$$
\dfrac{
\dfrac{
\dfrac{
\dfrac{
\dfrac{
\dfrac{}{\rhd\ a : \varnothing}\ \texttt{point}
}{\rhd\ a \vdash_0 a}\ \texttt{pd}
}{\rhd\ \begin{matrix} a, \mathsf{t} f : \varnothing \\ f : a \multimap \varnothing \end{matrix}}\ \texttt{shift}
}{\rhd\ \begin{matrix} a, \mathsf{t} f \\ f \end{matrix} \vdash_0 a}\ \texttt{pd}
}{\rhd\ \begin{matrix} a, \mathsf{t} f \\ f, g : a \multimap \varnothing \end{matrix}}\ \texttt{shift}
\qquad
\dfrac{
\dfrac{
\dfrac{}{\rhd\ b : \varnothing}\ \texttt{point}
}{\rhd\ b \vdash_0 b}\ \texttt{pd}
}{\rhd\ \begin{matrix} b, \mathsf{t} h : \varnothing \\ h : b \multimap \varnothing \end{matrix}}\ \texttt{shift}
}{
\dfrac{
\dfrac{\rhd\ \begin{matrix} a, b, \mathsf{t} f, \mathsf{t} g, \mathsf{t} h \\ f, g, h \end{matrix}}{b = \mathsf{t} f\ \rhd\ \begin{matrix} a, b, \mathsf{t} f, \mathsf{t} g, \mathsf{t} h \\ f, g, h \end{matrix}}\ (b{=}\mathsf{t} f)
}{b = \mathsf{t} f = \mathsf{t} g\ \rhd\ \begin{matrix} a, b, \mathsf{t} f, \mathsf{t} g, \mathsf{t} h \\ f, g, h \end{matrix}}\ (b{=}\mathsf{t} g)
}\ \texttt{sum}
$$

$$
\dfrac{
\dfrac{
\dfrac{
\dfrac{\vdots}{b = \mathsf{t} f = \mathsf{t} g\ \rhd\ \begin{matrix} a, b, \mathsf{t} f, \mathsf{t} g, \mathsf{t} h \\ f, g, h \end{matrix} \vdash_1 f}\ \texttt{pd}
}{b = \mathsf{t} f = \mathsf{t} g = \mathsf{t}^2 \alpha\ \rhd\ \begin{matrix} a, b, \mathsf{t} f, \mathsf{t} g, \mathsf{t} h, \mathsf{t}^2 \alpha : \varnothing \\ f, g, h, \mathsf{t} \alpha : a \multimap \varnothing \\ \alpha : f \multimap a \multimap \varnothing \end{matrix}}\ \texttt{shift}
}{\begin{matrix} b = \mathsf{t} f = \mathsf{t} g = \mathsf{t}^2 \alpha \\ g = \mathsf{t} \alpha \end{matrix}\ \rhd\ \begin{matrix} a, b, \mathsf{t} f, \mathsf{t} g, \mathsf{t} h, \mathsf{t}^2 \alpha : \varnothing \\ f, g, h, \mathsf{t} \alpha : a \multimap \varnothing \\ \alpha : f \multimap a \multimap \varnothing \end{matrix}}\ (g{=}\mathsf{t}\alpha)
}{\begin{matrix} b = \mathsf{t} f = \mathsf{t} g = \mathsf{t}^2 \alpha, a = \mathsf{t} h \\ g = \mathsf{t} \alpha \end{matrix}\ \rhd\ \begin{matrix} a, b, \mathsf{t} f, \mathsf{t} g, \mathsf{t} h, \mathsf{t}^2 \alpha : \varnothing \\ f, g, h, \mathsf{t} \alpha : a \multimap \varnothing \\ \alpha : f \multimap a \multimap \varnothing \end{matrix}}\ (a{=}\mathsf{t} h)
$$

where rules of the form $(x = y)$ are shorthands for $\texttt{glue-}(x = y)$.

**Example 3.11.2.** The opetopic set

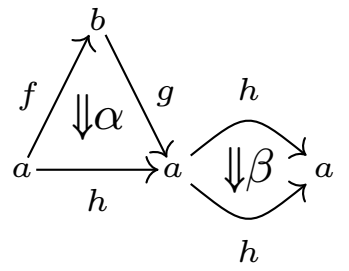

of example 3.7.2 on page 37 can be derived as follows.

$$\dfrac{\dfrac{}{\triangleright\ a : \varnothing}\ \texttt{point} \qquad \dfrac{}{\triangleright\ b : \varnothing}\ \texttt{point}}{\triangleright\ a, b}\ \texttt{sum}$$

$$\dfrac{\triangleright\ a, b}{\triangleright\ a, b \vdash_0 a : \varnothing}\ \texttt{pd}$$

$$\dfrac{\triangleright\ a, b \vdash_0 a : \varnothing}{\triangleright\ \begin{array}{l} a, b, \mathsf{t} f : \varnothing \\ f : a \multimap \varnothing \end{array}}\ \texttt{shift}$$

$$\dfrac{\triangleright\ \begin{array}{l} a, b, \mathsf{t} f : \varnothing \\ f : a \multimap \varnothing \end{array}}{b = \mathsf{t} f\ \triangleright\ \begin{array}{l} a, b, \mathsf{t} f \\ f \end{array}}\ \texttt{glue-}(b = \mathsf{t} f)$$

$$\dfrac{b = \mathsf{t} f\ \triangleright\ \begin{array}{l} a, b, \mathsf{t} f \\ f \end{array}}{b = \mathsf{t} f\ \triangleright\ \begin{array}{l} a, b, \mathsf{t} f \\ f \end{array} \vdash_0 \mathsf{t} f : \varnothing}\ \texttt{pd}$$

$$\dfrac{b = \mathsf{t} f\ \triangleright\ \begin{array}{l} a, b, \mathsf{t} f \\ f \end{array} \vdash_0 \mathsf{t} f : \varnothing}{b = \mathsf{t} f\ \triangleright\ \begin{array}{l} a, b, \mathsf{t} f, \mathsf{t} g : \varnothing \\ f, g : \mathsf{t} f \multimap \varnothing \end{array}}\ \texttt{shift}$$

$$\dfrac{b = \mathsf{t} f\ \triangleright\ \begin{array}{l} a, b, \mathsf{t} f, \mathsf{t} g : \varnothing \\ f, g : \mathsf{t} f \multimap \varnothing \end{array}}{b = \mathsf{t} f\ \triangleright\ \begin{array}{l} a, b, \mathsf{t} f, \mathsf{t} g \\ f, g \end{array} \vdash_0 a : \varnothing}\ \texttt{pd}$$

$$\dfrac{b = \mathsf{t} f\ \triangleright\ \begin{array}{l} a, b, \mathsf{t} f, \mathsf{t} g \\ f, g \end{array} \vdash_0 a : \varnothing}{b = \mathsf{t} f\ \triangleright\ \begin{array}{l} a, b, \mathsf{t} f, \mathsf{t} g, \mathsf{t} h : \varnothing \\ f, g, h : a \multimap \varnothing \end{array}}\ \texttt{shift}$$

$$\dfrac{b = \mathsf{t} f\ \triangleright\ \begin{array}{l} a, b, \mathsf{t} f, \mathsf{t} g, \mathsf{t} h : \varnothing \\ f, g, h : a \multimap \varnothing \end{array}}{b = \mathsf{t} f\ \triangleright\ \begin{array}{l} a, b, \mathsf{t} f, \mathsf{t} g, \mathsf{t} h \\ f, g, h \end{array} \vdash_1 g : \mathsf{t} f \multimap \varnothing}\ \texttt{pd}$$

$$\dfrac{b = \mathsf{t} f\ \triangleright\ \begin{array}{l} a, b, \mathsf{t} f, \mathsf{t} g, \mathsf{t} h \\ f, g, h \end{array} \vdash_1 g : \mathsf{t} f \multimap \varnothing}{b = \mathsf{t} f\ \triangleright\ \begin{array}{l} a, b, \mathsf{t} f, \mathsf{t} g, \mathsf{t} h \\ f, g, h \end{array} \vdash_1 g(\mathsf{t} f \leftarrow f) : a \multimap \varnothing}\ \texttt{graft}$$

$$\dfrac{b = \mathsf{t} f\ \triangleright\ \begin{array}{l} a, b, \mathsf{t} f, \mathsf{t} g, \mathsf{t} h \\ f, g, h \end{array} \vdash_1 g(\mathsf{t} f \leftarrow f) : a \multimap \varnothing}{b = \mathsf{t} f = \mathsf{t}\mathsf{t}\alpha\ \triangleright\ \begin{array}{l} a, b, \mathsf{t} f, \mathsf{t} g, \mathsf{t} h, \mathsf{t}\mathsf{t}\alpha : \varnothing \\ f, g, h, \mathsf{t}\alpha : a \multimap \varnothing \\ \alpha : g(\mathsf{t} f \leftarrow f) : a \multimap \varnothing \end{array}}\ \texttt{shift}$$

$$\dfrac{b = \mathsf{t} f = \mathsf{t}\mathsf{t}\alpha\ \triangleright\ \begin{array}{l} a, b, \mathsf{t} f, \mathsf{t} g, \mathsf{t} h, \mathsf{t}\mathsf{t}\alpha : \varnothing \\ f, g, h, \mathsf{t}\alpha : a \multimap \varnothing \\ \alpha : g(\mathsf{t} f \leftarrow f) : a \multimap \varnothing \end{array}}{\begin{array}{l} b = \mathsf{t} f = \mathsf{t}\mathsf{t}\alpha \\ g = \mathsf{t}\alpha \end{array}\ \triangleright\ \begin{array}{l} a, b, \mathsf{t} f, \mathsf{t} g, \mathsf{t} h, \mathsf{t}\mathsf{t}\alpha \\ f, g, h, \mathsf{t}\alpha \\ \alpha \end{array}}\ \texttt{glue-}(h = \mathsf{t}\alpha)$$

$$\dfrac{\begin{array}{l} b = \mathsf{t} f = \mathsf{t}\mathsf{t}\alpha \\ g = \mathsf{t}\alpha \end{array}\ \triangleright\ \begin{array}{l} a, b, \mathsf{t} f, \mathsf{t} g, \mathsf{t} h, \mathsf{t}\mathsf{t}\alpha \\ f, g, h, \mathsf{t}\alpha \\ \alpha \end{array}}{\begin{array}{l} b = \mathsf{t} f = \mathsf{t}\mathsf{t}\alpha, a = \mathsf{t} g \\ g = \mathsf{t}\alpha \end{array}\ \triangleright\ \begin{array}{l} a, b, \mathsf{t} f, \mathsf{t} g, \mathsf{t} h, \mathsf{t}\mathsf{t}\alpha \\ f, g, h, \mathsf{t}\alpha \\ \alpha \end{array}}\ \texttt{glue-}(a = \mathsf{t} g)$$

$$\dfrac{\begin{array}{l} b = \mathsf{t} f = \mathsf{t}\mathsf{t}\alpha, a = \mathsf{t} g \\ g = \mathsf{t}\alpha \end{array}\ \triangleright\ \begin{array}{l} a, b, \mathsf{t} f, \mathsf{t} g, \mathsf{t} h, \mathsf{t}\mathsf{t}\alpha \\ f, g, h, \mathsf{t}\alpha \\ \alpha \end{array}}{\begin{array}{l} b = \mathsf{t} f = \mathsf{t}\mathsf{t}\alpha, a = \mathsf{t} g = \mathsf{t} h \\ g = \mathsf{t}\alpha \end{array}\ \triangleright\ \begin{array}{l} a, b, \mathsf{t} f, \mathsf{t} g, \mathsf{t} h, \mathsf{t}\mathsf{t}\alpha \\ f, g, h, \mathsf{t}\alpha \\ \alpha \end{array}}\ \texttt{glue-}(a = \mathsf{t} h)$$

$$\dfrac{\begin{array}{l} b = \mathsf{t} f = \mathsf{t}\mathsf{t}\alpha, a = \mathsf{t} g = \mathsf{t} h \\ g = \mathsf{t}\alpha \end{array}\ \triangleright\ \begin{array}{l} a, b, \mathsf{t} f, \mathsf{t} g, \mathsf{t} h, \mathsf{t}\mathsf{t}\alpha \\ f, g, h, \mathsf{t}\alpha \\ \alpha \end{array}}{\begin{array}{l} b = \mathsf{t} f = \mathsf{t}\mathsf{t}\alpha, a = \mathsf{t} g = \mathsf{t} h \\ g = \mathsf{t}\alpha \end{array}\ \triangleright\ \begin{array}{l} a, b, \mathsf{t} f, \mathsf{t} g, \mathsf{t} h, \mathsf{t}\mathsf{t}\alpha \\ f, g, h, \mathsf{t}\alpha \\ \alpha \end{array} \vdash_1 h : a \multimap \varnothing}\ \texttt{pd}$$

$$\dfrac{\begin{array}{l} b = \mathsf{t} f = \mathsf{t}\mathsf{t}\alpha, a = \mathsf{t} g = \mathsf{t} h \\ g = \mathsf{t}\alpha \end{array}\ \triangleright\ \begin{array}{l} a, b, \mathsf{t} f, \mathsf{t} g, \mathsf{t} h, \mathsf{t}\mathsf{t}\alpha \\ f, g, h, \mathsf{t}\alpha \\ \alpha \end{array} \vdash_1 h : a \multimap \varnothing}{\begin{array}{l} b = \mathsf{t} f = \mathsf{t}\mathsf{t}\alpha, a = \mathsf{t} g = \mathsf{t} h = \mathsf{t}\mathsf{t}\beta \\ g = \mathsf{t}\alpha \end{array}\ \triangleright\ \begin{array}{l} a, b, \mathsf{t} f, \mathsf{t} g, \mathsf{t} h, \mathsf{t}\mathsf{t}\alpha, \mathsf{t}\mathsf{t}\beta : \varnothing \\ f, g, h, \mathsf{t}\alpha, \mathsf{t}\beta : a \multimap \varnothing \\ \alpha, \beta : h \multimap a \multimap \varnothing \end{array}}\ \texttt{shift}$$

$$\dfrac{\begin{array}{l} b = \mathsf{t} f = \mathsf{t}\mathsf{t}\alpha, a = \mathsf{t} g = \mathsf{t} h = \mathsf{t}\mathsf{t}\beta \\ g = \mathsf{t}\alpha \end{array}\ \triangleright\ \begin{array}{l} a, b, \mathsf{t} f, \mathsf{t} g, \mathsf{t} h, \mathsf{t}\mathsf{t}\alpha, \mathsf{t}\mathsf{t}\beta : \varnothing \\ f, g, h, \mathsf{t}\alpha, \mathsf{t}\beta : a \multimap \varnothing \\ \alpha, \beta : h \multimap a \multimap \varnothing \end{array}}{\begin{array}{l} b = \mathsf{t} f = \mathsf{t}\mathsf{t}\alpha, a = \mathsf{t} g = \mathsf{t} h = \mathsf{t}\mathsf{t}\beta \\ g = \mathsf{t}\alpha, f = \mathsf{t}\beta = h \end{array}\ \triangleright\ \begin{array}{l} a, b, \mathsf{t} f, \mathsf{t} g, \mathsf{t} h, \mathsf{t}\mathsf{t}\alpha, \mathsf{t}\mathsf{t}\beta \\ f, g, h, \mathsf{t}\alpha, \mathsf{t}\beta \\ \alpha, \beta \end{array}}\ \texttt{glue-}(h = \mathsf{t}\beta)$$

## 3.12. Python implementation

The system $\mathrm{OptSet}^!_m$ is implemented in module `opetopy.NamedOpetopicSetM` of [**Ho Thanh, 2018b**]. Its usage is very similar to that of `opetopy.NamedOpetope` and `opetopy.NamedOpetopicSet`, presented in sections 3.4 and 3.8 on page 28 and on page 38 respectively.

Figure 3.12.1. Derivation of example example 3.7.1 on page 36 using `opetopy.NamedOpetopicSetM`

```python
from NamedOpetopicSetM import Glue, Pd, Point, RuleInstance, Shift, Sum

# We first derive the f, g components
p1 = Shift(Pd(Point("a"), "a"), "f")
p1 = Shift(Pd(p1, "a"), "g")

# We then derive the h component
p2 = Shift(Pd(Point("b"), "b"), "h")

# We proceed to sum the two, glue some cells, and introduce alpha
example = Sum(p1, p2)  # type: RuleInstance
example = Glue(example, "b", "tf")
example = Glue(example, "b", "tg")
example = Shift(Pd(example, "f"), "alpha")
example = Glue(example, "b", "ttalpha")
example = Glue(example, "g", "talpha")
example = Glue(example, "a", "th")
```

CHAPTER 4

# Unnamed approach

## 4.1. The system for opetopes

The unnamed approach for opetopes relies on the calculus of higher addresses presented in section 2.2.3 on page 11 to identify cells, rather than names as in the named approach. For example, recall the opetope $\mathbf{3} \in \mathbb{O}_2$ from example 2.2.1 on page 12, drawn on the left, with its underlying $\mathfrak{Z}^0$-tree represented on the right:

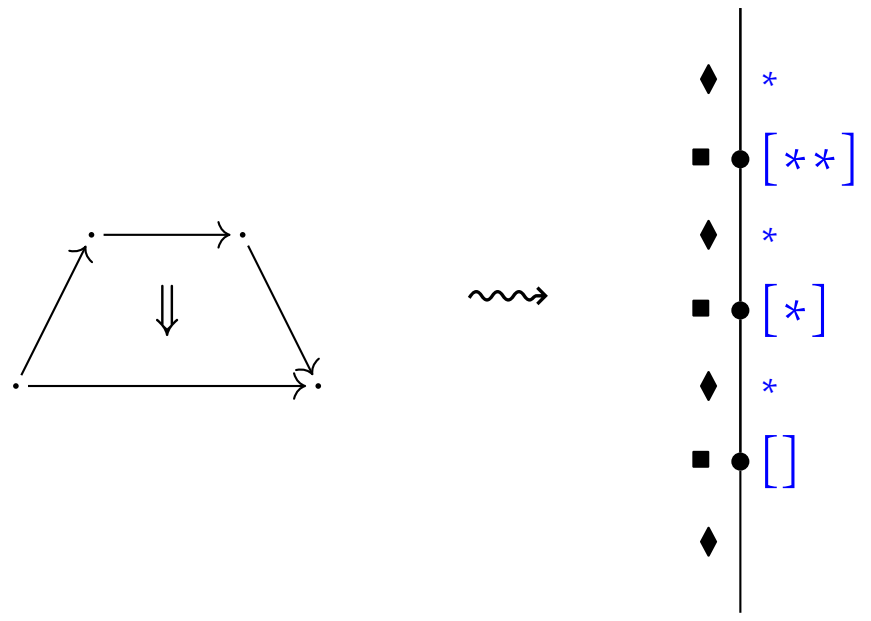

In the unnamed approach for opetopes presented in this section, $\mathbf{3}$ will be encoded as a mapping from its set of node addresses $\mathbf{3}^\bullet = \{[], [*], [**]\}$ to the set of 1-opetopes $\mathbb{O}_1$ as follows:

$$\mathbf{3} \qquad \rightsquigarrow \qquad \begin{cases} [] \leftarrow \blacksquare \\ [*] \leftarrow \blacksquare \\ [**] \leftarrow \blacksquare \end{cases}$$

The 1-opetope $\blacksquare$ can recursively be encoded by $\{* \leftarrow \blacklozenge$, which give a complete expression of $\mathbf{3}$:

$$\begin{cases} [] \leftarrow \{* \leftarrow \blacklozenge \\ [*] \leftarrow \{* \leftarrow \blacklozenge \\ [**] \leftarrow \{* \leftarrow \blacklozenge \end{cases}$$

This example will be treated in depth in example 4.3.2 on page 55. In a similar manner, the opetope $\omega$ on the left, whose tree is given in the middle, can be encoded as on the right:

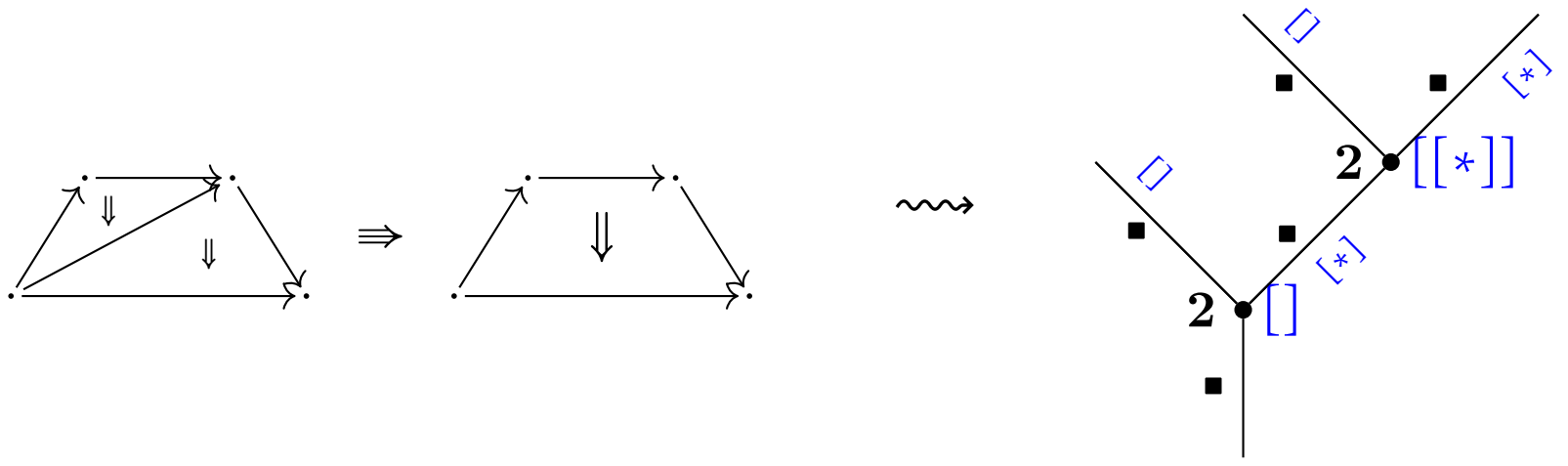

$$\rightsquigarrow \qquad \begin{cases} [] \leftarrow \mathbf{2} \\ [[*]] \leftarrow \mathbf{2} \end{cases}$$

Similarly to $\mathbf{3}$, the opetope $\mathbf{2}$ can be expressed by $\begin{cases} [] \leftarrow \blacksquare \\ [*] \leftarrow \blacksquare \end{cases}$, and recall that $\blacksquare$ can be expressed by $\{* \leftarrow \blacklozenge$. We thus have a complete encoding of $\omega$ as:

$$\omega \qquad \rightsquigarrow \qquad \begin{cases} [] \leftarrow \begin{cases} [] \leftarrow \{* \leftarrow \blacklozenge \\ [*] \leftarrow \{* \leftarrow \blacklozenge \end{cases} \\ [[*]] \leftarrow \begin{cases} [] \leftarrow \{* \leftarrow \blacklozenge \\ [*] \leftarrow \{* \leftarrow \blacklozenge \end{cases} \end{cases}$$

This example is fully treated in example 4.3.3 on page 55.

4.1.1. **Preopetopes.** Now we define the set $\mathbb{P}_n$ of $n$-preopetopes by the following grammar:

$$
\begin{array}{rcll}
\mathbb{P}_{-1} & ::= & \varnothing & \text{by convention} \\
\mathbb{P}_0 & ::= & \blacklozenge & \\
\mathbb{P}_n & ::= & \begin{cases} \mathbb{A}_{n-1} \leftarrow \mathbb{P}_{n-1} \\ \vdots \\ \mathbb{A}_{n-1} \leftarrow \mathbb{P}_{n-1} \end{cases} & n \geq 1 \qquad (4.1.1) \\
 & \mid & \{\!\{ \mathbb{P}_{n-2} & n \geq 2 \qquad (4.1.2)
\end{array}
$$

where the set $\mathbb{A}_n$ of $n$-addresses is defined in section 2.2.3 on page 11. In line (4.1.1), we require further that there is at least one $(n-1)$-address, and that all addresses are distinct. Further, a preopetope of this form is considered as a set of expressions $\mathbb{A}_{n-1} \leftarrow \mathbb{P}_{n-1}$ rather than a list: for instance, the following two $(n+1)$-preopetopes are equal

$$
\begin{cases} [p_1] \leftarrow \mathbf{p}_1 \\ [p_2] \leftarrow \mathbf{p}_2 \end{cases} = \begin{cases} [p_2] \leftarrow \mathbf{p}_2 \\ [p_1] \leftarrow \mathbf{p}_1 \end{cases}
$$

for any distinct $n$-addresses $[p_1], [p_2] \in \mathbb{A}_n$, and any $\mathbf{p}_1, \mathbf{p}_2 \in \mathbb{P}_n$. Here are examples of a 1, 2, 3, and 4-preopetope respectively:

$$
\{ * \leftarrow \blacklozenge \qquad \begin{cases} [] \leftarrow \{ * \leftarrow \blacklozenge \\ [*****] \leftarrow \{ * \leftarrow \blacklozenge \end{cases} \qquad \begin{cases} [[*]] \leftarrow \{\!\{ \blacklozenge \\ [[**][*][]] \leftarrow \{ [] \leftarrow \{ * \leftarrow \blacklozenge \end{cases} \qquad \{\!\{ \{\!\{ \blacklozenge
$$

As mentioned in the introduction of this chapter, the first corresponds to the unique 1-opetope $\blacksquare$. We will see that the second does not correspond to an actual opetope, as it is impossible for a 2-opetope to only contain addresses $[]$ and $[*****]$ (it would at least need addresses $[*]$, $[**]$, $[***]$, and $[****]$), and neither does the third, as it does not have a root node (corresponding to address $[]$). Finally, the last preopetope does represent an opetope.

An $n$-preopetope $\mathbf{p}$ is *degenerate* if it is of the form of line (4.1.2), it is *non-degenerate* otherwise; we write $\dim \mathbf{p} = n$ for its *dimension*. There is a unique 1-preopetope $\{ * \leftarrow \blacklozenge$, which we write $\blacksquare$. If we have a non-degenerate $n$-preopetope of the form

$$
\mathbf{p} = \begin{cases} [p_1] \leftarrow \mathbf{q}_1 \\ \vdots \\ [p_k] \leftarrow \mathbf{q}_k \end{cases} \tag{4.1.3}
$$

we call $[p_1], \ldots, [p_k]$ the *source addresses* of $\mathbf{p}$ (or just *sources*), we write $\mathbf{p}^\bullet := \{[p_1], \ldots, [p_k]\} \subseteq \mathbb{A}_{n-1}$ the set of addresses of $\mathbf{p}$, and $\mathsf{s}_{[p_i]} \mathbf{p} := \mathbf{q}_i$ the $[p_i]$-*source of* $\mathbf{p}$. Assume $n \geq 2$. A *leaf of* $\mathbf{p}$ is an $(n-1)$-address of the form $[p[q]]$ such that $[p] \in \mathbf{p}^\bullet$, $[q] \in (\mathsf{s}_{[p]} \mathbf{p})^\bullet$, and $[p[q]] \notin \mathbf{p}^\bullet$. We write $\mathbf{p}^| \subseteq \mathbb{A}_{n-1}$ for the set of leaves of $\mathbf{p}$. By convention, if $\mathbf{p}$ is degenerate, we set $\mathbf{p}^\bullet := \varnothing$ and $\mathbf{p}^| := \{[]\}$, and furthermore, $\blacklozenge^\bullet = \blacklozenge^| := \varnothing$.

**Definition 4.1.4** (Improper grafting)**.** Let $\mathbf{p}$ be an $n$-preopetope as in equation (4.1.3), and $\mathbf{q}$ be an $(n-1)$-preopetope. If $[r] \in \mathbf{p}^|$ is a leaf address of $\mathbf{p}$, write

$$
\mathbf{p} \underset{[r]}{\tilde{\circ}} \mathbf{q} := \begin{cases} [p_1] \leftarrow \mathbf{q}_1 \\ \vdots \\ [p_k] \leftarrow \mathbf{q}_k \\ [r] \leftarrow \mathbf{q} \end{cases}
$$

and call $\mathbf{p} \tilde{\circ}_{[r]} \mathbf{q}$ the improper grafting of $\mathbf{q}$ on $\mathbf{p}$ at address $[r]$. By convention, the grafting operation is associative on the right.

For example,

$$
\begin{cases} [] \leftarrow \begin{cases} [] \leftarrow \{ * \leftarrow \blacklozenge \\ [*] \leftarrow \{ * \leftarrow \blacklozenge \end{cases} \\ [[*]] \leftarrow \begin{cases} [] \leftarrow \{ * \leftarrow \blacklozenge \\ [*] \leftarrow \{ * \leftarrow \blacklozenge \end{cases} \end{cases} = \left\{ [] \leftarrow \begin{cases} [] \leftarrow \{ * \leftarrow \blacklozenge \\ [*] \leftarrow \{ * \leftarrow \blacklozenge \end{cases} \right. \underset{[[*]]}{\tilde{\circ}} \begin{cases} [] \leftarrow \{ * \leftarrow \blacklozenge \\ [*] \leftarrow \{ * \leftarrow \blacklozenge \end{cases}
$$

which, together with the introduction of this chapter, means that graphically,

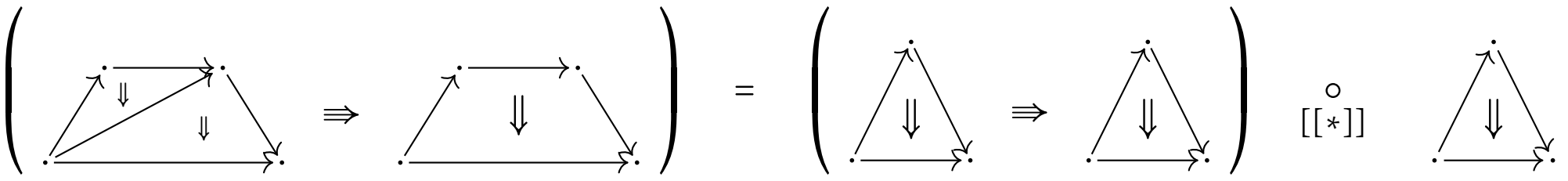

The denomination "improper" is motivated by the fact that $\mathbf{p}$ and $\mathbf{q}$ do not have the same dimension. A definition of "proper" grafting is given for polynomial trees in section 2.1 on page 9. Any preopetope can be obtained by iterated improper grafting as follows.

**Lemma 4.1.5.** *Let $\mathbf{p}$ be an $n$-preopetope as in equation* (4.1.3) *on the preceding page, and assume that whenever $1 \le i < j \le k$, we have either $[p_i] \sqsubseteq [p_j]$, or $[p_i]$ and $[p_j]$ are $\sqsubseteq$-incomparable (in particular, this condition is satisfied if the $[p_i]$'s are lexicographically sorted). Then*

$$\mathbf{p} := \left( \cdots \left( \left\{ [p_1] \leftarrow \mathbf{q}_1 \right) \underset{[p_2]}{\tilde{\circ}} \mathbf{q}_2 \cdots \right) \underset{[p_k]}{\tilde{\circ}} \mathbf{q}_k .$$

*Proof.* The condition on the sequence $[p_1], \ldots, [p_k]$ guarantees that the successive improper graftings are well-defined. ☐

4.1.2. **Inference rules.** We now introduce a typing system for preopetopes in order to characterize those corresponding to opetopes, which is formally shown in theorem 4.2.4 on page 54. We will deal with sequents of the form

$$\Gamma \vdash \mathbf{p} \longrightarrow \mathbf{t},$$

where $\mathbf{p}$ is an $n$-preopetope, $\mathbf{t}$ is an $(n-1)$-preopetope, and the context $\Gamma$ is a finite set of pairs consisting of addresses $[l] \in \mathbf{p}^{|}$ and $[q] \in \mathbf{t}^{\bullet}$, denoted by $\frac{[l]}{[q]}$. The preopetope $\mathbf{p}$ is the real object of interest, and as we will see in subsequent results, we may think of $\mathbf{t}$ as the "target" of $\mathbf{p}$, while $\Gamma$ establishes a bijection between the leaves of $\mathbf{p}$ and the nodes of its target.

The operation of substitution, which consists in replacing a node by a pasting scheme in an opetope, can be defined as follows in our formalism.

**Definition 4.1.6** (Substitution)**.** Let $\mathbf{t}$ be a non degenerate $n$-preopetope, and $\Upsilon \vdash \mathbf{q} \longrightarrow \mathbf{u}$ be a sequent, where $\mathbf{q}$ is an $n$-preopetope as well. Write $\mathbf{t}$ as

$$\mathbf{t} = \begin{cases} [t_1] \leftarrow \mathbf{w}_1 \\ \vdots \\ [t_l] \leftarrow \mathbf{w}_l \end{cases}$$

For $[t_i] \in \mathbf{t}^{\bullet}$, we define $\mathbf{t} \, \square_{[t_i]} \, \mathbf{q}$, the *substitution* by $\mathbf{q}$ in $\mathbf{t}$ at $[t_i]$, as follows:

(1) if $l = 1$ and $\mathbf{q}$ is degenerate, then $\mathbf{t} \, \square_{[t_1]} \, \mathbf{q} = \mathbf{q}$;
(2) if $l \ge 2$ and $\mathbf{q}$ is degenerate, then

$$\mathbf{t} \underset{[t_i]}{\square} \mathbf{q} := \begin{cases} \chi[t_1] \leftarrow \mathbf{w}_1 \\ \vdots \\ \chi[t_{i-1}] \leftarrow \mathbf{w}_{i-1} \\ \chi[t_{i+1}] \leftarrow \mathbf{w}_{i+1} \\ \vdots \\ \chi[t_l] \leftarrow \mathbf{w}_l \end{cases} \qquad \text{where} \qquad \chi[t_j] := \begin{cases} [t_i \cdots] & \text{if } [t_j] = [t_i[\,]\cdots] \\ [t_j] & \text{otherwise,} \end{cases}$$

(3) if $l \ge 2$, and $\mathbf{q}$ is not degenerate, write it as

$$\mathbf{q} = \begin{cases} [q_1] \leftarrow \mathbf{v}_1 \\ \vdots \\ [q_k] \leftarrow \mathbf{v}_k \end{cases}$$

and define

$$\mathbf{t} \underset{[t_i]}{\square} \mathbf{q} := \begin{cases} \chi[t_1] \leftarrow \mathbf{w}_1 \\ \vdots \\ \chi[t_{i-1}] \leftarrow \mathbf{w}_{i-1} \\ [t_i q_1] \leftarrow \mathbf{v}_1 \\ \vdots \\ [t_i q_k] \leftarrow \mathbf{v}_k \\ \chi[t_{i+1}] \leftarrow \mathbf{w}_{i+1} \\ \vdots \\ \chi[t_l] \leftarrow \mathbf{w}_l \end{cases} \qquad \text{where} \qquad \chi[t_j] := \begin{cases} [t_i a \cdots] & \text{if } [t_j] = [t_i[b]\cdots] \\ & \text{for some } \frac{[a]}{[b]} \in \Upsilon, \\ [t_j] & \text{otherwise.} \end{cases}$$

This operation relies on the context $\Upsilon$, which we leave implicit. By convention, $\square$ is associative on the left.

Refer to section 4.3 on page 54 for examples of applications of this construction. We now state the inference rules of our unnamed system Opt$^?$ in figure 4.1.1.

Figure 4.1.1. The Opt$^?$ system.

**Introduction of points:**

$$\frac{}{\vdash \blacklozenge \longrightarrow \varnothing}\ \texttt{point}$$

**Introduction of degeneracies:**

$$\frac{\Gamma \vdash \mathbf{p} \longrightarrow \mathbf{t}}{\frac{[]}{[]} \vdash \{\!\{\mathbf{p} \longrightarrow \{[] \leftarrow \mathbf{p}}\ \texttt{degen}$$

If $\dim \mathbf{p} = n$, remark that $\dim\left(\{\!\{\mathbf{p}\right) = n + 2$, and that $\dim\left(\{[] \leftarrow \mathbf{p}\right) = n + 1$.

**Shift to the next dimension:** Assume $\mathbf{p}^\bullet = \{[p_1], \ldots, [p_k]\}$.

$$\frac{\Gamma \vdash \mathbf{p} \longrightarrow \mathbf{t}}{\frac{[[p_1]]}{[p_1]}, \ldots, \frac{[[p_k]]}{[p_k]} \vdash \{[] \leftarrow \mathbf{p} \longrightarrow \mathbf{p}}\ \texttt{shift}$$

As in the previous rule, if $\dim \mathbf{p} = n$, then $\dim\left(\{[] \leftarrow \mathbf{p}\right) = n + 1$.

**Grafting:** Assume $\dim \mathbf{p} = n \geq 2$, $\dim \mathbf{q} = n - 1$, $\mathbf{q}^\bullet = \{[s_1], \ldots, [s_l]\}$, and $\mathsf{s}_{[q]}\, \mathsf{s}_{[p]}\, \mathbf{p} = \mathbf{u}$.

$$\frac{\Gamma, \frac{[p[q]]}{[r]} \vdash \mathbf{p} \longrightarrow \mathbf{t} \qquad \Upsilon \vdash \mathbf{q} \longrightarrow \mathbf{u}}{\Gamma', \frac{[p[q][s_1]]}{[rs_1]}, \ldots, \frac{[p[q][s_l]]}{[rs_l]} \vdash \mathbf{p} \underset{[p[q]]}{\tilde{\circ}} \mathbf{q} \longrightarrow \mathbf{t} \underset{[r]}{\square} \mathbf{q}}\ \texttt{graft}$$

where $\Gamma'$ is given by pairs of the form

(1) $\frac{[a]}{[rxr']}$, where $\frac{[a]}{[r[y]r']} \in \Gamma$ and $\frac{[x]}{[y]} \in \Upsilon$,

(2) $\frac{[a]}{[b]}$, where $\frac{[a]}{[b]} \in \Gamma$ is not as above (i.e. $[b]$ not of the form $[r[y]r']$ for some $\frac{[x]}{[y]} \in \Upsilon$).

In large derivation trees, we will sometimes refer to this rule as `graft`-$[p[q]]$ for clarity, or simply as $[p[q]]$ in order to make notations lighter.

**Remark 4.1.7.** The transformation of context defined in rule `graft` in figure 4.1.1 might look complicated and unintuitive. And indeed it is, but we try to explain its purpose. Let $\Gamma \vdash \mathbf{p} \longrightarrow \mathbf{t}$ be a derivable sequent in Opt$^?$, with $\mathbf{p} \in \mathbb{P}_n$. As proved in lemma 4.1.9 on the facing page, $\Gamma$ exhibits a bijection between $\mathbf{p}^|$ and $\mathbf{t}^\bullet$. In theorem 4.2.4 on page 54, we prove that $\mathbf{p}$ corresponds uniquely to an $n$-opetope $\omega = [\![\mathbf{p}]\!]^{\text{poly}}$, that $\mathsf{t}\omega = [\![\mathbf{t}]\!]^{\text{poly}}$, and that $\Gamma$ corresponds to the readdressing function $\wp_\omega : \omega^| \longrightarrow (\mathsf{t}\omega)^\bullet$ (see section 2.1.1 and appendix A.1 on page 10 and on page 71).

We now turn our attention to the named approach of chapter 3 on page 15. Applying theorem 3.2.22 on page 25, we know that $\omega$ corresponds to a unique sequent (modulo variable renaming)

$$E \triangleright \Upsilon \vdash_n x : \mathsf{s}\, x \multimapdotinv \mathsf{s}^2\, x \multimapdotinv \cdots \multimapdotinv \varnothing$$

where $x \in \mathbb{V}$. More precisely, considered as a tree, $\omega$ is encoded by the term $\mathsf{s}\, x$, and by proposition 3.2.11 on page 22, $\mathsf{t}\,\omega$ is encoded by $\mathsf{s}^2\, x$. In lemma 3.2.10 on page 22, we show that $\wp_\omega$ exhibits a bijection

$$\left\{ \&_{\mathsf{s}\, x} b \mid b \in \mathbb{V}_{n-2}, b \in (\mathsf{s}^2\, x)^\bullet \right\} \xrightarrow{\wp_\omega} \left\{ \&_{\mathsf{s}^2\, x} b \mid b \in \mathbb{V}_{n-2}, b \in (\mathsf{s}^2\, x)^\bullet \right\}.$$

Say that a *node* of the term $\mathsf{s}^2\, x$ is a variable $b \in (\mathsf{s}^2\, x)^\bullet$, while a *leaf* of $\mathsf{s}\, x$ is a variable that can be used for grafting (see rule `graft` in figure 3.1.1 on page 17), i.e. a variable $b \in (\mathsf{s}^2\, x)^\bullet$. Then left hand side can be considered as the set of leaf addresses of $\mathsf{s}\, x$, while the right hand side is its set of node addresses of $\mathsf{s}^2\, x$, and $\wp_\omega$ maps the address of $b \in \mathbb{V}_{n-2}$ as a leaf of $\mathsf{s}\, x$ to the address of $b$ as a node of $\mathsf{s}^2\, x$. But here, the function $\wp_\omega$ is unnecessary: this correspondence is already established by the name of the variables!

In $\text{Opt}^?$ however, such bookkeeping is paramount since there is no such things as names. Based on the previous discussion, we conclude that $\Gamma$ is precisely the desired correspondence.

We now prove basic properties of derivable sequents in $\text{Opt}^?$. In proof trees, we may sometimes omit irrelevant information: for instance, if contexts and targets are not important, the `shift` rule may be written as

$$\frac{\mathbf{p}}{\{[] \leftarrow \mathbf{p}} \ \texttt{shift}$$

**Lemma 4.1.8.** *If* $\Gamma \vdash \mathbf{p} \longrightarrow \mathbf{t}$ *is a derivable sequent, then* $\dim \mathbf{p} = \dim \mathbf{t} + 1$.

*Proof.* We proceed by induction. The only non trivial case is `graft`. Remark that $\dim\left(\mathbf{p} \,\tilde{\circ}_{[p[q]]}\, \mathbf{q}\right) = \dim \mathbf{p} = \dim \mathbf{q} + 1 = \dim\left(\mathbf{t} \,\square_{[r]}\, \mathbf{q}\right) + 1$. ☐

**Lemma 4.1.9.** *Let* $\Gamma \vdash \mathbf{p} \longrightarrow \mathbf{t}$ *be a derivable sequent with* $\dim \mathbf{p} \geq 2$. *Then* $\Gamma$ *establishes a bijection between* $\mathbf{p}^|$ *and* $\mathbf{t}^\bullet$ *(i.e. as a set of pairs,* $\Gamma$ *is the graph of a bijective function).*

*Proof.* The fact that $\Gamma$ is a relation from $\mathbf{p}^|$ to $\mathbf{t}^\bullet$ (i.e. that whenever $\frac{[a]}{[b]} \in \Gamma$ we have $[a] \in \mathbf{p}^|$ and $[b] \in \mathbf{t}^\bullet$) is clear from the inference rules. It is also clear that $\Gamma$ is a function (i.e. that whenever $\frac{[a]}{[b]}, \frac{[a']}{[b']} \in \Gamma$ we have $[a] \neq [a']$). Finally, the fact that it is a bijection is clear in the case of `degen` and `shift`, and is a routine verification in the case of `graft`. ☐

**Proposition 4.1.10.** *If* $\Gamma \vdash \mathbf{p} \longrightarrow \mathbf{t}$ *is derivable, then* $\mathbf{t}$ *is too.*

*Proof.* The only non obvious case is (as always) `graft`, where we have to show that $\mathbf{t} \,\square_{[r]}\, \mathbf{q}$ is derivable. Since the sequent $\Gamma, \frac{[p[q]]}{[r]} \vdash \mathbf{p} \longrightarrow \mathbf{t}$ has a non empty context, $\mathbf{p}$ and $\mathbf{t}$ are non-degenerate. Write $\mathbf{t}$ and $\mathbf{q}$ as in definition 4.1.6 on page 49. Up to reindexing, assume that $\chi[t_j] = [t_j]$ if and only if $j < i$. Assume moreover that the sequence $[t_1], \ldots, [t_{i-1}]$ is lexicographically ordered, and likewise for $[t_{i+1}], \ldots, [t_l]$. For $j > i$ write $[t_j] = [t_i[b_j]x_j]$ and $\wp[t_j] = [t_i a_j x_j]$, so that $\Upsilon = \left\{ \frac{[a_j]}{[b_j]} \right\}_{i<j\leq l}$. Then the proof tree of $\mathbf{t} \,\square_{[t_i]}\, \mathbf{q}$ is sketched as follows:

(1) If $[t_i] = []$, then necessarily $i = 1$, and $\mathbf{t} \,\square_{[t_i]}\, \mathbf{q}$ is derived as follows:

$$\dfrac{\dfrac{\vdots \qquad \qquad \vdots}{\dfrac{\dfrac{\dfrac{\mathbf{q} \qquad \mathbf{v}_2}{\mathbf{q} \,\tilde{\circ}_{[a_2x_2]}\, \mathbf{v}_2}\ \texttt{graft-}[a_2x_2] \qquad \begin{matrix}\vdots \\ \mathbf{v}_3\end{matrix}}{\left(\mathbf{q} \,\tilde{\circ}_{[a_2x_2]}\, \mathbf{v}_2\right) \tilde{\circ}_{[a_3x_3]}\, \mathbf{v}_3}\ \texttt{graft-}[a_3x_3]}{\vdots}}{}$$

$$\dfrac{\left(\cdots\left(\mathbf{q} \,\tilde{\circ}_{[a_2x_2]}\, \mathbf{v}_2\right)\cdots\right) \tilde{\circ}_{[a_{l-1}x_{l-1}]}\, \mathbf{v}_{l-1} \qquad \begin{matrix}\vdots \\ \mathbf{v}_l\end{matrix}}{\left(\cdots\left(\mathbf{q} \,\tilde{\circ}_{[a_2x_2]}\, \mathbf{v}_2\right)\cdots\right) \tilde{\circ}_{[a_lx_l]}\, \mathbf{v}_l}\ \texttt{graft-}[a_lx_l]$$

and by definition, $\left(\cdots\left(\mathbf{q} \,\tilde{\circ}_{[a_2x_2]}\, \mathbf{v}_2\right)\cdots\right) \tilde{\circ}_{[a_lx_l]}\, \mathbf{v}_l = \mathbf{t} \,\square_{[t_i]}\, \mathbf{q}$.

(2) If $[t_i] \neq []$, then necessarily $i > 1$, and moreover $[t_1] = []$. The process goes similarly: we first derive $\begin{cases}[t_1] \leftarrow \mathbf{v}_1 \\ \vdots \\ [t_{i-1}] \leftarrow \mathbf{v}_{i-1}\end{cases}$, then graft the sources of $\mathbf{q}$, and lastly graft the remaining $\mathbf{v}_j$'s, where $j > i$. □

**Proposition 4.1.11.** *Let* $\Gamma_1 \vdash \mathbf{p} \longrightarrow \mathbf{t}_1$ *and* $\Gamma_2 \vdash \mathbf{p} \longrightarrow \mathbf{t}_2$ *be two derivable sequents. Then* $\Gamma_1 = \Gamma_2$ *and* $\mathbf{t}_1 = \mathbf{t}_2$*.*

*Proof.* If $\mathbf{p}$ is $\blacklozenge$, of the form $\{\!\{\mathbf{q}$, or of the form $\{[] \leftarrow \mathbf{q}$, then the result is clear, as both sequents have been obtained by `point`, `degen`, or `shift` respectively. Otherwise, both sequents have been obtained by an instance of the `graft` rule, so $\mathbf{p}$ has at least two addresses: $\#\,\mathbf{p}^\bullet \geq 2$. Assume $\#\,\mathbf{p}^\bullet = 2$, and write $\mathbf{p} = \begin{cases}[] \leftarrow \mathbf{q}_1 \\ [[q]] \leftarrow \mathbf{q}_2\end{cases}$, for $[q] \in \mathbf{q}_1^\bullet$. Note that, by induction, the sequents around $\mathbf{q}_1$ and $\mathbf{q}_2$ are uniquely determined. Then $\mathbf{p}$ has necessarily been derived as

$$\dfrac{\dfrac{\overset{\vdots}{\mathbf{q}_1}}{\{[] \leftarrow \mathbf{q}_1}\,\texttt{shift} \qquad \overset{\vdots}{\mathbf{q}_2}}{\mathbf{p}}\,\texttt{graft-}[[q]]$$

whence $\Gamma_1 = \Gamma_2$ and $\mathbf{t}_1 = \mathbf{t}_2$ in this case. Assume now $\#\,\mathbf{p}^\bullet > 2$, and that the two contexts in the proposition statement have been derived as

$$\dfrac{\overset{\vdots}{\mathbf{p}_1} \quad \overset{\vdots}{\mathbf{q}_1}}{\Gamma_1 \vdash \mathbf{p} \longrightarrow \mathbf{t}_1}\,\texttt{graft-}[l_1] \qquad \dfrac{\overset{\vdots}{\mathbf{p}_2} \quad \overset{\vdots}{\mathbf{q}_2}}{\Gamma_2 \vdash \mathbf{p} \longrightarrow \mathbf{t}_2}\,\texttt{graft-}[l_2]$$

Note that, by induction again, the sequents around $\mathbf{p}_i$ and $\mathbf{q}_i$, $i = 1, 2$, are uniquely determined. Then $\mathbf{p}_1$ is $\mathbf{p}$ with address $[l_1]$ removed, and $[l_2] \in \mathbf{p}_1^\bullet$ is an address "close to the leaves", in that $\forall [p] \neq [l_1] \in \mathbf{p}_1^\bullet$ we have $[l_1] \not\sqsupseteq [p]$. Likewise for $\mathbf{p}_2$. If $[l_1] = [l_2]$, then obviously the result holds, as both sequents have the same derivation tree. Otherwise, let $\mathbf{a}$ be $\mathbf{p}$ with addresses $[l_1]$ and $[l_2]$ removed. Then $\mathbf{a}$ is a "common ancestor" to $\mathbf{p}_1$ and $\mathbf{p}_2$:

$$\dfrac{\overset{\vdots}{\mathbf{a}} \quad \overset{\vdots}{\mathbf{q}_2}}{\mathbf{p}_1}\,\texttt{graft-}[l_2] \qquad \dfrac{\overset{\vdots}{\mathbf{a}} \quad \overset{\vdots}{\mathbf{q}_1}}{\mathbf{p}_2}\,\texttt{graft-}[l_1]$$

Again, the sequent around $\mathbf{a}$ is uniquely determined by induction. From there, $\Gamma_1 = \Gamma_2$, and $\mathbf{t}_1 = \mathbf{t}_2$ follow by routine verifications. □

We denote by $\mathbb{P}_n^{\checkmark}$ the set of derivable $n$-preopetopes, i.e. the $n$-preopetopes $\mathbf{p}$ such that there exists a derivable sequent $\Gamma \vdash \mathbf{p} \longrightarrow \mathbf{t}$. By proposition 4.1.11, this sequent is uniquely determined by $\mathbf{p}$, so let $\mathsf{t}\,\mathbf{p} = \mathbf{t}$ be the *target* of $\mathbf{p}$, and $\wp_{\mathbf{p}} : \mathbf{p}^| \xrightarrow{\cong} \mathbf{t}^\bullet$ be the bijection described by $\Gamma$. As such, the sequent around a derivable opetope $\mathbf{p}$ can be reconstructed as $\wp_{\mathbf{p}} \vdash \mathbf{p} \longrightarrow \mathsf{t}\,\mathbf{p}$.

**Remark 4.1.12.** Our syntax is closely related to the one given for multitopes [**Hermida et al., 2002**, Section 3], called here Hmp. Briefly, in Hmp, the unique 0 and 1-opetopes are respectively denoted ★ and # and, given an $n$-opetope $\mathbf{p}$, the notation $[\mathbf{p}]$ (resp. $\ulcorner\mathbf{p}\urcorner$) is used for the corresponding degenerated (resp. shifted) $(n+2)$ (resp. $(n+1)$) opetope. The nodes of an opetope come equipped with a canonical order (just as in our system we could require preopetopes to be always sorted according to the lexicographical order $\leq$), which apparently dispenses from using addresses. In Hmp, an inductive definition of opetopes is given, in the same spirit as our sequent calculus: in particular, typing conditions involving targets when grafting opetopes (grafting is simply application in Hmp) are involved. However, no explicit definition at the level of the syntax is given for computing targets (the description given resorts to multicategorical composition), and it is not clear to us how to define it without considering addresses and maintaining more information, as we do with our sequent calculus.

## 4.2. Equivalence with polynomial opetopes

We now show theorem 4.2.4 on page 54, stating that the elements of $\mathbb{P}_n^{\checkmark}$ are in bijective correspondence with the set $\mathbb{O}_n$ of polynomial $n$-opetopes. To this end, we now construct a bijection $[\![-]\!]^? : \mathbb{O}_n \longrightarrow \mathbb{P}_n^{\checkmark}$. If $n = 0, 1$, then both sets are singletons, so that $[\![-]\!]^?$ is trivially defined. Assume $n \geq 2$, and that $[\![-]\!]^?$ is defined for $k < n$. Recall that $\mathbb{O}_n = \operatorname{tr} \mathfrak{Z}^{n-2}$. We distinguish three cases:

(1) (Degenerate case) for $\phi \in \mathbb{O}_{n-2}$, let $[\![\mathsf{I}_\phi]\!]^? := \{\!\{[\![\phi]\!]^?$ ;
(2) (Corolla case) for $\psi \in \mathbb{O}_{n-1}$, let $[\![\mathsf{Y}_\psi]\!]^? := \{[] \leftarrow [\![\psi]\!]^?$ ;
(3) (Graft case) for $\omega \in \mathbb{O}_n$ that is neither degenerate nor a corolla, write $\omega = \nu \circ_{[l]} \mathsf{Y}_\psi$, for $\nu \in \mathbb{O}_n$, $\psi \in \mathbb{O}_{n-1}$, and $[l] \in \nu^|$, then let $[\![\omega]\!]^?$ be defined by the following proof tree:

$$\frac{[\![\nu]\!]^? \qquad [\![\psi]\!]^?}{[\![\omega]\!]^?} \, \texttt{graft-}[l]$$

or in other words, $[\![\omega]\!]^? = [\![\nu]\!]^? \,\tilde{\circ}_{[l]}\, [\![\psi]\!]^?$.

We prove that the last case makes sense in the following three lemmas.

**Lemma 4.2.1.** *For $\omega \in \mathbb{O}_n$, we have $\omega^\bullet = [\![\omega]\!]^{?\bullet}$, and if $\omega$ is non degenerate, then $\omega^| = [\![\omega]\!]^{?|}$.*

*Proof.* We proceed by opetopic induction (see remark 2.2.3 on page 12).

(1) (Degenerate case) Assume $\omega = \mathsf{I}_\phi$, for some $\phi \in \mathbb{O}_{n-2}$. Then $\omega^\bullet = \varnothing = \left(\{\!\{[\![\phi]\!]^?\right)^\bullet$. For leaves: $\omega^| = \varnothing = \left(\{\!\{[\![\phi]\!]^?\right)^| = [\![\omega]\!]^{?|}$.
(2) (Corolla case) Assume $\omega = \mathsf{Y}_\psi$, for some $\psi \in \mathbb{O}_{n-1}$, so that we have $[\![\omega]\!]^? = \{[] \leftarrow [\![\psi]\!]^?$. Then $\omega^\bullet = \{[]\} = \left(\{[] \leftarrow [\![\psi]\!]^?\right)^\bullet = [\![\omega]\!]^{?\bullet}$. By induction, $\psi^\bullet = [\![\psi]\!]^{?\bullet}$, so the leaf addresses of $\omega$ and $\{[] \leftarrow [\![\psi]\!]^?$ are both of the form $[[q]]$, where $[q] \in \psi^\bullet$, hence $\omega^| = [\![\omega]\!]^{?|}$.
(3) (Graft case) Assume $\omega = \nu \circ_{[l]} \mathsf{Y}_\psi$, for some $\nu \in \mathbb{O}_n$, $\psi \in \mathbb{O}_{n-1}$, and $[l] \in \nu^|$. Then, by induction,

$$\omega^\bullet = \nu^\bullet + \{[l]\} = [\![\nu]\!]^{?\bullet} + \{[l]\} = \left([\![\nu]\!]^? \underset{[l]}{\tilde{\circ}} [\![\psi]\!]^?\right)^\bullet = [\![\omega]\!]^{?\bullet}.$$

Similarly, by induction,

$$\begin{aligned} \omega^| &= \nu^| - \{[l]\} + \{[l[q]] \mid [q] \in \psi^\bullet\} \\ &= [\![\nu]\!]^{?|} - \{[l]\} + \left\{[l[q]] \mid [q] \in [\![\psi]\!]^{?\bullet}\right\} \\ &= \left([\![\nu]\!]^? \underset{[l]}{\tilde{\circ}} [\![\psi]\!]^?\right)^| = [\![\omega]\!]^{?|}. \end{aligned}$$

and we conclude. □

**Lemma 4.2.2.** *For $\omega \in \mathbb{O}_n$ and $[p] \in \omega^\bullet$, we have $[\![\mathsf{s}_{[p]}\,\omega]\!]^? = \mathsf{s}_{[p]}\,[\![\omega]\!]^?$.*

*Proof.* We proceed by opetopic induction.

(1) (Degenerate case) If $\omega$ is degenerate, then $\omega^\bullet = \varnothing$, so there is nothing to check in this case.
(2) (Corolla case) Assume $\omega = \mathsf{Y}_\psi$, for some $\psi \in \mathbb{O}_{n-1}$. Then $\omega^\bullet = \{[]\}$, and $[\![\mathsf{s}_{[]}\,\omega]\!]^? = [\![\psi]\!]^? = \mathsf{s}_{[]}\,\{[] \leftarrow [\![\psi]\!]^? = \mathsf{s}_{[]}\,[\![\omega]\!]^?$.
(3) (Graft case) Assume $\omega = \nu \circ_{[l]} \mathsf{Y}_\psi$, for some $\nu \in \mathbb{O}_n$, $\psi \in \mathbb{O}_{n-1}$, and $[l] \in \nu^|$. Let $[p] \in \omega^\bullet$. If $[p] = [l]$, then

$$[\![\mathsf{s}_{[l]}\,\omega]\!]^? = [\![\psi]\!]^? = \mathsf{s}_{[l]}\left([\![\nu]\!]^? \underset{[l]}{\tilde{\circ}} [\![\psi]\!]^?\right) = \mathsf{s}_{[l]}\,[\![\omega]\!]^?.$$

Otherwise, we have

$$[\![\mathsf{s}_{[p]}\,\omega]\!]^? = [\![\mathsf{s}_{[p]}\,\nu]\!]^? = \mathsf{s}_{[p]}\,[\![\nu]\!]^? = \mathsf{s}_{[p]}\left([\![\nu]\!]^? \underset{[l]}{\tilde{\circ}} [\![\psi]\!]^?\right) = \mathsf{s}_{[p]}\,[\![\omega]\!]^?.$$

□

**Lemma 4.2.3.** *For $\omega \in \mathbb{O}_n$, we have $[\![\mathsf{t}\,\omega]\!]^? = \mathsf{t}\,[\![\omega]\!]^?$, and $\wp_\omega = \wp_{[\![\omega]\!]^?}$.*

*Proof.* We proceed by opetopic induction.

(1) (Degenerate case) Assume $\omega = \mathsf{I}_\phi$, for some $\phi \in \mathbb{O}_{n-2}$. Then $[\![\mathsf{t}\,\omega]\!]^? = [\![\mathsf{Y}_\phi]\!]^? = \big\{[] \leftarrow [\![\phi]\!]^? = \mathsf{t}\big\{\!\!\big\{[\![\phi]\!]^? = \mathsf{t}\,[\![\omega]\!]^?$. Since $\omega$ and $[\![\omega]\!]^?$ are both degenerate (as opetope and preopetope, respectively), $\wp_\omega$ and $\wp_{[\![\omega]\!]^?}$ both map $[]$ (the unique leaf of $\omega$ and $[\![\omega]\!]^?$) to $[]$ (the unique node of $\mathsf{t}\,\omega$ and $\mathsf{t}\,[\![\omega]\!]^?$).

(2) (Corolla case) Assume $\omega = \mathsf{Y}_\psi$, for some $\psi \in \mathbb{O}_{n-1}$, so that we have $[\![\omega]\!]^? = \big\{[] \leftarrow [\![\psi]\!]^?$. Then $[\![\mathsf{t}\,\omega]\!]^? = [\![\psi]\!]^? = \mathsf{t}\big\{[] \leftarrow [\![\psi]\!]^? = \mathsf{t}\,[\![\omega]\!]^?$. Moreover, we have $\omega^| = [\![\omega]\!]^{?|} = \{[[p]] \mid [p] \in \psi^\bullet\}$, and by definition, $\wp_\omega[[p]] = [p] = \wp_{[\![\omega]\!]^?}[p]$.

(3) (Graft case) Assume $\omega = \nu \circ_{[l]} \mathsf{Y}_\psi$, for some $\nu \in \mathbb{O}_n$, $\psi \in \mathbb{O}_{n-1}$, and $[l] \in \nu^|$. Then,

$$
\begin{aligned}
[\![\mathsf{t}\,\omega]\!]^? &= \left[\!\!\left[(\mathsf{t}\,\nu) \underset{\wp_\nu[l]}{\square} \psi\right]\!\!\right]^? \\
&= [\![\mathsf{t}\,\nu]\!]^? \underset{\wp_{[\![\nu]\!]^?}[l]}{\square} [\![\psi]\!]^? && \text{by ind. } \wp_\nu = \wp_{[\![\nu]\!]^?} \\
&= (\mathsf{t}\,[\![\nu]\!]^?) \underset{\wp_{[\![\nu]\!]^?}[l]}{\square} [\![\psi]\!]^? && \text{by induction} \\
&= \mathsf{t}\left([\![\nu]\!]^? \underset{[l]}{\tilde{\circ}} [\![\psi]\!]^?\right) \\
&= \mathsf{t}\,[\![\omega]\!]^?
\end{aligned}
$$

Lastly, we prove $\wp_\omega = \wp_{[\![\omega]\!]^?}$. Let $[u] \in \omega^|$.

(a) If $[u] = [l[q]]$ for $[q] \in \psi^\bullet$, then $\wp_\omega[l[q]] = \wp_\nu[l] \cdot [q]$ (see line (A.2.3)), where $U = \mathsf{t}\,\nu$, $T = \psi$, and $[u] = \wp_\nu[l])$. On the other hand, $\wp_{[\![\omega]\!]^?}[l[q]] = \wp_{[\![\nu]\!]^?}[l] \cdot [q]$.

(b) If $[u] \in \nu^|$ is of the form $\wp_\nu[u] = [p[\wp_\psi[v]]p']$, for some $[v] \in \psi^|$, then $\wp_\omega[u] = [pvp']$. On the other hand, by induction, $\wp_{[\![\nu]\!]^?}[u] = \wp_\nu[u] = [p[\wp_\psi[v]]p'] = [p[\wp_{[\![\psi]\!]^?}[v]]p']$, and by definition we have $\wp_{[\![\omega]\!]^?}[u] = [pvp']$.

(c) If $[u] \in \nu^|$ is not of the form above, then $\wp_\omega[u] = \wp_\nu[u] = \wp_{[\![\nu]\!]^?}[u] = \wp_{[\![\omega]\!]^?}[u]$. ☐

We finally prove that $[\![-]\!]^?$ is a bijection by constructing its inverse $[\![-]\!]^{\text{poly}} : \mathbb{P}_n^{\checkmark} \longrightarrow \mathbb{O}_n$. If $n = 0, 1$, then both sets are singletons, so that $[\![-]\!]^{\text{poly}}$ is trivially defined. Assume $n \geq 2$, and that $[\![-]\!]^{\text{poly}}$ is defined for $k < n$. We distinguish three cases:

(1) (Degenerate case) if $\mathbf{q} \in \mathbb{P}_{n-2}^{\checkmark}$, then $\big[\!\!\big[\big\{\!\!\big\{\mathbf{q}\big]\!\!\big]^{\text{poly}} := \mathsf{I}_{[\![\mathbf{q}]\!]^{\text{poly}}}$;

(2) (Corolla case) if $\mathbf{q} \in \mathbb{P}_{n-1}^{\checkmark}$, then $\big[\!\!\big[\big\{[] \leftarrow \mathbf{q}\big]\!\!\big]^{\text{poly}} := \mathsf{Y}_{[\![\mathbf{q}]\!]^{\text{poly}}}$;

(3) (Graft case) for $\mathbf{p} \in \mathbb{P}_n^{\checkmark}$, $\mathbf{q} \in \mathbb{P}_{n-1}^{\checkmark}$, and $[l] \in \mathbf{p}^|$, define

$$
\left[\!\!\left[\mathbf{p} \underset{[l]}{\tilde{\circ}} \mathbf{q}\right]\!\!\right]^{\text{poly}} := [\![\mathbf{p}]\!]^{\text{poly}} \underset{[l]}{\circ} \mathsf{Y}_{[\![\mathbf{q}]\!]^{\text{poly}}}.
$$

**Theorem 4.2.4.** *The functions $[\![-]\!]^?$ and $[\![-]\!]^{\text{poly}}$ are mutually inverse. Moreover, for $n \geq 2$ and $\omega \in \mathbb{O}_n$, we have*

*(1) $\omega^\bullet = [\![\omega]\!]^{?\bullet}$, and if $\omega$ is non degenerate, $\omega^| = [\![\omega]\!]^{?|}$;*

*(2) for $[p] \in \omega^\bullet$, $[\![\mathsf{s}_{[p]}\,\omega]\!]^? = \mathsf{s}_{[p]}\,[\![\omega]\!]^?$, and $[\![\mathsf{t}\,\omega]\!]^? = \mathsf{t}\,[\![\omega]\!]^?$;*

*(3) $\wp_\omega = \wp_{[\![\omega]\!]^?}$.*

## 4.3. Examples

In this section, we give example derivations in system $\textsc{Opt}^?$.

**Example 4.3.1** (The arrow)**.** The unique 1-opetope $\blacksquare = \big\{* \leftarrow \blacklozenge$ is derived by

$$
\dfrac{\dfrac{}{\vdash \blacklozenge \longrightarrow \varnothing}\ \texttt{point}}{\vdash \big\{* \leftarrow \blacklozenge \longrightarrow \blacklozenge}\ \texttt{shift}
$$

**Example 4.3.2** (Opetopic integers)**.** As in example 3.3.2 on page 26, we provide the derivation of opetopic integers $\mathbf{n} \in \mathbb{O}_2$. The derivation of $\mathbf{0}$ is given on the left, whereas for $n \geq 1$, the opetope $\mathbf{n}$ is derived as on the right:

$$\dfrac{\dfrac{}{\vdash \blacklozenge}\ \texttt{point}}{\frac{[]}{[]} \vdash \left\{\!\!\left\{\blacklozenge \longrightarrow \blacksquare\right.\right.}\ \texttt{degen}$$

$$\dfrac{\dfrac{\dfrac{\dfrac{\dfrac{\vdots}{\vdash \blacksquare \longrightarrow \blacklozenge}}{\frac{[*]}{*} \vdash \left\{[] \leftarrow \blacksquare \longrightarrow \blacksquare\right.}\ \texttt{shift} \qquad \dfrac{\vdots}{\vdash \blacksquare \longrightarrow \blacklozenge}}{\frac{[**]}{*} \vdash \begin{cases}[] \leftarrow \blacksquare \\ [*] \leftarrow \blacksquare\end{cases} \longrightarrow \blacksquare}\ [*] \qquad \dfrac{\vdots}{\vdash \blacksquare \longrightarrow \blacklozenge}}{\frac{[***]}{*} \vdash \begin{cases}[] \leftarrow \blacksquare \\ [*] \leftarrow \blacksquare \\ [**] \leftarrow \blacksquare\end{cases} \longrightarrow \blacksquare}\ [**]$$

$$\vdots$$

$$\dfrac{\frac{[*^{n-1}]}{*} \vdash \begin{cases}[] \leftarrow \blacksquare \\ \vdots \\ [*^{n-2}] \leftarrow \blacksquare\end{cases} \longrightarrow \blacksquare \qquad \dfrac{\vdots}{\vdash \blacksquare \longrightarrow \blacklozenge}}{\frac{[*^{n}]}{*} \vdash \begin{cases}[] \leftarrow \blacksquare \\ \vdots \\ [*^{n-1}] \leftarrow \blacksquare\end{cases} \longrightarrow \blacksquare}\ [*^{n-1}]$$

where $*^n = \underbrace{*\cdots*}_{n}$, and where there is a total of $n-1$ instances of the `graft` rule.

**Example 4.3.3** (A classic)**.** The 3-opetope on the left is derived as on the right:

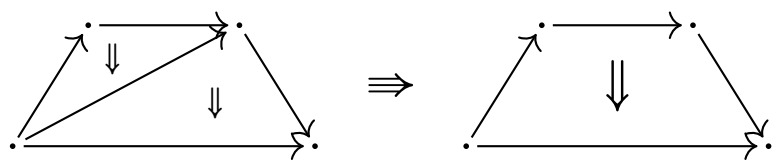

$$\dfrac{\dfrac{\dfrac{\vdots}{\frac{[**]}{*} \vdash \mathbf{2} \longrightarrow \blacksquare}}{\frac{[[*]]}{[*]} \vdash \left\{[] \leftarrow \mathbf{2} \longrightarrow \mathbf{2}\right.}\ \texttt{shift} \qquad \dfrac{\vdots}{\frac{[**]}{*} \vdash \mathbf{2} \longrightarrow \blacksquare}}{\frac{[[]]}{[]}, \frac{[[*][]]}{[*]}, \frac{[[*][*]]}{[**]} \vdash \begin{cases}[] \leftarrow \mathbf{2} \\ [[*]] \leftarrow \mathbf{2}\end{cases} \longrightarrow \mathbf{3}}\ \texttt{graft-}[[*]]$$

**Example 4.3.4** (A degenerate case)**.** The 3-opetope on the left is derived as on the right:

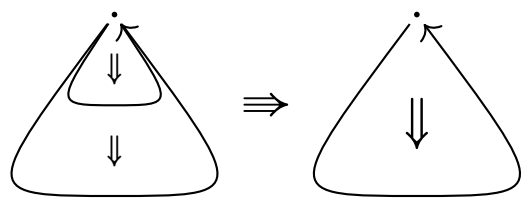

$$\dfrac{\dfrac{\dfrac{\vdots}{\frac{[*]}{*} \vdash \mathbf{1} \longrightarrow \blacksquare}}{\frac{[[]]}{[]} \vdash \left\{[] \leftarrow \mathbf{1} \longrightarrow \mathbf{1}\right.}\ \texttt{shift} \qquad \dfrac{\vdots}{\frac{[]}{[]} \vdash \mathbf{0} \longrightarrow \blacksquare}}{\vdash \begin{cases}[] \leftarrow \mathbf{1} \\ [[]] \leftarrow \mathbf{0}\end{cases} \longrightarrow \mathbf{0}}\ \texttt{graft-}[[]]$$

**Example 4.3.5** (Another degenerate case)**.** The 3-opetope on the left is derived as on the right:

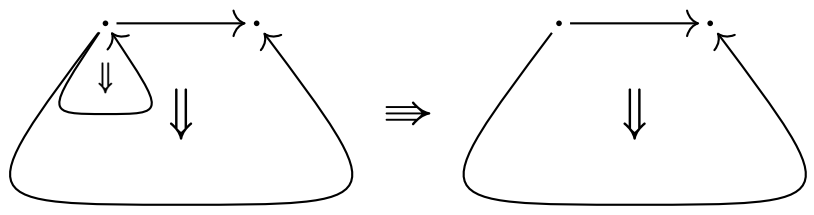

$$\dfrac{\dfrac{\dfrac{\vdots}{\frac{[**]}{*} \vdash \mathbf{2} \longrightarrow \blacksquare}}{\frac{[[]]}{[]}, \frac{[[*]]}{[*]} \vdash \left\{[] : \mathbf{2} \longrightarrow \mathbf{2}\right.}\ \texttt{shift} \qquad \dfrac{\vdots}{\vdash \mathbf{0} \longrightarrow \blacksquare}}{\frac{[[]]}{[]} \vdash \begin{cases}[] \leftarrow \mathbf{2} \\ [[*]] \leftarrow \mathbf{0}\end{cases} \longrightarrow \mathbf{1}}\ \texttt{graft-}[[*]]$$

**Example 4.3.6** (A 4-opetope)**.** The 4-opetope

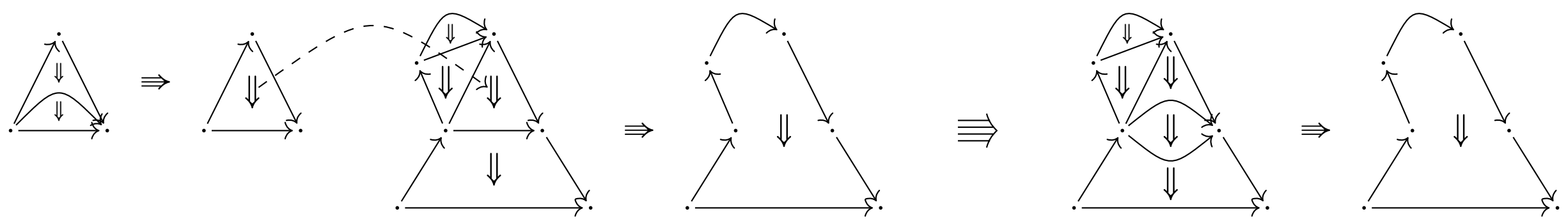

is derived by

$$
\dfrac{
\dfrac{
\dfrac{
\dfrac{
\dfrac{
\dfrac{
\begin{array}{c}\vdots\\ \frac{[***]}{*}\vdash\mathbf{3}\longrightarrow\blacksquare\end{array}
}{
\frac{[[]]}{[]},\frac{[[*]]}{[*]},\frac{[[**]]}{[**]}\vdash\{[]\leftarrow\mathbf{3}\longrightarrow\mathbf{3}
}\,\mathsf{s}
\qquad
\begin{array}{c}\vdots\\ \frac{[**]}{*}\vdash\mathbf{2}\longrightarrow\blacksquare\end{array}
}{
\frac{[[]]}{[]},\frac{[[**]]}{[***]},\frac{[[*][]]}{[*]},\frac{[[*][*]]}{[**]}\vdash\begin{cases}[]\leftarrow\mathbf{3}\\ [[*]]\leftarrow\mathbf{2}\end{cases}\longrightarrow\mathbf{4}
}\,[[*]]
\qquad
\begin{array}{c}\vdots\\ \frac{[**]}{*}\vdash\mathbf{2}\longrightarrow\blacksquare\end{array}
}{
\frac{[[]]}{[]},\frac{[[**]]}{[****]},\frac{[[*][]]}{[*]},\frac{[[*][*][]]}{[**]},\frac{[[*][*][*]]}{[***]}\vdash\begin{cases}[]\leftarrow\mathbf{3}\\ [[*]]\leftarrow\mathbf{2}\\ [[*][*]]\leftarrow\mathbf{2}\end{cases}\longrightarrow\mathbf{5}
}\,[[*][*]]
\qquad
\begin{array}{c}\vdots\\ \frac{[*]}{*}\vdash\mathbf{1}\longrightarrow\blacksquare\end{array}
}{
\frac{[[]]}{[]},\frac{[[**]]}{[****]},\frac{[[*][]]}{[*]},\frac{[[*][*][*]]}{[***]},\frac{[[*][*][][]]}{[**]}\vdash\begin{cases}[]\leftarrow\mathbf{3}\\ [[*]]\leftarrow\mathbf{2}\\ [[*][*]]\leftarrow\mathbf{2}\\ [[*][*][]]\leftarrow\mathbf{1}\end{cases}\longrightarrow\mathbf{5}
}\,[[*][*][]]
}{
\frac{[[]]}{[]},\frac{[[[*]]]}{[[*]]},\frac{[[[*][*]]]}{[[*][*]]},\frac{[[[*][*][]]]}{[[*][*][]]}\vdash\left\{[]\leftarrow\begin{cases}[]\leftarrow\mathbf{3}\\ [[*]]\leftarrow\mathbf{2}\\ [[*][*]]\leftarrow\mathbf{2}\\ [[*][*][]]\leftarrow\mathbf{1}\end{cases}\right.\longrightarrow\begin{cases}[]\leftarrow\mathbf{3}\\ [[*]]\leftarrow\mathbf{2}\\ [[*][*]]\leftarrow\mathbf{2}\\ [[*][*][]]\leftarrow\mathbf{1}\end{cases}
}\,\mathsf{s}
\qquad
\dfrac{
\dfrac{
\begin{array}{c}\vdots\\ \frac{[*]}{*}\vdash\mathbf{1}\longrightarrow\blacksquare\end{array}
}{
\frac{[[]]}{[]}\vdash\{[]\leftarrow\mathbf{1}\longrightarrow\mathbf{1}
}\,\mathsf{s}
\qquad
\begin{array}{c}\vdots\\ \frac{[**]}{*}\vdash\mathbf{2}\longrightarrow\blacksquare\end{array}
}{
\frac{[[][]]}{[]},\frac{[[][*]]}{[*]}\vdash\begin{cases}[]\leftarrow\mathbf{1}\\ [[]]\leftarrow\mathbf{2}\end{cases}\longrightarrow\mathbf{2}
}\,[[]]
}{
\frac{[[]]}{[]},\frac{[[[*][*]]]}{[[*][][*]]},\frac{[[[*][*][]]]}{[[*][][*][]]},\frac{[[[*]][[]]]}{[[*]]},\frac{[[[*]][[[]]]]}{[[*][]]}\vdash\begin{cases}[]\leftarrow\begin{cases}[]\leftarrow\mathbf{3}\\ [[*]]\leftarrow\mathbf{2}\\ [[*][*]]\leftarrow\mathbf{2}\\ [[*][*][]]\leftarrow\mathbf{1}\end{cases}\\ [[[*]]]\leftarrow\begin{cases}[]\leftarrow\mathbf{1}\\ [[]]\leftarrow\mathbf{2}\end{cases}\end{cases}\longrightarrow\begin{cases}[]\leftarrow\mathbf{3}\\ [[*]]\leftarrow\mathbf{1}\\ [[*][]]\leftarrow\mathbf{2}\\ [[*][][*]]\leftarrow\mathbf{2}\\ [[*][][*][]]\leftarrow\mathbf{1}\end{cases}
}\,[[[*]]]
$$

where rule s is a shorthand for `shift`.

## 4.4. Deciding opetopes

We now present in algorithm 4.4.1 the isOpetope function that, given a preopetope $\omega \in \mathbb{P}$, decides if $\omega \in \mathbb{P}^{\checkmark}$, as proved in proposition 4.4.1 on the next page. This algorithm tries to deconstruct $\omega$ by finding the last rule instance if its potential proof tree, and recursively checking the validity of the premises. We emphasize that this algorithm, while straightforward, is extremely inefficient.

**Algorithm 4.4.1** Well formation algorithm

```
procedure isOpetope(ω ∈ ℙ)                                   ▷ Returns a boolean
    if ω = ⧫ then
        return true
    else if ω = ⦃ϕ then
        return isOpetope(ϕ)
    else
        while ω has an address of the form [p[q]] do
            if [p] ∉ ω• or not isOpetope(s_[p] ω) then
                return false
            else if [q] ∉ (s_[p] ω)• or not isOpetope(s_[p[q]] ω) then
                return false
            else if t s_[p] ω ≠ s_[p[q]] ω then
                return false
            else
                Remove address [p[q]] from ω
            end if
        end while
        if ω is of the form {[] ← ψ then
            return isOpetope(ψ)
        else
            return false
        end if
    end if
end procedure
```

**Proposition 4.4.1.** *For $\omega \in \mathbb{P}$, the execution* IsOpetope($\omega$) *returns* ***true*** *if and only if* $\omega \in \mathbb{P}^{\checkmark}$.

*Proof.* This algorithm tries to deconstruct the potential proof tree of $\omega$ in system $\textsc{Opt}^{?}$:

(1) condition at line (2) removes an instance of the `point` rule;
(2) condition at line (4) removes an instance of `degen`;
(3) each iteration of the **while** loop at line (7) removes an instance of `graft`;
(4) finally the condition at line (18) removes an instance of `shift`.

If the algorithm encounters an expression that is not the conclusion of any instance of any rule of $\textsc{Opt}^{?}$, it returns **false**. Otherwise, if all branches of the proof tree lead to ♦, it returns **true**. □

## 4.5. Python implementation

The derivation system $\textsc{Opt}^{?}$ and all required syntactic constructs are implemented in module `opetopy.UnnamedOpetope` of [**Ho Thanh, 2018b**]. In particular, the four derivation rules are represented by functions of the same name: `point`, `degen`, `shift`, and `graft`. As an example, we review the implementation of `shift` in figure 4.5.1.

Figure 4.5.1. Implementation of $\textsc{Opt}^{?}$'s `shift` rule in `opetopy.UnnamedOpetope.shift`

```python
# This function takes a sequent (opetopy.UnnamedOpetope.Sequent) and returns a sequent. A
↪   sequent seq is structured as follows: for seq = Γ ⊢ s ⟶ t we have seq.context = Γ,
↪   seq.source = s, and seq.target = t.
def shift(seq: Sequent) -> Sequent:
    # We let n be the dimension of the preopetope s = seq.source
    n = seq.source.dimension
    # We construct a new context ctx dimension n + 1
    ctx = Context(n + 1)
    # We let ctx = { [a]/a | a ∈ s• }
    for a in seq.source.nodeAddresses():
        ctx += (a.shift(), a)
    # We return the sequent ctx ⊢ {[] ← s ⟶ s
    return Sequent(
        ctx,
        Preopetope.fromDictOfPreopetopes({
            Address.epsilon(n): seq.source
        }),
        seq.source
    )
```

To construct proof trees, those rules are further abstracted in classes `Point`, `Degen`, `Shift`, as well as `Graft`. We review the implementation of some examples in figures 4.5.2 to 4.5.4 on pages 58–59.

FIGURE 4.5.2. Derivation of the arrow sequent (see example 4.3.1 on page 54) using `opetopy.UnnamedOpetope`

```
from opetopy.UnnamedOpetope import *

# The arrow sequent is obtained by an application of the point rule followed by an application of the shift rule.
arrow = Shift(Point())
# Note that the function opetopy.UnnamedOpetope.Arrow can be used to concisely get the proof tree of ■.
```

FIGURE 4.5.3. Derivation of some opetopic integers (see example 4.3.2 on page 55) using `opetopy.UnnamedOpetope`

```
from opetopy.UnnamedOpetope import *

# The opetopic integer 0 is obtained by an instance of the point rule, followed by an application of the degen rule.
opetopic_integer_0 = Degen(Point())
# The opetopic integer 1 is obtained by applying rule shift to the arrow ■ as defined in the previous figure.
opetopic_integer_1 = Shift(arrow)
# The opetopic integer 2 is defined by 2 = 1 õ[*] ■. The address [*] is obtained with the convenient UnnamedOpetope.address function (as opposed to using the UnnamedOpetope.Address class).
opetopic_integer_2 = Graft(
    opetopic_integer_1,
    arrow,
    address(['*']))
# Likewise, 3 = 2 õ[**] ■.
opetopic_integer_3 = Graft(
    opetopic_integer_2,
    arrow,
    address(['*', '*']))
# Note that the function opetopy.UnnamedOpetope.OpetopicInteger can be used to get the proof tree of an arbitrary opetopic integer.
```

FIGURE 4.5.4. Derivation of example 4.3.3 on page 55 using `opetopy.UnnamedOpetope`

```
from opetopy.UnnamedOpetope import *

# Recall that in this example, the final opetope ω is defined by ω := ({[] ← 2) õ_[[*]] 2.
example_classic = Graft(
    Shift(opetopic_integer_2),
    opetopic_integer_2,
    address([['*']]))
```

## 4.6. The system for opetopic sets

We now present $\text{OptSet}^?$, a derivation system for opetopic set that is *controlled* by $\text{Opt}^?$ (unlike system $\text{OptSet}^!$ that is *based* on $\text{Opt}^!$, see figures 3.1.1 and 3.5.1 on page 17 and on page 30).

As always, contexts are considered as sets, so that the order in which the typings are written is irrelevant, even though those typings might be interdependent. We rely on two types of judgment that can be understood as follows:

(1) $\Gamma$ context means that $\Gamma$ is a well-formed context,
(2) $\Gamma \vdash \mathbf{P}$ means that in context $\Gamma$, $\mathbf{P}$ is a well-formed pasting diagram.

We now state the inference rules in figure 4.6.1 on the next page, and simultaneously assign a *shape* $x^\natural$ to any variable $x$ in a derivable context.

FIGURE 4.6.1. The OPTSET$^?$ system.

**Points:**

$$\frac{\Gamma \ \mathsf{context}}{\Gamma, x : \blacklozenge \ \mathsf{context}}\ \mathtt{point}$$

for $x$ a fresh name. Such a cell $x$ has no source, no target, and its shape is given by $x^\natural := \blacklozenge$.

**Degenerate pasting diagrams:**

$$\frac{\Gamma, x : T \ \mathsf{context}}{\Gamma, x : T \vdash \{\{ x}\ \mathtt{degen}$$

The shape of this pasting diagram is $\left(\{\{ x\right)^\natural := \{\{ x^\natural$.

**Non degenerate pasting diagrams:** If there exists $\mathbf{p} \in \mathbb{P}^\checkmark$ a non degenerate opetope

$$\mathbf{p} = \begin{cases} [p_1] \leftarrow \psi_1 \\ \vdots \\ [p_k] \leftarrow \psi_k \end{cases}$$

and variables $x_i : T_i$, for $1 \leq i \leq k$, such that

(1) the cell $x_i : T_i$ is such that $x_i^\natural = \psi_i$,
(2) **(Inner)** whenever $[p_j] = [p_i[q]]$ we have $\mathsf{t}\, x_j = \mathsf{s}_{[q]}\, x_i$ (the latter notation is defined in rule `shift`),

then:

$$\frac{\Gamma, x_1 : T_1, \ldots, x_k : T_k \ \mathsf{context}}{\Gamma, x_1 : T_1, \ldots, x_k : T_k \vdash \begin{cases} [p_1] \leftarrow x_1 \\ \vdots \\ [p_k] \leftarrow x_k \end{cases}}\ \mathtt{graft}$$

Denote this pasting diagram (i.e. the big expression on the right hand side of $\vdash$ in the conclusion of the above rule) by $\mathbf{P}$. Its shape is given by $\mathbf{P}^\natural := \mathbf{p}$, and let $\mathsf{s}_{[p_i]}\, \mathbf{P} := x_i$. Forming a pasting diagram in this manner is essentially an unbiased (or non binary) grafting, whence the name of the rule.

**Shift to the next dimension:** If we have a pasting diagram $\mathbf{P}$ of shape $\mathbf{P}^\natural = \mathbf{p}$, a cell $x : \mathbf{Q} \longrightarrow a$, such that

(1) $x^\natural = \mathsf{t}\, \mathbf{p}$,
(2) **(Glob1)** if $\mathbf{p}$ is non degenerate, we have $\mathsf{t}\,\mathsf{s}_{[]}\, \mathbf{P} = \mathsf{t}\, x$,
(3) **(Glob2)** if $\mathbf{p}$ is non degenerate, for all leaves $[p[q]]$ of $\mathbf{p}$ we have $\mathsf{s}_{[q]}\, \mathsf{s}_{[p]}\, \mathbf{P} = \mathsf{s}_{\wp_{\mathbf{p}}[p[q]]}\, x$,
(4) **(Degen)** if $\mathbf{p}$ is degenerate, we have $\mathbf{Q} = \{[] \leftarrow a$,

then:

$$\frac{\Gamma, x : \mathbf{Q} \longrightarrow a \vdash \mathbf{P}}{\Gamma, x : \mathbf{Q} \longrightarrow a, y : \mathbf{P} \longrightarrow x \ \mathsf{context}}\ \mathtt{shift}$$

for $y$ a fresh name. The shape of $y$ is given by $y^\natural := \mathbf{P}^\natural$, its source is $\mathsf{s}\, y := \mathbf{P}$, for $[p]$ an address in $\mathbf{P}$, its $[p]$-source is $\mathsf{s}_{[p]}\, y := \mathsf{s}_{[p]}\, \mathbf{P}$, and its target is $\mathsf{t}\, y := x$.

## 4.7. EQUIVALENCE WITH OPETOPIC SETS

Let $\Upsilon$ and $\Gamma = (x_1 : T_1, \ldots, x_k : T_k)$ be two derivable contexts in OPTSET$^?$. Akin to the standard definition, a substitution $\sigma : \Upsilon \longrightarrow \Gamma$ is a sequence of expressions $(\sigma_1, \ldots, \sigma_k)$ such that for $1 \leq i \leq k$ we have

$$\sigma_i : T_i[\sigma_1/x_1]\cdots[\sigma_{i-1}/x_{i-1}] \in \Upsilon.$$

**Lemma 4.7.1.** *In the setting above, we have $\sigma_i^\natural = x_i^\natural$.*

*Proof.* We proceed by induction. By definition, $\sigma_1 = x_1$, and thus $\sigma_1^\natural = x_1^\natural$. By induction, take $1 \leq i \leq k$ and assume $\sigma_j^\natural = x_j^\natural$ for $j < i$. Since the shape of a pasting diagram $\mathbf{P} = \begin{cases} [p_1] \leftarrow y_1 \\ \vdots \end{cases}$ only depends on the shape of the

$y_j$'s, we have

$$x_i^\natural = T_i^\natural = (T_i[\sigma_1/x_1]\cdots[\sigma_{i-1}/x_{i-1}])^\natural = \sigma_i^\natural.$$

□

Let $\mathcal{C}\mathrm{tx}^?$ be the syntactic category of our sequent calculus, i.e. the category whose objects are derivable contexts, and morphisms are substitutions as defined above. We now construct the *unnamed stratification functor* $S^? : (\mathcal{C}\mathrm{tx}^?)^{\mathrm{op}} \longrightarrow \mathcal{F}\mathrm{in}\widehat{\mathbb{O}}$. For $\Gamma \in \mathcal{C}\mathrm{tx}^?$ and $\omega \in \mathbb{O}$, let

$$S^?\Gamma_\omega = \left\{ x \in \Gamma \mid x^\natural = [\![\omega]\!]^? \right\}.$$

If $x^\natural \neq \blacklozenge$, then the type $X$ of $x$ is of the form $\mathbf{P} \longrightarrow z$, and we let $\mathsf{t}\,x := z$. This is well defined as by construction of $\Gamma$ we have $z^\natural = \mathsf{t}\,x^\natural$. For $[p] \in \omega^\bullet$, we let $\mathsf{s}_{[p]}\,x := \mathsf{s}_{[p]}\,\mathbf{P}$. Again, this is well defined as $(\mathsf{s}_{[p]}\,\mathbf{P})^\natural = \mathsf{s}_{[p]}\,\mathbf{P}^\natural = \mathsf{s}_{[p]}\,x^\natural$. From there, the opetopic identities clearly hold, and $S^?\Gamma$ is a finite opetopic set. To abbreviate notations, and for $\mathsf{f} : \psi \longrightarrow \omega$ a morphism in $\mathbb{O}$, we let $\mathsf{f} = S^?\Gamma\mathsf{f} : S^?\Gamma_\omega \longrightarrow S^?\Gamma_\psi$.

Take now $\sigma = (\sigma_1, \ldots, \sigma_k) : \Upsilon \longrightarrow \Gamma$ a substitution, where as before $\Gamma = (x_1 : T_1, \ldots, x_k : T_k)$. We define a morphism $S^?\sigma : S^?\Gamma \longrightarrow S^?\Upsilon$. For $x_i$ a variable of $\Gamma$, and $\omega \in \mathbb{O}$ such that $[\![\omega]\!]^? = x_i^\natural$, there is a corresponding cell $x_i \in S^?\Gamma_\omega$, and we let $S^?\sigma(x_i) := \sigma_i$. This is well-defined since by lemma lemma 4.7.1 on the facing page, we have $\sigma_i \in S^?\Upsilon_\omega$.

**Lemma 4.7.2.** *The mapping $S^?\sigma$ is a morphism of opetopic sets (i.e. a natural transformation) $S^?\Gamma \longrightarrow S^?\Upsilon$.*

*Proof.* Assume $\omega \neq \blacklozenge$, so that the type of $x_i$ is $\mathbf{P} \longrightarrow x_j$ for some $j < i$, and the type of $\sigma_i$ is $\mathbf{P}[\sigma_1/x_1]\cdots[\sigma_{i-1}/x_{i-1}] \longrightarrow \sigma_j$. Then $S^?\sigma(\mathsf{t}\,x_i) = \sigma_j = \mathsf{t}(S^?\sigma(x_i))$. If $[p] \in \omega^\bullet$, then $\mathsf{s}_{[p]}\,x_i = x_l$, for some $l < i$, and

$$S^?\sigma(\mathsf{s}_{[p]}\,x_i) = \sigma_l = \mathsf{s}_{[p]}\,(\mathbf{P}[\sigma_1/x_1]\cdots[\sigma_{i-1}/x_{i-1}]) = \mathsf{s}_{[p]}(S^?\sigma(x_i)).$$

Hence, $S^?\sigma$ is a morphism of opetopic sets $S^?\Gamma \longrightarrow S^?\Upsilon$. □

**Theorem 4.7.3.** *The functor $S^? : (\mathcal{C}\mathrm{tx}^?)^{\mathrm{op}} \longrightarrow \mathcal{F}\mathrm{in}\widehat{\mathbb{O}}$ is an equivalence of categories.*

*Proof.* Let $\Gamma, \Upsilon \in \mathcal{C}\mathrm{tx}^?$, with $\Gamma = (x_1 : T_1, \ldots, x_k : T_k)$, and $f : S^?\Gamma \longrightarrow S^?\Upsilon$. Then $f = S^?\,(f(x_1), \ldots, f(x_k))$, showing that $S^?$ is fully faithful. We now prove that it is essentially surjective. Take $X \in \mathcal{F}\mathrm{in}\widehat{\mathbb{O}}$, and enumerate its cells as $x_1 \in X_{\omega_1}, \ldots, x_k \in X_{\omega_k}$, such that whenever $i < j$ we have $\dim \omega_i \leq \dim \omega_j$. To each $0 \leq i \leq k$ we associate a derivable context $\Gamma^{(i)} = (\overline{x_1} : T_1, \ldots, \overline{x_i} : T_i)$ such that $\overline{x_i}^\natural = [\![\omega_i]\!]^?$, as follows:

(1) If $\omega_i = \blacklozenge$, let $\Gamma^{(i)}$ be given by the following proof tree (so that implicitly, $T_i = \blacklozenge$):

$$\dfrac{\Gamma^{(i-1)}\ \ \mathsf{context}}{\Gamma^{(i-1)}, \overline{x_i} : \blacklozenge\ \ \mathsf{context}}\ \texttt{point}$$

(2) Assume $\omega_i \neq \blacklozenge$ is not degenerate. Then, by induction, for $x_j = \mathsf{t}\,x_i$ we have $\overline{\mathsf{t}\,x_i}^\natural = \mathsf{t}\,\omega_i$, and for all address $[p_j] \in \omega_i^\bullet$, we have $\overline{\mathsf{s}_{[p_j]}\,x_i}^\natural = \mathsf{s}_{[p_j]}\,[\![\omega_i]\!]^?$. From this, $\Gamma^{(i)}$ is given by the following proof tree:

$$\dfrac{\dfrac{\Gamma^{(i-1)}\ \ \mathsf{context}}{\Gamma^{(i-1)} \vdash \begin{cases} [p_1] \leftarrow \overline{\mathsf{s}_{[p_1]}\,x_i} \\ \vdots \end{cases}}\ \texttt{graft}}{\Gamma^{(i-1)}, \overline{x_i} : \begin{cases} [p_1] \leftarrow \overline{\mathsf{s}_{[p_1]}\,x_i} \\ \vdots \end{cases} \longrightarrow \overline{\mathsf{t}\,x_i}\ \ \mathsf{context}}\ \texttt{shift}$$

(3) Assume $\omega_i$ is degenerate, Then $\Gamma^{(i)}$ is given by the following proof tree:

$$\dfrac{\dfrac{\Gamma^{(i-1)}\ \ \mathsf{context}}{\Gamma^{(i-1)} \vdash \{\!\{\overline{\mathsf{t}\,\mathsf{t}\,x_i}}\ \texttt{degen}}{\Gamma^{(i-1)}, \overline{x_i} : \{\!\{\overline{\mathsf{t}\,\mathsf{t}\,x_i} \longrightarrow \overline{\mathsf{t}\,x_i}\ \ \mathsf{context}}\ \texttt{shift}$$

Finally, the mapping $x_i \longmapsto \overline{x_i}$ exhibits an isomorphism $X \longrightarrow S^?\Gamma^{(k)}$. □

The category $\mathcal{C}\mathrm{tx}^?$ has finite limits, induced from finite colimits in $\mathcal{F}\mathrm{in}\widehat{\mathbb{O}}$ through $S^?$. We conclude this section with a result similar to theorem 3.6.17 on page 36, stating that opetopic sets essentially are "models of the algebraic theory $\mathcal{C}\mathrm{tx}^?$".

**Theorem 4.7.4.** *We have an equivalence* $\widehat{\mathbb{O}} \simeq \mathcal{M}\mathrm{od}(\mathcal{C}\mathrm{tx}^?, \mathcal{S}\mathrm{et})$, *where* $\mathcal{M}\mathrm{od}$ *is defined in theorem 3.6.17 on page 36*

*Proof.* We proceed as in the proof of theorem 3.6.17 on page 36. □

## 4.8. Examples

In this section, we give example derivations in system $\textsc{OptSet}^?$. For clarity, we do not repeat the type of previously typed variables in proof trees.

**Example 4.8.1.** We now derive the following opetopic set, which is not representable:

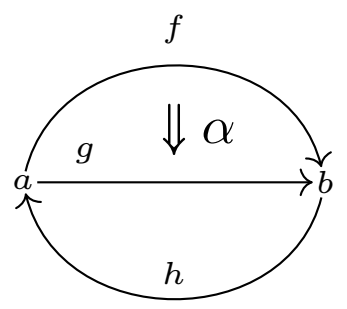

We decide to introduce all points first:

$$\dfrac{\dfrac{}{a : \blacklozenge \ \ \mathsf{context}}\ \texttt{point}}{a, b : \blacklozenge \ \ \mathsf{context}}\ \texttt{point}$$

Then we introduce $f$, by first specifying its source pasting diagram with the `graft` rule, parameterized by opetope $\blacksquare = \big\{ * \leftarrow \blacklozenge$, and then applying the `shift` rule:

$$\dfrac{\dfrac{\vdots \atop a, b \ \ \mathsf{context}}{a, b \vdash \big\{ * \leftarrow a}\ \texttt{graft}}{a, b, f : \big\{ * \leftarrow a \longrightarrow b \ \ \mathsf{context}}\ \texttt{shift}$$

We proceed similarly for $g$ and $h$:

$$\dfrac{\dfrac{\dfrac{\dfrac{\vdots \atop a, b, f \ \ \mathsf{context}}{a, b, f \vdash \big\{ * \leftarrow a}\ \texttt{graft}}{a, b, f, g : \big\{ * \leftarrow a \longrightarrow b \ \ \mathsf{context}}\ \texttt{shift}}{a, b, f, g \vdash \big\{ * \leftarrow b}\ \texttt{graft}}{a, b, f, g, h : \big\{ * \leftarrow b \longrightarrow a \ \ \mathsf{context}}\ \texttt{shift}$$

Lastly, we introduce $\alpha$, first by specifying its source with the `graft` rule, parameterized by opetope $\mathbf{1} = \big\{ [] \leftarrow \blacksquare$ (see the opetopic integers defined in example 4.3.2 on page 55), and applying the `shift` rule:

$$\dfrac{\dfrac{\vdots \atop a, b, f, g, h \ \ \mathsf{context}}{a, b, f, g, h \vdash \big\{ [] \leftarrow f}\ \texttt{graft}}{a, b, f, g, h, \alpha : \big\{ [] \leftarrow f \longrightarrow g \ \ \mathsf{context}}\ \texttt{shift}$$

**Example 4.8.2** (A classic, maximally folded)**.** We derive the following opetopic set

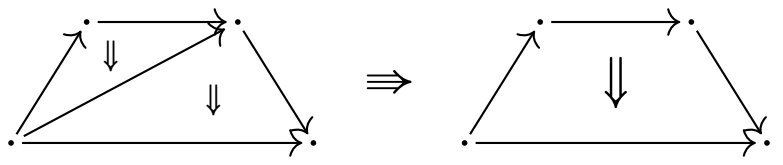

where all 0-cells are the same cell $a$, all 1-cells are $f$, the 2-cells on the left are $\alpha$, the 2-cell on the right is $\beta$, and the 3-cell is $A$. Note that those identifications make the opetopic set not representable. We first derive $a$ and $f$:

$$\dfrac{\dfrac{\dfrac{}{a : \blacklozenge \ \ \mathsf{context}}\ \texttt{point}}{a \vdash \big\{ * \leftarrow a}\ \texttt{graft}}{a, f : \big\{ * \leftarrow a \longrightarrow a \ \ \mathsf{context}}\ \texttt{shift}$$

Then we introduce $\alpha$, by first specifying its source pasting diagram with the `graft` rule, parameterized by opetope $\mathbf{2} = \begin{cases}[] \leftarrow \blacksquare \\ [*] \leftarrow \blacksquare\end{cases}$ (see the opetopic integers defined in example 4.3.2 on page 55), and applying the `shift` rule:

$$\dfrac{\dfrac{\begin{array}{c}\vdots\\ a, f \ \ \mathsf{context}\end{array}}{a, f \vdash \begin{cases}[] \leftarrow f\\ [*] \leftarrow f\end{cases}}\ \texttt{graft}}{a, f, \alpha : \begin{cases}[] \leftarrow f\\ [*] \leftarrow f\end{cases} \longrightarrow f \ \ \mathsf{context}}\ \texttt{shift}$$

Likewise, we introduce $\beta$, where the `graft` rule is parameterized by $\mathbf{3}$:

$$\dfrac{\dfrac{\begin{array}{c}\vdots\\ a, f, \alpha \ \ \mathsf{context}\end{array}}{a, f, \alpha \vdash \begin{cases}[] \leftarrow f\\ [*] \leftarrow f\\ [**] \leftarrow f\end{cases}}\ \texttt{graft}}{a, f, \alpha, \beta : \begin{cases}[] \leftarrow f\\ [*] \leftarrow f\\ [**] \leftarrow f\end{cases} \longrightarrow f \ \ \mathsf{context}}\ \texttt{shift}$$

Lastly, we introduce $A$, where the `graft` rule is parameterized by $\begin{cases}[] \leftarrow \mathbf{2}\\ [[*]] \leftarrow \mathbf{2}\end{cases}$:

$$\dfrac{\dfrac{\begin{array}{c}\vdots\\ a, f, \alpha, \beta \ \ \mathsf{context}\end{array}}{a, f, \alpha, \beta \vdash \begin{cases}[] \leftarrow \alpha\\ [[*]] \leftarrow \alpha\end{cases}}\ \texttt{graft}}{a, f, \alpha, \beta, A : \begin{cases}[] \leftarrow \alpha\\ [[*]] \leftarrow \alpha\end{cases} \longrightarrow \beta \ \ \mathsf{context}}\ \texttt{shift}$$

**Example 4.8.3.** The opetopic set

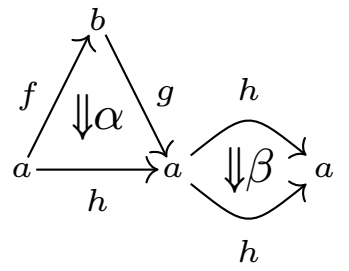

is straightforwardly derived as

$$\dfrac{\dfrac{\dfrac{\dfrac{\dfrac{\dfrac{}{a : \blacklozenge \ \ \mathsf{context}}\ \texttt{point}}{a, b : \blacklozenge \ \ \mathsf{context}}\ \texttt{point}}{a, b \vdash \{* \leftarrow a}\ \texttt{graft}}{a, b, f : \{* \leftarrow a \longrightarrow b \ \ \mathsf{context}}\ \texttt{shift}}{a, b, f \vdash \{* \leftarrow b}\ \texttt{graft}}{a, b, f, g : \{* \leftarrow b \longrightarrow a \ \ \mathsf{context}}\ \texttt{shift}$$

$$\dfrac{\dfrac{\dfrac{\dfrac{\dfrac{\dfrac{\vdots}{a, b, f, g \vdash \{* \leftarrow a}\ \texttt{graft}}{a, b, f, g, h : \{* \leftarrow a \longrightarrow a}\ \texttt{shift}}{a, b, f, g, h \vdash \begin{cases}[] \leftarrow g\\ [*] \leftarrow f\end{cases}}\ \texttt{graft}}{a, b, f, g, h, \alpha : \begin{cases}[] \leftarrow g\\ [*] \leftarrow f\end{cases} \longrightarrow h \ \ \mathsf{context}}\ \texttt{shift}}{a, b, f, g, h, \alpha \vdash \{[] \leftarrow h}\ \texttt{graft}}{a, b, f, g, h, \alpha, \beta : \{[] \leftarrow h \longrightarrow h}\ \texttt{shift}$$

## 4.9. Python implementation

System $\mathrm{OptSet}^?$ is implemented in module `opetopy.UnnamedOpetopicSet` of [**Ho Thanh, 2018b**]. The rules are represented by functions `point`, `degen`, `graft`, and `shift`, and are further encapsulated in rule instance classes `Point`, `Degen`, `Graft`, and `shift`.

FIGURE 4.9.1. Derivation of example 4.8.2 on page 62 using `opetopy.NamedOpetope`

```python
from opetopy.UnnamedOpetopicSet import *
from opetopy.UnnamedOpetope import address, Arrow, OpetopicInteger
from opetopy.UnnamedOpetope import Graft as OptGraft
from opetopy.UnnamedOpetope import Shift as OptShift
# We first derive the unique point a.
classic = Point("a")
# We then derive f by firstly specifying a pasting diagram with the graft rule. It is constructed with a proof tree of its shape opetope (in system OPT^?), and an address-to-variable mapping.
classic = Graft(pastingDiagram(
    Arrow(),
    {
        address([], 0): "a" # * ← a
    }), classic)
# We then derive f.
classic = Fill("a", "f", classic)
# In a similar way, we derive α of shape 2.
classic = Graft(pastingDiagram(
    OpetopicInteger(2),
    {
        address([], 1): "f", # [] ← f
        address(['*']): "f" # [*] ← f
    }), classic)
classic = Fill("f", "alpha", classic)
# In a similar way, we derive β of shape 3.
classic = Graft(pastingDiagram(
    OpetopicInteger(3),
    {
        address([], 1): "f", # [] ← f
        address(['*']): "f", # [*] ← f
        address(['*', '*']): "f" # [**] ← f
    }), classic)
classic = Fill("f", "beta", classic)
# We now take a break to derive ω = Y_2 ∘_[[*]] Y_2, the shape of A, in system OPT^?.
omega = OptGraft(OptShift(OpetopicInteger(2)),
                 OpetopicInteger(2),
                 address([['*']]))
# Finally, we derive A of shape ω.
classic = Graft(pastingDiagram(
    omega,
    {
        address([], 2): "alpha", # [] ← α
        address([['*']]): "alpha" # [[*]] ← f
    }), classic)
classic = Fill("beta", "A", classic)
```

CHAPTER 5

# Conclusion

We have introduced two syntactic presentations for opetopes, as well as for opetopic sets, and formally related them to preexisting definitions. Together with an adequate formulation of opetopic higher categories, it is our hope that this work will be used productively for mechanical proofs of coherence in opetopic $\omega$-categories or opetopic $\omega$-groupoids.

APPENDIX A

# Polynomial monads and the Baez–Dolan $(-)^+$ construction

From a polynomial monad (definition A.1.1, corollary A.1.3, and theorem A.1.9 on pages 71–72) $M$ with set of operations $B$, the Baez–Dolan construction gives a new polynomial monad $M^+$ having $B$ as its set of colors. In this chapter, we study the monad structure of $M^+$ in depth.

## A.1. Polynomial monads

A polynomial monad is classically defined as a monoid object in $\mathcal{P}\mathrm{oly}\mathcal{E}\mathrm{nd}(I)$, for some set $I$. However, there are two other approaches that we present in this section. The first one presents polynomial monads as algebras over the mouthful "free polynomial monad monad" [**Gambino and Kock, 2013, Kock, 2011**], and the other using so-called "partial laws" [**Kock et al., 2010**].

**Definition A.1.1** (Polynomial monad, classical definition)**.** A polynomial monad $M$ over a set $I$ is a monoid object in the category $\mathcal{P}\mathrm{oly}\mathcal{E}\mathrm{nd}(I)$ of polynomial endofunctors over $I$, where the monoidal product $\odot$ is the composition of polynomial functors [**Gambino and Kock, 2013**]. The structure morphisms $\eta : \mathrm{id} \longrightarrow M$ and $\mu : M \odot M \longrightarrow M$ indeed give the functor $M : \mathcal{S}\mathrm{et}/I \longrightarrow \mathcal{S}\mathrm{et}/I$ the structure of a cartesian monad.

Let $\mathcal{P}\mathrm{oly}\mathcal{M}\mathrm{nd}(I)$ be the subcategory of $\mathcal{P}\mathrm{oly}\mathcal{E}\mathrm{nd}(I)$ spanned by polynomial monads and morphisms of monads. We emphasize that in $\mathcal{P}\mathrm{oly}\mathcal{M}\mathrm{nd}(I)$, morphisms are expected to be identities on the set $I$ of colors. There is an obvious forgetful functor $U : \mathcal{P}\mathrm{oly}\mathcal{M}\mathrm{nd}(I) \longrightarrow \mathcal{P}\mathrm{oly}\mathcal{E}\mathrm{nd}(I)$ that turns out to have a left adjoint that we now present. Let $F \in \mathcal{P}\mathrm{oly}\mathcal{E}\mathrm{nd}(I)$ be as on the left, and define a new polynomial functor $F^\star$ as on the right,

$$I \xleftarrow{s} E \xrightarrow{p} B \xrightarrow{t} I, \qquad I \xleftarrow{s} \mathrm{tr}^{|} F \xrightarrow{u} \mathrm{tr}\, F \xrightarrow{r} I,$$

where $\mathrm{tr}\, F$ is the set of $F$-tree (see section 2.1 on page 9), and $\mathrm{tr}^{|} F$ is the set of $F$-trees equipped with a marked leaf.

**Theorem A.1.2** (Free monad [**Kock et al., 2010, Kock, 2011**])**.** *The polynomial functor $F^\star$ is canonically a polynomial monad. Moreover, the adjunction $(-)^\star \dashv U$ is monadic.*

Let us abuse notations and write $(-)^\star$ for the associated monad on the category $\mathcal{P}\mathrm{oly}\mathcal{E}\mathrm{nd}(I)$. Write its laws by $\mu^\star : (-)^{\star\star} \longrightarrow (-)^\star$ for multiplication (or *unbiased grafting*), and $\mathsf{Y} : \mathrm{id} \longrightarrow (-)^\star$ for the unit (or *introduction of corollas*).

**Corollary A.1.3** (Polynomial monad, algebraic definition)**.** *A polynomial monad is a $(-)^\star$-algebra, i.e. the data of a polynomial endofunctor $M$ together with a morphism $\sup : M^\star \longrightarrow M$ as*

$$\begin{array}{ccccccc}
I & \longleftarrow & \mathrm{tr}^{|} M & \longrightarrow & \mathrm{tr}\, M & \longrightarrow & I \\
\| & & \big\downarrow \wp & \lrcorner & \big\downarrow \mathsf{t} & & \| \\
I & \longleftarrow & E & \longrightarrow & B & \longrightarrow & I.
\end{array} \tag{A.1.4}$$

*such that the following two diagrams commute:*

$$\begin{array}{ccc}
M^{\star\star} & \xrightarrow{\sup^\star} & M^\star \\
{\scriptstyle \mu^\star}\big\downarrow & & \big\downarrow{\scriptstyle \sup} \\
M^\star & \xrightarrow{\sup} & M,
\end{array}
\qquad
\begin{array}{ccc}
M & \xrightarrow{\mathsf{Y}} & M^\star \\
& \searrow\!\!\!\!= & \big\downarrow{\scriptstyle \sup} \\
& & M.
\end{array} \tag{A.1.5}$$

The classical monad structure maps $\eta : \mathrm{id} \longrightarrow M$ and $\mu : M \odot M \longrightarrow M$ can be retrieved by noting that both id and $M \odot M$ are canonically subobjects of $M^\star$.

**Corollary A.1.6** (Contraction associativity formula)**.** *Let $T, U \in \operatorname{tr} M$, and $[l] \in T^|$ such that the grafting $T \circ_{[l]} U$ is defined. Then*

$$\mathsf{t}(T \underset{[l]}{\circ} U) = \mathsf{t}(\mathsf{Y}_{\mathsf{t}T} \underset{[\wp_T[l]]}{\circ} \mathsf{Y}_{\mathsf{t}U}). \tag{A.1.7}$$

*Further, for $[r] \in U^|$, we have a leaf $[l] \cdot [r] = [lr] \in (T \circ_{[l]} U)^|$, and*

$$\wp_{T \circ_{[l]} U}[lr] = \wp_V (\wp_T[l] \cdot \wp_U[r]) \quad \text{with} \quad V = \mathsf{Y}_{\mathsf{t}T} \underset{[\wp_T[l]]}{\circ} \mathsf{Y}_{\mathsf{t}U}. \tag{A.1.8}$$

*Proof.* Since $M$ is a $(-)^\star$-algebra, the left diagram of (A.1.5) commutes. □

**Theorem A.1.9** (Polynomial monads via partial laws)**.** *Take a polynomial endofunctor $M \in \mathcal{P}\mathrm{oly}\mathcal{E}\mathrm{nd}(I)$, and write it as:*

$$I \xleftarrow{s} E \xrightarrow{p} B \xrightarrow{t} I.$$

*A $(-)^\star$-algebra structure on $M$ (i.e. a structure of polynomial monad on $M$) is equivalent to the following data:*

1. **(Unit)** *a map $\eta : I \longrightarrow B$;*
2. **(Partial multiplication)** *a map $\square : E \times_I B \longrightarrow B$, where for $(e, b) \in E \times_I B$ and $a = p(e)$ we write $a \square_e b$ for $\square(e, b)$, and say that $a \square_e b$ is an* admissible expression *(or just* admissible*); in the sequel, we assume that all expressions of this kind are admissible;*
3. **(Partial reindexing)** *for $a \square_e b$ admissible, an isomorphism*

$$\chi_{a \square_e b} : E(a) + E(b) - \{e\} \xrightarrow{\cong} E(a \underset{e}{\square} b);$$

*such that:*

1. **(Trivial)** *for $i \in I$, we have $t(\eta(i)) = i$, and $E(\eta(i))$ is a singleton whose unique element $e$ is such that $s(e) = i$;*
2. **(Left unit)** *for $i \in I$, $b \in B$, and $e$ the unique element of $E(\eta(i))$, we have $\eta(i) \square_e b = b$, and $\chi_{\eta(i) \square_e b} : E(b) \longrightarrow E(b)$ is the identity;*
3. **(Right unit)** *for $i \in I$, $b \in B$, and $e \in E(b)$, we have $b \square_e \eta(i) = b$, and $\chi_{b \square_e \eta(i)}$ is given by*

$$\begin{aligned} E(b) + E(\eta(i)) - \{e\} &\longrightarrow E(b) \\ x \in E(b) &\longmapsto x \\ y \in E(\eta(i)) &\longmapsto e; \end{aligned}$$

4. **(Disjoint multiplication)** *for $a \square_e b$ and $a \square_f c$ admissible expressions, and $e \neq f$, we have*

$$(a \underset{e}{\square} b) \underset{\chi_{a \square_e b}(f)}{\square} c \quad = \quad (a \underset{f}{\square} c) \underset{\chi_{a \square_f c}(e)}{\square} b,$$

*and the following diagram commutes:*

$$\begin{array}{ccc} E(a) + E(b) + E(c) - \{e, f\} & \xrightarrow{\chi_{a \square_e b}} & E(a \underset{e}{\square} b) + E(c) - \{f\} \\ {\scriptstyle \chi_{a \square_f c}} \downarrow & & \downarrow {\scriptstyle \chi_C} \\ E(a \underset{f}{\square} c) + E(b) - \{e\} & \xrightarrow{\chi_B} & E((a \underset{e}{\square} b) \underset{\chi_{a \square_e b} f}{\square} c); \end{array} \tag{A.1.10}$$

*with*

$$B = (a \underset{e}{\square} b) \underset{\chi_{a \square_e b}(f)}{\square} c, \qquad C = (a \underset{f}{\square} c) \underset{\chi_{a \square_f c} e}{\square} b;$$

5. **(Nested multiplication)** *for $a \square_e b$ and $b \square_f c$ admissible, we have*

$$(a \underset{e}{\square} b) \underset{\chi_{a \square_e b} f}{\square} c \quad = \quad a \underset{e}{\square} (b \underset{f}{\square} c),$$

*and the following diagram commutes:*

$$\begin{array}{ccc} E(a)+E(b)+E(c)-\{e,f\} & \xrightarrow{\chi_{a \,\square_e\, b}} & E(a \underset{e}{\square} b)+E(c)-\{f\} \\ \chi_{b \,\square_f\, c}\big\downarrow & & \big\downarrow \chi_{a\,\square_e (b \,\square_f\, c)} \\ E(a \underset{f}{\square} c)+E(b)-\{e\} & \xrightarrow{\chi_B} & E((a \underset{e}{\square} b) \underset{\chi_{a \,\square_e\, b} f}{\square} c). \end{array} \tag{A.1.11}$$

*with*

$$B = (a \underset{e}{\square} b) \underset{\chi_{a \,\square_e\, b} f}{\square} c.$$

*Moreover with the data above, the components of the structure map $M^\star \longrightarrow M$ (as in diagram (A.1.4)) are inductively given by:*

*(1) (Target) for $i \in I$, $\mathsf{t}\,\mathsf{I}_i = \eta(i)$; for $b \in B$, $\mathsf{t}\,\mathsf{Y}_b = b$; for $T \in \operatorname{tr} M^+$, $[l] \in T^|$, and $b \in B$ such that the grafting $T \circ_{[l]} \mathsf{Y}_b$ is defined:*

$$\mathsf{t}(T \underset{[l]}{\circ} \mathsf{Y}_b) = (\mathsf{t}\,T) \underset{\wp_T [l]}{\square} b,$$

*where $\wp_T$ is defined next;*

*(2) (Readdressing) for $b \in B$,*

$$\begin{aligned} \wp_{\mathsf{Y}_b} : \mathsf{Y}_b^| &\longrightarrow E(b) \\ [e] &\longmapsto e; \end{aligned}$$

*for $T \in \operatorname{tr} M^+$, $[l] \in T^|$, and $b \in B$ such that the grafting $T \circ_{[l]} \mathsf{Y}_b$ is defined, the readdressing $\wp_{T \circ_{[l]} \mathsf{Y}_b}$ is given by*

$$\begin{aligned} (T \underset{[l]}{\circ} \mathsf{Y}_b)^| = T^| + \mathsf{Y}_b^| - \{[l]\} &\longrightarrow E((\mathsf{t}\,T) \underset{\wp_T [l]}{\square} b) \\ [p] \in T^| - \{[l]\} &\longmapsto \chi_{(\mathsf{t}\,T) \,\square_{\wp_T [l]}\, b}(\wp_T [l]) \\ [e] \in \mathsf{Y}_b^| &\longmapsto \chi_{(\mathsf{t}\,T) \,\square_{\wp_T [l]}\, b}(e) \end{aligned}$$

*Proof (sketch).* Conditions **(Disjoint multiplication)** and **(Nested multiplication)** ensure that $\mathsf{t}\,T$ does not depend on the chosen decomposition, and thus define a structure map $\sup : M^\star \longrightarrow M$ making the left square of (A.1.5) commute. Conditions **(Trivial)**, **(Left unit)**, and **(Right unit)** ensure that the right triangle of (A.1.5) commutes too, and thus that $M$ is a $(-)^\star$-algebra, that is, a polynomial monad. ☐

## A.2. A complete $(-)^+$ construction

Let $M$ be a polynomial monad and define $M^+$ as being the following polynomial functor:

$$B \xleftarrow{\ \mathsf{s}\ } \operatorname{tr}^\bullet M \xrightarrow{\ u\ } \operatorname{tr} M \xrightarrow{\ \mathsf{t}\ } B, \tag{A.2.1}$$

where for $T \in \operatorname{tr} M$, the fiber $u^{-1}T$ (which we shall also denote by $T^\bullet$) is the set of node addresses in $T$, and for $[p] \in T^\bullet$, $\mathsf{s}[p] := \mathsf{s}_{[p]} T$.

**Definition A.2.2** (Monad structure on $M^+$)**.** We endow $M^+$ with a monad structure by means of theorem A.1.9 on the preceding page. The maps $\eta^+ : I \longrightarrow B$ and $\square : \operatorname{tr}^\bullet M \times_B \operatorname{tr} M$ are defined as follows:

(1) **(Unit)** for $b \in B$, let $\eta^+(b) := \mathsf{Y}_b$;

(2) **(Partial multiplication)** for $U, T \in \operatorname{tr} M$, $[p] \in U^\bullet$ such that $\mathsf{s}_{[p]} U = \mathsf{t}\,T$, write (recall that $\circ$ is associative on the right)

$$U \quad = \quad X \underset{[p]}{\circ} \mathsf{Y}_{\mathsf{t}\,T} \bigcirc_{[e_i]} Y_i$$

where $\{[e_i]\}_i \subseteq \mathsf{Y}_{\mathsf{t}\,T}^|$ (equivalently, $\{e_i\}_i \subseteq E(\mathsf{t}\,T)$), and define

$$U \underset{[p]}{\square} T \quad := \quad X \underset{[p]}{\circ} T \bigcirc_{\chi_T^{-1} e_i} Y_i;$$

(3) (Partial reindexing) for $U \,\square_{[p]}\, T$ admissible, define $\chi_{U \square_{[p]} T}$ by

$$U^\bullet + T^\bullet - \{[p]\} \longrightarrow (U \underset{[p]}{\square} T)^\bullet$$

$$[q] \in T^\bullet \longmapsto [pq] \tag{A.2.3}$$

$$[p[e_i]p'] \in U^\bullet \longmapsto [p] \cdot \chi_T^{-1} e_i \cdot [p'] \tag{A.2.4}$$

$$[q] \in U^\bullet \text{ not as above} \longmapsto [q] \tag{A.2.5}$$

With this definition, the $(-)^\star$-algebra structure map $m^+ : (M^+)^\star \longrightarrow M^+$ of $M^+$ (see (A.1.5)) is given as follows. The target map $\mathsf{t} : \operatorname{tr} M^+ \longrightarrow \operatorname{tr} M$ is defined by induction: for $b \in B$, $\mathsf{t}\,\mathsf{I}_b = \mathsf{Y}_b$, for $T \in \operatorname{tr} M$, $\mathsf{t}\,\mathsf{Y}_T = T$, and for $U, v \in \operatorname{tr} M^+$, $[l] \in U^|$ such that $U \circ_{[l]} V$ is defined,

$$\mathsf{t}(U \underset{[l]}{\circ} V) = (\mathsf{t}\,U) \underset{\wp_U[l]}{\square} (\mathsf{t}\,V). \tag{A.2.6}$$

The readdressing map $\wp : \operatorname{tr}^| M^+ \longrightarrow \operatorname{tr}^\bullet M$ also admit a somewhat simple description. Let $U \in \operatorname{tr} M^+$, and $[p[q]] \in U^|$. Then $U$ decomposes as $U = V \circ_{[p]} W$, for some $V, W \in \operatorname{tr} M^+$, and $[q] \in (\mathsf{s}_{[]}\, W)^\bullet$. Then,

$$\wp_U[p[q]] = (\wp_V[p]) \cdot (\wp_W[[q]]). \tag{A.2.7}$$

We unfold this equality further in the special case of "level 2 trees", i.e. those trees $V \in \operatorname{tr} M^+$ whose leaf addresses are of the form $[[a][b]]$, for $[a] \in (\mathsf{s}_{[]}\, V)^\bullet$, and $[b] \in (\mathsf{s}_{[a]}\, V)^\bullet$. Trees of this form are the operations of $M^+ \odot M^+$. We can decompose $V$ as

$$V \quad = \quad T \underset{[p_i]}{\bigcirc} U_i,$$

where $[p_i]$ ranges over $T^\bullet$. Write $\{[q_{i,1}], \ldots, [q_{i,k_i}]\} = U_i^\bullet$, and let $[[p_i][q_{i,j}]] \in V^|$. Then $[p_i]$ is a node address of $T$, so we may decompose it as $[p_i] = [[r_1] \cdots [r_k]]$, and

$$\wp_V[[p_i][q_{i,j}]] = (\wp^{-1}_{\mathsf{s}_{[]} V}[r_1]) \cdot (\wp^{-1}_{\mathsf{s}_{[[r_1]]} V}[r_2]) \cdot \cdots \cdot (\wp^{-1}_{\mathsf{s}_{[[r_1]\cdots[r_{l-1}]]} V}[r_l]) \cdot [q_{i,j}]. \tag{A.2.8}$$

The rest of this section is dedicated to prove that A.2.2 indeed defines a monad. We thus check the conditions of theorem A.1.9 on page 72.

*Proof of* **(Trivial)**. For $b \in B$ we have $\mathsf{t}\,\eta^+(b) = \mathsf{t}\,\mathsf{Y}_b = b$. On the other hand, $(\eta^+(b))^\bullet = \{[]\}$ and $\mathsf{s}_{[]}\,\eta^+(b) = \mathsf{s}_{[]}\,\mathsf{Y}_b = b$. ☐

*Proof of* **(Left unit)**. Let $b \in B$ and $T \in \operatorname{tr} M$ be such that $\eta^+(b) \,\square_{[]}\, T$ is admissible. Then, by definition, $\eta^+(b) \,\square_{[]}\, T = \mathsf{Y}_b \,\square_{[]}\, T = T$. Then, $\eta^+(b)^\bullet + T^\bullet - \{[]\} = T^\bullet$, and $\chi_{\eta^+(b) \square_{[]} T}$ maps $[q] \in T^\bullet$ to $[q]$, so it is indeed the identity. ☐

*Proof of* **(Right unit)**. Let $b \in B$, $T \in \operatorname{tr} M$, and $[p] \in T^\bullet$ be such that $T \,\square_{[p]}\, \eta^+(b)$ is admissible. By definition, $T \,\square_{[p]}\, \eta^+(b) = T \,\square_{[p]}\, \mathsf{Y}_b = T$. Then, $\chi_{T \square_{[p]} \eta^+(b)}$ is given by

$$T^\bullet + \{[]\} - \{[p]\} \longrightarrow T^\bullet$$

$$[p[e]p'] \in T^\bullet \longmapsto [p] \cdot \chi^{-1}_{\eta^+(b)} e \cdot [p'] = [p[e]p']$$

$$[p'] \in T^\bullet \text{ not as above} \longmapsto [p']$$

$$[] \in \eta^+(b)^\bullet \longmapsto [p]$$

as indeed $\chi_{\eta^+(b)}$ maps $[e] \in (\mathsf{Y}_b)^|$ to $e \in E(b)$. ☐

*Proof of* **(Disjoint multiplication)**. Let $A, B, C \in \operatorname{tr} M$, $[e], [f] \in A^\bullet$ be different, and such that $A \,\square_{[e]}\, B$ and $A \,\square_{[f]}\, C$ are admissible. Without loss of generality, we distinguish two cases: one where $[e] \sqsubseteq [f]$, for $\sqsubseteq$ the prefix order, and one where $[e]$ and $[f]$ are $\sqsubseteq$-incomparable.

(1) Assume $[e] \sqsubseteq [f]$, so that $[f] = [e[q]r]$ for some $e$ and $r$, and write $A$ as (recall that $\circ$ is associative on the right)

$$A \quad = \quad X \underset{[e]}{\circ} \mathsf{Y}_{\mathsf{t}B} \underset{[q]}{\circ} Y \underset{[r]}{\circ} \mathsf{Y}_{\mathsf{t}C} \underset{[v_i]}{\bigcirc} Z_i,$$

where $q \in E(\mathsf{t}B)$ and $\{v_i\}_i \subseteq E(\mathsf{t}C)$. Then

$$A \underset{[e]}{\square} B \;=\; X \underset{[e]}{\circ} B \underset{\chi_B^{-1}q}{\circ} Y \underset{[r]}{\circ} \mathsf{Y}_{\mathsf{t}C} \underset{[v_i]}{\bigcirc} Z_i,$$

and $\chi_{A \square_{[e]} B}[f] = \chi_{A \square_{[e]} B}[e[q]r] = [e] \cdot \chi_B^{-1} q \cdot [r]$. Thus,

$$(A \underset{[e]}{\square} B) \underset{[e]\cdot\chi_B^{-1}q\cdot[r]}{\square} C \;=\; X \underset{[e]}{\circ} B \underset{\chi_B^{-1}q}{\circ} Y \underset{[r]}{\circ} C \underset{\chi_C^{-1}v_i}{\bigcirc} Z_i,$$

and the reindexing $\chi_{(A \square_{[e]} B) \square_{\chi_{A \square_{[e]} B}[f]} C}\chi_{A \square_{[e]} B}$ is given by

$$\begin{array}{llllll}
[p] \in B^\bullet & \longmapsto & [ep] & \longmapsto & [ep] \\
 & & [p] \in C^\bullet & \longmapsto & [e] \cdot \chi_B^{-1} q \cdot [rp] \\
[e[q]s] \subsetneq [f] & \longmapsto & [e] \cdot \chi_B^{-1} q \cdot [s] & \longmapsto & [e] \cdot \chi_B^{-1} q \cdot [s] \\
[f[v_i]s] & \longmapsto & [e] \cdot \chi_B^{-1} q \cdot [r] \cdot [[v_i]s] & \longmapsto & [e] \cdot \chi_B^{-1} q \cdot [r] \cdot \chi_C^{-1} v_i \cdot [s] \\
[p] \in A^\bullet \text{ not as above} & \longmapsto & [p] & \longmapsto & [p]
\end{array}$$

On the other hand, we have

$$A \underset{[f]}{\square} C \;=\; X \underset{[e]}{\circ} \mathsf{Y}_{\mathsf{t}B} \underset{[q]}{\circ} Y \underset{[r]}{\circ} C \underset{\chi_C^{-1}v_i}{\bigcirc} Z_i.$$

The reindexing gives $\chi_{A \square_{[f]} C}[e] = [e]$, and

$$\begin{aligned}
(A \underset{[f]}{\square} C) \underset{[e]}{\square} B \;&=\; X \underset{[e]}{\circ} B \underset{\chi_B^{-1}q}{\circ} Y \underset{[r]}{\circ} C \underset{\chi_C^{-1}v_i}{\bigcirc} Z_i \\
&=\; (A \underset{[e]}{\square} B) \underset{[e]\cdot\chi_B^{-1}q\cdot[r]}{\square} C.
\end{aligned}$$

The reindexing $\chi_{(A \square_{[f]} C) \square_{[e]} B}\chi_{A \square_{[f]} C}$ is given by

$$\begin{array}{llllll}
 & & [p] \in B^\bullet & \longmapsto & [ep] \\
[p] \in C^\bullet & \longmapsto & [fp] & \longmapsto & [e] \cdot \chi_B^{-1} q \cdot [rp] \\
[e[q]s] \subsetneq [f] & \longmapsto & [e[q]s] & \longmapsto & [e] \cdot \chi_B^{-1} q \cdot [s] \\
[f[v_i]s] & \longmapsto & [f] \cdot \chi_C^{-1} v_i \cdot [s] & \longmapsto & [e] \cdot \chi_B^{-1} q \cdot [r] \cdot \chi_C^{-1} v_i \cdot [s] \\
[p] \in A^\bullet \text{ not as above} & \longmapsto & [p] & \longmapsto & [p]
\end{array}$$

We see that the square (A.1.10) commutes in the case $[e] \sqsubseteq [f]$.

(2) Assume $[e]$ and $[f]$ are $\sqsubseteq$-incomparable, and write $A$ as

$$A \;=\; (X \underset{[e]}{\circ} \mathsf{Y}_{\mathsf{t}B} \underset{[v_i]}{\bigcirc} Y_i) \underset{[f]}{\circ} \mathsf{Y}_{\mathsf{t}C} \underset{[w_j]}{\bigcirc} Z_j.$$

Then

$$A \underset{[e]}{\square} B \;=\; (X \underset{[e]}{\circ} B \underset{\chi_B^{-1}v_i}{\bigcirc} Y_i) \underset{[f]}{\circ} \mathsf{Y}_{\mathsf{t}C} \underset{[w_j]}{\bigcirc} Z_j,$$

the reindexing gives $\chi_{A \square_{[e]} B}[f] = [f]$,

$$(A \underset{[e]}{\square} B) \underset{[f]}{\square} C \;=\; (X \underset{[e]}{\circ} B \underset{\chi_B^{-1}v_i}{\bigcirc} Y_i) \underset{[f]}{\circ} C \underset{\chi_C^{-1}w_j}{\bigcirc} Z_j,$$

and the complete reindexing $\chi_{(A \square_{[e]} B) \square_{[f]} C}\chi_{A \square_{[e]} B}$ is given by

$$\begin{array}{llllll}
[p] \in B^\bullet & \longmapsto & [ep] & \longmapsto & [ep] \\
 & & [p] \in C^\bullet & \longmapsto & [fp] \\
[e[v_i]s] & \longmapsto & [e] \cdot \chi_B^{-1} v_i \cdot [s] & \longmapsto & [e] \cdot \chi_B^{-1} v_i \cdot [s] \\
[f[w_j]s] & \longmapsto & [f[w_j]s] & \longmapsto & [f] \cdot \chi_C^{-1} w_j \cdot [s] \\
[p] \in A^\bullet \text{ not as above} & \longmapsto & [p] & \longmapsto & [p]
\end{array}$$

On the other hand,

$$A \underset{[f]}{\square} C \;=\; (X \underset{[e]}{\circ} \mathsf{Y}_{\mathsf{t}B} \underset{[v_i]}{\bigcirc} Y_i) \underset{[f]}{\circ} C \underset{\chi_C^{-1}w_j}{\bigcirc} Z_j,$$

we have $\chi_{A \,\square_{[f]}\, C}[e] = [e]$,

$$\begin{aligned} (A \underset{[f]}{\square} C) \underset{[e]}{\square} B &= (X \underset{[e]}{\circ} B \bigcirc_{\chi_B^{-1} v_i} Y_i) \underset{[f]}{\circ} C \bigcirc_{\chi_C^{-1} w_j} Z_j \\ &= (A \underset{[e]}{\square} B) \underset{[f]}{\square} C, \end{aligned}$$

and further

$$\begin{array}{lclcl} & & [p] \in B^\bullet & \longmapsto & [ep] \\ [p] \in C^\bullet & \longmapsto & [fp] & \longmapsto & [fp] \\ e[v_i]s] & \longmapsto & [e[v_i]s] & \longmapsto & [e] \cdot \chi_B^{-1} v_i \cdot [s] \\ [f[w_j]s] & \longmapsto & [f] \cdot \chi_C^{-1} w_j \cdot [s] & \longmapsto & [f] \cdot \chi_C^{-1} w_j \cdot [s] \\ [p] \in A^\bullet \text{ not as above} & \longmapsto & [p] & \longmapsto & [p] \end{array}$$

so that the square (A.1.10) commutes in the case where $[e]$ and $[f]$ are $\sqsubseteq$-incomparable too. Finally, the monad structure of $M^+$ given in definition A.2.2 on page 73 satisfies condition **(Disjoint multiplication)** of theorem A.1.9 on page 72. □

*Proof of* **(Nested multiplication)**. Let $A, B, C \in \operatorname{tr} M$, $[e] \in A^\bullet$, $[f] \in B^\bullet$, such that $A \,\square_{[e]}\, B$ and $B \,\square_{[f]}\, C$ are admissible. Write

$$A = (X \underset{[e]}{\circ} \mathsf{Y}_{\mathsf{t} B} \bigcirc_{[v_i]} Y_i), \qquad \text{and} \qquad B = Z \underset{[f]}{\circ} \mathsf{Y}_{\mathsf{t} C} \bigcirc_{[w_j]} T_j.$$

Then,

$$A \underset{[e]}{\square} B = X \underset{[e]}{\circ} B \bigcirc_{\chi_B^{-1} v_i} Y_i,$$

we have $\chi_{A \,\square_{[e]}\, B}[f] = [ef]$, and

$$(A \underset{[e]}{\square} B) \underset{[ef]}{\square} C = X \underset{[e]}{\circ} (Z \underset{[f]}{\circ} \mathsf{Y}_{\mathsf{t} C} \bigcirc_{[w_j]} T_j) \bigcirc_{\alpha(\chi_B^{-1} v_i)} Y_i,$$

where

$$\alpha(\chi_B^{-1} v_i) = \begin{cases} [f] \cdot \chi_C^{-1} w_j \cdot [r] & \text{if } \chi_B^{-1} v_i \text{ of the form } [f[w_j]r], \\ \chi_B^{-1} v_i & \text{otherwise.} \end{cases}$$

Remark that $\alpha(\chi_B^{-1} v_i) = \chi^{-1}_{B \,\square_{[f]}\, C} v_i$. The reindexing $\chi_{(A \,\square_{[e]}\, B) \,\square_{[ef]}\, C} \chi_{A \,\square_{[e]}\, B}$ is given by:

$$\begin{array}{lclcl} & & [p] \in C^\bullet & \longmapsto & [efp] \\ [f[w_j]r] \in B^\bullet & \longmapsto & [ef[w_j]r] & \longmapsto & [ef] \cdot \chi_C^{-1} w_j \cdot [r] \\ [p] \in B^\bullet,\ [f] \not\sqsubseteq [p] & \longmapsto & [ep] & \longmapsto & [ep] \\ [e[v_i]r] \in A^\bullet & \longmapsto & [e] \cdot \chi_B^{-1} v_i \cdot [r] & \longmapsto & [e] \cdot \chi^{-1}_{B \,\square_{[f]}\, C} v_i \cdot [r] \end{array}$$

On the other hand, we have

$$\begin{aligned} B \underset{[f]}{\square} C &= Z \underset{[f]}{\circ} C \bigcirc_{\chi_C^{-1} w_j} T_j, \\ A \underset{[e]}{\square} (B \underset{[f]}{\square} C) &= X \underset{[e]}{\circ} (Z \underset{[f]}{\circ} \mathsf{Y}_{\mathsf{t} C} \bigcirc_{[w_j]} T_j) \bigcirc_{\chi^{-1}_{B \,\square_{[f]}\, C} v_i} Y_i \\ &= (A \underset{[e]}{\square} B) \underset{[ef]}{\square} C \end{aligned}$$

and the reindexing is given by

$$\begin{array}{lclcl} [p] \in C^\bullet & \longmapsto & [fp] & \longmapsto & [efp] \\ [f[w_j]r] \in B^\bullet & \longmapsto & [f] \cdot \chi_C^{-1} w_j \cdot [r] & \longmapsto & [ef] \cdot \chi_C^{-1} w_j \cdot [r] \\ [p] \in B^\bullet,\ [f] \not\sqsubseteq [p] & \longmapsto & [p] & \longmapsto & [ep] \\ & & [e[v_i]r] \in A^\bullet & \longmapsto & [e] \cdot \chi^{-1}_{B \,\square_{[f]}\, C} v_i \cdot [r] \end{array}$$

We thus see that the square (A.1.11) commutes, and that the monad structure of $M^+$ given in definition A.2.2 on page 73 satisfies condition **(Nested multiplication)** of theorem A.1.9 on page 72. □